\documentclass[english,11pt]{article}%
\usepackage{bbm,ulem}
\usepackage{graphicx}
\usepackage{subcaption}
\usepackage{makeidx}
\usepackage{amssymb}
\usepackage[utf8]{inputenc}
\usepackage{amsfonts}
\usepackage{amsmath}
\usepackage{amssymb}
\usepackage{amsthm}
\usepackage{url}
\usepackage[colorlinks,citecolor=blue]{hyperref}
\usepackage{mathabx}
\usepackage{wrapfig}
\usepackage{dsfont}
\usepackage{resizegather}
\usepackage{yfonts}
\usepackage{authblk}
\usepackage{geometry}
\usepackage[dvipsnames]{xcolor}

\DeclareMathOperator*{\argmin}{arg\,min}

\providecommand{\U}[1]{\protect\rule{.1in}{.1in}}

\usepackage{xcolor}

\newtheorem{prop}{Proposition}[section]
\newtheorem{cor}[prop]{Corollary}
\newtheorem{defi}[prop]{Definition}

\newtheorem{rmk}[prop]{Remark}
\newtheorem{lem}[prop]{Lemma}

\newtheorem{theo}[prop]{Theorem}
\newtheorem{examp}[prop]{Example}

\newcommand{\EE}{\mathbb{E}}

\newcommand{\LL}{\mathbb{L}}

\newcommand{\PP}{\mathbb{P}}

\newcommand{\RR}{\mathbb{R}}

\newcommand{\TT}{\mathbb{T}}
\newcommand{\UU}{\mathbb{U}}
\newcommand{\VV}{\mathbb{V}}

\newcommand{\WW}{\mathbb{W}}

\newcommand{\ZZ}{\mathbb{Z}}

\newcommand{\Aa}{ {\cal A }}

\newcommand{\Da}{ {\cal D }}
\newcommand{\La}{ {\cal L }}

\newcommand{\Ka}{ {\cal K }}

\newcommand{\Ea}{ {\cal E }}
\newcommand{\Sa}{ {\cal S }}

\newcommand{\Va}{ {\cal V }}

\newcommand{\Fa}{ {\cal F }}
\newcommand{\Ga}{ {\cal G }}

\newcommand{\Ma}{ {\cal M }}
\newcommand{\Ta}{ {\cal T}}
\newcommand{\Ha}{ {\cal H }}
\newcommand{\Ja}{ {\cal J }}
\newcommand{\Pa}{ {\cal P }}
\newcommand{\Za}{ {\cal Z }}

\newcommand{\Wa}{ {\cal W }}

\newcommand{\point}{\mbox{\LARGE .}}
\newcommand{\cqfd}{\hfill\blbx \\}
\def\blbx{\hbox{\vrule height 5pt width 5pt depth 0pt}\medskip}

\def \PP{\mathbb{P}}
\def \RR{\mathbb{R}}

\def \EE{\mathbb{E}}

\def \LL{\mathbb{L}}
\def \JJ{\mathbb{J}}
\def \ZZ{\mathbb{Z}}
\def \WW{\mathbb{W}}

\newcommand{\cchi}{\protect\raisebox{2pt}{$\chi$}}

\numberwithin{equation}{section}

\newcommand{\vertiii}[1]{{\left\vert\kern-0.25ex\left\vert\kern-0.25ex\left\vert #1
    \right\vert\kern-0.25ex\right\vert\kern-0.25ex\right\vert}}

\usepackage{amsopn}

\begin{document}

  \title{A Contraction Theory for Sinkhorn and Schrödinger Bridges via Log-Sobolev Inequalities}

\author{P. Del Moral\thanks{Centre de Recherche Inria Bordeaux Sud-Ouest, Talence, 33405, France. {\footnotesize E-Mail:\,} \texttt{\footnotesize pierre.del\_moral@inria.fr}}}

\maketitle
  \begin{abstract}

 We develop a quantitative contraction framework for Schrödinger and Sinkhorn bridges based on transportation–cost inequalities and Riccati matrix difference equations. Our approach combines logarithmic Sobolev and Talagrand-type inequalities to obtain explicit entropy and Wasserstein contraction bounds for Sinkhorn bridge measures, entropic optimal transport plans, and the associated Markov transport maps. A key feature of the analysis is the interplay between transport-cost inequalities and matrix Riccati difference equations arising in filtering and stochastic control.
 
The results are established under local regularity assumptions on the reference transition, formulated in terms of curvature, Lipschitz continuity, and Fisher-information bounds. Within this general setting, we derive quantitative stability and convergence estimates for Schrödinger bridges and Sinkhorn iterates that are robust with respect to the choice of reference measure.
  
As a main application, we specialize the theory to linear-Gaussian reference transitions, where the Gaussian structure permits sharp constants, refined exponential decay rates, and continuity estimates for Schrödinger bridges, Sinkhorn iterates, barycentric projections, conditional covariances, and proximal sampler semigroups. In this setting, we recover and extend several known contraction results for entropic and Wasserstein distances, and obtain new quantitative bounds that improve previously available rates. Our results provide a unified probabilistic framework for stability, regularity, and convergence of Sinkhorn algorithms. We illustrate the impact of our results on regularized entropic transport, proximal samplers, diffusion-based generative models, as well as on diffusion flow-matching.
\end{abstract}

\newpage
{\scriptsize
\tableofcontents}
\newpage

\section{Introduction}
\subsection{Some background}
Consider probability measures $\mu$ and $\eta$ on $\RR^d$.
 The Schr\" odinger bridge from $\mu$ to 
$\eta$ with respect to a reference probability distribution $ P$ on $(\RR^d\times\RR^d)$ is defined by the bridge entropic projection map
\begin{equation}\label{ob-eq}
P\mapsto P_{\mu,\eta}:= \argmin_{Q\,\in\, \Pi(\mu,\eta)}\Ha(Q~|~P).
 \end{equation}
In the above display, $\Ha(Q|P)= 
Q(\log{\left({dQ}/{dP}\right)})
$ is the relative entropy of $Q$ w.r.t. $P$ and  $\Pi(\mu,\eta)$ is the convex subset of probability measures $Q$ on $
(\RR^d\times\RR^d)$ with marginals $\mu$ (first coordinate) and $\eta$ (second coordinate), see Section~\ref{transport-sec} for a more precise description.

Existence of a solution is guaranteed provided there exists $Q\in \Pi(\mu,\eta)$ such that 
$\Ha(Q~|~P)<\infty$.
This classical fact goes back to Csisz\'ar~\cite{csiszar-2}, see also Section 6 of the Lecture Notes by Nutz~\cite{nutz}, and the survey article by L\'eonard~\cite{leonard}.

In contrast with classical optimal transport, which may be difficult to compute in high dimensions, Schr\" odinger bridges can be efficiently approximated using the celebrated Sinkhorn algorithm (also known as the iterative proportional fitting procedure)~\cite{cuturi,genevay-cuturi-2,genevay,sinkhorn,sinkhorn-2,sinkhorn-3}.   For clarity, we recall its basic structure; further details and proofs can be found in Appendix~\ref{review-sinkhorn}.
Consider a given reference measure of the form
$$
P(d(x,y))=(\mu\times \Ka)(d(x,y)):=\mu(dx)~\Ka(x,dy),
$$
for some Markov transition $\Ka(x,dy)$ on $\RR^d$.
Starting from the reference measure $\Pa_0=P$ Sinkhorn recursion alternates between two entropic projections
\begin{equation}\label{sinhorn-entropy-form}
\Pa_{2n+1}:= \argmin_{\Pa\in \Pi(\point,\eta)}\Ha(\Pa~|~\Pa_{2n})\quad \mbox{\rm and}\quad
\Pa_{2(n+1)}:= \argmin_{\Pa\in  \Pi(\mu,\point)}\Ha(\Pa~|~\Pa_{2n+1}).
\end{equation}
where $\Pi(\point,\eta)$, respectively $ \Pi(\mu,\point)$, stands for the distributions  $Q$ on $
(\RR^d\times\RR^d)$ with marginal $\eta$ w.r.t. the second coordinate and respectively  $\mu$ w.r.t. the first coordinate. 
Let  $\Ka_{n}$ be the Markov transitions on $\RR^d$ defined by the disintegration formulae
\begin{equation}\label{sinhorn-entropy-form-Sch}
\begin{array}{rclcrcl}
\Pa_{2n}(d(x,y))&=&\mu(dx)~\Ka_{2n}(x,dy)\in \Pi(\mu,\pi_{2n})& \mbox{\rm with}&
 \pi_{2n}&:=&\mu \Ka_{2n},\\
 &&&&&&\\
\Pa_{2n+1}(d(x,y))&=&\eta(dy)~\Ka_{2n+1}(y,dx)\in \Pi(\pi_{2n+1},\eta)&\mbox{\rm with}&
\pi_{2n+1}&:=&\eta \Ka_{2n+1}.
\end{array}
\end{equation}
In the above display, we have used the   integral transport notation defined in (\ref{ref-prodm}).

With these notations at hand, the collection of Sinkhorn transitions $(\Ka_n)_{n\geq 0}$  are defined sequentially starting from the reference measure 
$\Pa_0=P=\mu\times \Ka_0$ with $\Ka_0=\Ka$. At every time step $n\geq 0$, given the distribution $\Pa_{2n}$ we choose the transition
$\Ka_{2n+1}$ as the $\Pa_{2n}$-conditional distribution of the first coordinate given the second. Given the distribution $\Pa_{2n+1}$ we choose the transition
 $\Ka_{2(n+1)}$ as the $\Pa_{2n+1}$-conditional distribution of the second coordinate given the first, and so on. Equivalently, the kernel $\Ka_{2n+1}$ is defined in terms of $\Ka_{2n}$ by the conjugate formula
 $$
 \mu(dx)~\Ka_{2n}(x,dy)=(\mu\Ka_{2n})(dy)~\Ka_{2n+1}(y,dx)
 $$
 Similarly, the kernel $\Ka_{2(n+1)}$ is defined in terms of $\Ka_{2n+1}$ by the conjugate formula
 $$
\eta(dy)~\Ka_{2n+1}(y,dx)=(\eta\Ka_{2n+1})(dx)~\Ka_{2(n+1)}(x,dy)
 $$
 
This Gibbs-type structure leads to the fixed-point relations
     \begin{equation}\label{def-pin}
     \begin{array}{l}
 \mu=\mu S_{2n+1}\quad\mbox{\rm and}\quad
\eta=\eta S_{2n}\quad\mbox{\rm with}\quad
(S_{2n},S_{2n+1}):=(\Ka_{2n-1}\Ka_{2n},\Ka_{2n}\Ka_{2n+1})\\
\\
\Longrightarrow\quad \pi_{2n-1}~\Ka_{2n}=\eta
\quad \mbox{\rm and}\quad
\pi_{2n}~\Ka_{2n+1}=\mu.
 \end{array}
\end{equation} 
In the above display, we have used  the composition of Markov operators defined in (\ref{int-compo}); that is, for any $n\geq 0$ we have
$$
(\Ka_n\Ka_{n+1})(x,dz):=\int\Ka_n(x,dy)
\Ka_{n+1}(y,dz).
$$
Note that the fixed point equations (\ref{def-pin}) also imply the time varying evolution equations
\begin{equation}\label{transport-g}
\pi_{2(n-1)}S_{2n}=\pi_{2n}\quad\mbox{\rm and}\quad
\pi_{2n-1}S_{2n+1}=\pi_{2n+1}
\end{equation}
In view of  (\ref{def-pin}), when $n\rightarrow\infty$, we expect Sinkhorn distributions $\pi_{2n}$ to converge towards $\eta$ while $\pi_{2n+1}$ converge towards $\mu$. Using (\ref{def-pin}) and the Markov transport equations (\ref{transport-g})  we check the entropy inequalities
\begin{equation}\label{e12}
  \begin{array}{ccccc}
 \Ha(\pi_{2n}~|~\eta)&\leq& 
   \Ha(\mu~|~\pi_{2n-1})&\leq&  \Ha(\pi_{2(n-1)}~|~\eta),\\
   &&&&\\
     \Ha(\pi_{2n+1}~|~\mu)&\leq &    \Ha(\eta~|~\pi_{2n})&\leq &
         \Ha(\pi_{2n-1}~|~\mu).
 \end{array}
 \end{equation}
 We also have the bridge entropy formulae
\begin{equation}
\begin{array}{l}
\begin{array}{rcccl}
 \Ha(\eta~|~\pi_{2n})&\leq&  \Ha(P_{\mu,\eta}~|~\Pa_{2n})&=&   \Ha(P_{\mu,\eta}|~\Pa_{2n-1})-\Ha(\mu~|~\pi_{2n-1})\\
 &&&&\\
 \Ha(\mu~|~\pi_{2n+1})&\leq&   \Ha(P_{\mu,\eta}~|~\Pa_{2n+1})&=& 
    \Ha(P_{\mu,\eta}~|~\Pa_{2n})-   \Ha(\eta~|~\pi_{2n})
    \end{array}\\
    \\
    \Longrightarrow\qquad
    \Ha(P_{\mu,\eta}~|~\Pa_{2(n+1)})\leq  \Ha(P_{\mu,\eta}~|~\Pa_{2n+1})\leq 
      \Ha(P_{\mu,\eta}~|~\Pa_{2n})\leq  \Ha(P_{\mu,\eta}|~\Pa_{2n-1}).
     \end{array}\label{ent-r1}
\end{equation}
 A detailed proof of  (\ref{e12}) and (\ref{ent-r1}) is provided in Appendix~\ref{review-sinkhorn}.
Thus one expects convergence
$$
\lim_{n\rightarrow\infty} \Ha(P_{\mu,\eta}~|~\Pa_{n})=0\quad \mbox{\rm and}\quad
\lim_{n\rightarrow\infty} \Ha(\eta~|~\pi_{2n})=0=\lim_{n\rightarrow\infty} \Ha(\mu~|~\pi_{2n+1}).
$$
\subsection{Schr\" odinger bridges}

\subsubsection*{Reference flexibility}
A key property of Schrödinger bridges is the flexibility in choosing the reference measure. As shown for instance in Proposition 4.5 in~\cite{ajay-dp-26} (see also Theorem 2.1 in~\cite{nutz} or Section 3.2 and Section 6 in~\cite{adm-24} for the linear-Gaussian case), the bridges  $P_{\mu,\eta}$ and $P^{\flat}_{\eta,\mu}$ may be computed using any Sinkhorn bridge as reference measure. More precisely, for any $n\geq 0$ 
\begin{equation}\label{def-entropy-pb-v2-2}
P_{\mu,\eta}=\argmin_{Q\,\in\, \Pi(\mu,\eta)}\Ha(Q~|~\Pa_{2n})\quad
\mbox{\rm and}\quad
\quad P_{\eta,\mu}^{\flat}=\argmin_{Q\,\in\, \Pi(\eta,\mu)}\Ha(Q~|~\Pa_{2n+1}^{\flat}).
\end{equation}
In particular, taking $n=0$ shows that any result proved for the bridge $P_{\mu,\eta}$ with reference kernel $\Ka_0$ automatically transfers to its conjugate 
 $ P_{\eta,\mu}^{\flat}$ with reference 
$\Ka_1$ after exchanging the roles of $\mu$ and $\eta$.

\subsubsection*{Bridge identities}

Recall that $P_{\mu,\pi_{2n}}$ denotes the Schr\" odinger bridge (\ref{ob-eq}) from $\mu$ to 
$\pi_{2n}$ with respect to the reference distribution $ P$, and similarly
$P_{\pi_{2n+1},\eta}$ is the  bridge from $\pi_{2n+1}$ to 
$\eta$ with respect to the reference distribution $ P$.
The Sinkhorn flow $\pi_n$ defined in (\ref{sinhorn-entropy-form-Sch}) is directly connected to the sequence of Sinkhorn bridges $\Pa_n$ discussed in (\ref{sinhorn-entropy-form-Sch})  through the identities 
\begin{equation}\label{bridge-form-intro}
\Pa_{2n}=P_{\mu,\pi_{2n}}
\quad\mbox{\rm and}\quad\Pa_{2n+1}=P_{\pi_{2n+1},\eta}.
 \end{equation}

For any probability measure $Q$ on $(\RR^d\times\RR^d)$, denote by $Q^{\flat}$ its conjugate measure  $Q^{\flat}(d(x,y)):=Q(d(y,x))$. This terminology is consistent with the dual formulation (\ref{dual-transition-r-d}).
 In this notation, the conjugate bridge of $P_{\mu,\eta}$ satisfies
\begin{equation}\label{ref-back-bridge}
(P_{\mu,\eta})^{\flat}:= \argmin_{Q^{\flat}\,\in\, \Pi(\eta,\mu)}\Ha(Q^{\flat}~|~P^{\flat})=
P^{\flat}_{\eta,\mu}:=(P^{\flat})_{\eta,\mu},
\end{equation}
that is, the conjugate  of the bridge 
from $\mu$ to $\eta$ with reference $P$ coincides with 
the bridge from $\eta$ to $\mu$ with conjugate reference  $P^{\flat}$.

Rewriting (\ref{bridge-form-intro}) using the disintegration  (\ref{sinhorn-entropy-form-Sch}) we obtain
\begin{equation}\label{bridge-form}
\Pa_{2n}=P_{\mu,\pi_{2n}}=P_{\mu,\mu \Ka_{2n}}
\quad\mbox{\rm and}\quad
\Pa_{2n+1}^{\flat}=P^{\flat}_{\eta,\pi_{2n+1}}=P^{\flat}_{\eta,\eta \Ka_{2n+1}}.
 \end{equation}
The fixed point equations (\ref{def-pin}) further imply that
\begin{equation}\label{bridge-form-2}
 P_{\mu,\eta}=P_{\mu,\pi_{2n-1}~\Ka_{2n}}
\quad
\mbox{\rm and}\quad
 P_{\eta,\mu}^{\flat}=P^{\flat}_{\eta,\pi_{2n}\Ka_{2n+1}}.
 \end{equation}
\subsubsection*{Entropy identities}
We now derive identities that will be central to the stability analysis.
First, observe that for any probability measure $\mu$ and any Markov kernel $\Ka$,
\begin{equation}\label{ref-P-mu-muK}
P=\mu\times \Ka\Longrightarrow P_{\mu,\mu\Ka}=P.
\end{equation}
In the above display, we have used  the product notation $\mu\times \Ka$  defined in (\ref{prod-nott}). Next we discuss some consequences of  (\ref{def-entropy-pb-v2-2}).

\begin{itemize}
\item Consider the reference measure $Q$ given by
$$
Q:=\eta\times\Ka_{2n+1}=\Pa_{2n+1}^{\flat}.
$$
Using (\ref{def-entropy-pb-v2-2}), (\ref{ref-P-mu-muK}) and the fixed point relation (\ref{def-pin}) we obtain the equalities
$$
Q_{\eta,\pi_{2n}\Ka_{2n+1}}=Q_{\eta,\mu}=P_{\eta,\mu}^{\flat}
\quad\mbox{\rm and}\quad
 Q=Q_{\eta,\eta\Ka_{2n+1}}=\Pa_{2n+1}^{\flat}.
$$
Thus:
\begin{equation}\label{bridge-form-v2}
\Ha(P_{\mu,\eta}~|~\Pa_{2n+1})=\Ha(P_{\eta,\mu}^{\flat}~|~\Pa_{2n+1}^{\flat})=\Ha(Q_{\eta,\pi_{2n}\Ka_{2n+1}}~|~Q_{\eta,\eta\Ka_{2n+1}}).
 \end{equation}
\item In the same vein, consider the reference measure  $R$ given by
$$
R:=\mu\times\Ka_{2n}=\Pa_{2n}
$$  
Using (\ref{def-entropy-pb-v2-2}),  (\ref{ref-P-mu-muK})  and the fixed point relation  (\ref{def-pin}) we obtain the equalities
$$
R_{\mu,\pi_{2n-1}\Ka_{2n}}=
R_{\mu,\eta}=P_{\mu,\eta}
\quad\mbox{\rm and bywe have}\quad
 R=R_{\mu,\mu\Ka_{2n}}=\Pa_{2n}.
$$
Thus:
\begin{equation}\label{bridge-form-2-V2}
\Ha(P_{\mu,\eta}~|~\Pa_{2n})=
\Ha(R_{\mu,\pi_{2n-1}\Ka_{2n}}~|~R_{\mu,\mu\Ka_{2n}}).
 \end{equation} 
 \end{itemize}
 The entropy identities (\ref{bridge-form-v2}) and (\ref{bridge-form-2-V2}) reveal a direct link between  the stability of Sinkhorn bridges  $\Pa_{2n+1}$ and $\Pa_{2n}$ and the stability of the marginal flow $\pi_{2n}$ and $\pi_{2n-1}$ (towards the invariant measures $\eta$ and $\mu$ as $n\rightarrow\infty$).
 This connection will be exploited repeatedly to derive contraction bounds and exponential convergence estimates for Sinkhorn iterations in later sections.

\subsection{Motivation and related work}

The  fixed point equations (\ref{def-pin}) and the evolution equations (\ref{transport-g}) show that Sinkhorn semigroups belong to the class of time varying Gibbs-type Markov processes sharing a common invariant target measure at each individual step. This type of time varying Markov semigroups arise naturally in a variety of areas in applied mathematics including on nonlinear filtering, physics and molecular chemistry, see for instance~\cite{dpa,douc,saloff-zuniga,saloff-zuniga-2,saloff-zuniga-3}  and references therein. 
The stability analysis of such time varying processes is a broad and notoriously difficult subject. The articles~\cite{douc,saloff-zuniga,saloff-zuniga-2} provide several successful probabilistic and analytical  tools, including coupling methods, spectral analysis as well as functional inequalities such as Nash or log-Sobolev inequalities.  These tools  can be effective in favourable settings, but they interact poorly with the nonlinear conjugacy transformations that govern the Sinkhorn dynamics.

From a different angle, entropic optimal transport techniques, including Schr\"odinger bridges and Sinkhorn algorithm
 have become state-of-the-art tools in generative modeling and machine learning, see for instance~\cite{bortoli-heng,cuturi,kolouri,peyre} and references therein. 
 
Although the Sinkhorn algorithm has an elementary iterative structure, explicit formulas for its updates are rarely available outside finite or Gaussian models. In discrete state spaces, each step corresponds to alternating row and column normalization. In the linear-Gaussian setting, however, Schr\"odinger bridges can be fully characterized through Riccati matrix recursions, yielding closed-form expressions for Sinkhorn transitions see~\cite{adm-24}, Example~\ref{examp-lg-intro-2}, as well as Appendix~\ref{sec-a-stm} on page~\pageref{lin-gauss-sec}. 
  
Much of the existing quantitative theory concerns either finite state spaces~\cite{borwein,sinkhorn-2,sinkhorn-3,soules} or compact settings with bounded costs, often using projective metrics \`a la Birkhoff-Hilbert~\cite{chen-2016,deligiannidis,franklin-1989,marino}.

In contrast, only a few works treat non-compact spaces with unbounded costs~\cite{adm-25,chiarini,durmus} on the exponential convergence of Sinkhorn iterates on non-compact spaces and unbounded cost functions that apply to Gaussian entropic optimal transport models. 

The  articles~\cite{chiarini,durmus} present new quantitative contraction rates for target marginal distributions  with an asymptotically positive log-concavity profile and a reference measure associated with the heat kernel transition on a time horizon $t$.
In this context, the parameter $t$ can be seen as a regularization parameter of the Monge-Kantorovich problem (cf. for instance Section~\ref{sec-ent-maps} and Section~\ref{reg-transp-sec} in the present article). 
The study~\cite{durmus} provides  exponential quantitative estimates for sufficiently large regularization parameter. 
The main conditions in~\cite{durmus} are expressed in terms of the integrated convexity profiles of marginal potentials. Extending ideas from~\cite{conforti-ptrf}  the proofs are based on the propagation of these convexity profiles based on coupling diffusion  by reflection techniques introduced in~\cite{eberle2016reflection}. 

These exponential decays have been recently refined to  apply to all values of the regularization parameter in the more recent article~\cite{chiarini}. In the context of symmetric quadratic costs, the results of~\cite{chiarini} provide entropy bounds for Schrödinger bridges but they require semi-concavity estimates for generally unknown bridge kernels. This restricts their applicability and excludes important cases in which the reference measure arises from Ornstein–Uhlenbeck-type models \eqref{r-OU} or general denoising diffusions \eqref{ref-cnf}.

Article~\cite{eckstein} is also based on an extension of Hilbert’s projective metric for spaces of integrable functions of bounded growth. Building on the results of~\cite{chen-2016,deligiannidis} and extending them to unbounded cost functions,~\cite{eckstein} establishes exponential decay rates for these generalized Hilbert projective metrics with respect to suitable cones of functions. These results hold under tail assumptions requiring the marginal distributions to be sufficiently light relative to the growth of the cost function, and they further imply exponential convergence in total variation distance.   These results also imply the exponential convergence of the total variation distance. In normed finite-dimensional spaces, the main assumption in ~\cite{eckstein} is that the logarithm of the reference transition density grows faster than a polynomial of order $p$ while the tails of both marginals decay faster than $e^{-r^{p+q}}$ for some $q>0$. This framework excludes Gaussian models.

Reference~\cite{adm-25} also develops a semigroup contraction analysis based on Lyapunov techniques
  to prove the exponential convergence of Sinkhorn algorithm on weighted Banach spaces. These Lyapunov approaches also apply  in a variety of situations, ranging from polynomial growth potentials and heavy tailed marginals on general normed
spaces to more sophisticated boundary state space models, including semi-circle transi-
tions, Beta, Weibull, exponential marginals as well as semi-compact models. Last but
not least, the approach in~\cite{adm-25} also allows to consider statistical finite mixture of the above
models, including kernel-type density estimators of complex data distributions arising
in generative modeling.
  
   It is clearly out of the scope of this article to review all the contributions in this field, we simply refer to the book~\cite{peyre} and the review article~\cite{ajay-dp-26} and the  references therein.

\subsection{Log-concavity and Riccati equations}

 In practice, the choice of the reference distribution $P=\mu\times \Ka$ in (\ref{ob-eq})  and (\ref{sinhorn-entropy-form}) is critical. In the context of $t$-entropically regularized transport problems with quadratic costs, the transition $\Ka$
is typically taken to be the heat semigroup (\ref{r-HSG}) at time $t$ (see for instance (\ref{entoo}) and (\ref{ref-tcost})). 
Heat kernels also arise naturally in proximal samplers (Section~\ref{prox-s}), whereas 
diffusion-based generative models (Section~\ref{dgm-sec}) commonly employ Ornstein-Uhlenbeck reference 
transitions. Even more general linear-Gaussian kernels of the form \eqref{ref-cnf} appear throughout the denoising diffusion literature.

As expected for linear-Gaussian transitions $\Ka$ and Gaussian marginals $(\mu,\eta)$ 
Sinkhorn transitions $\Ka_n$ remain linear-Gaussian and can be calculated sequentially 
using least squares and linear regression methods (a.k.a. Bayes' rule).
In this regime, the conditional mean and covariance updates produced by Sinkhorn projections \eqref{sinhorn-entropy-form} coincide with the classical Kalman update. The analysis~\cite{adm-24} develops this filtering/Riccati viewpoint to obtain explicit Schrödinger bridge solutions via discrete algebraic Riccati equations (DARE), together with refined convergence estimates for the Sinkhorn algorithm across a broad class of linear-Gaussian models.

As in Kalman filtering theory, the covariance flows associated with Sinkhorn iterations satisfy offline Riccati matrix difference equations. Their stability properties and stationary solutions are well understood; see, for instance, \cite{adm-24,dh-23}. Section~\ref{basics-sec} reviews Riccati flows in the Sinkhorn context, providing closed-form descriptions of fixed points and sharp exponential convergence estimates.

For more general marginal measures $(\mu,\eta)$ Sinkhorn iterations are based on nonlinear conditional/conjugate transformations and exact finite-dimensional solutions cannot be computed. 
Stability analyses for these models, such as those in~\cite{chiarini,durmus,lee} typically rely on controlling the Hessians of log-densities (cf. Corollary 9.3.2 and Corollary 5.7.2 in~\cite{bgl}, as well as Examples~\ref{nabla2-logsob}, \ref{nabla2-logsob-v2}, and \ref{examp-logsob}). Hessian bounds combined with Brascamp-Lieb and Cramér-Rao inequalities connect the log-concavity of Schrödinger bridges and Sinkhorn transitions to their conditional covariances. Iterating these inequalities along the Sinkhorn recursion yields Riccati matrix difference equations (see~\cite{adm-24,chewi,durmus} and Appendix~\ref{sec-a-stm}).

  \subsection{Main contributions}

The article develops a quantitative contraction framework for Schrödinger and Sinkhorn bridges based on transportation-cost inequalities and Riccati matrix difference equations. Our approach combines logarithmic Sobolev and Talagrand-type inequalities to obtain explicit entropy and Wasserstein contraction bounds for Sinkhorn bridge measures, entropic optimal transport plans, and the associated Markov transport maps. A key feature of the analysis is the interplay between transport inequalities and matrix Riccati difference equations arising in filtering and stochastic control.

One of the main results is a novel entropic continuity theorem that applies to {\it general classes of cost functions and to general reference measures} (see Theorem~\ref{th1-intro} and Theorem~\ref{lem-Int-log-sob-talagrand}). A key feature of this result is that the exponential stability of Sinkhorn’s algorithm follows directly from the quantitative entropic estimates developed here (see Theorem~\ref{th2-intro}, as well as the contraction bounds \eqref{f1-c}, \eqref{f2-c}, and Theorem~\ref{theo-imp}).

The regularity assumptions on the target marginals are intentionally flexible: we require only that one measure be log-concave at infinity and that the other satisfy a quadratic transportation inequality. When both probability measures are log-concave at infinity, we obtain even sharper exponential rates based on second-order linear recurrence inequalities (Theorems~\ref{th2-intro} and \ref{th3-intro}). 

As a main application, we specialize the theory to linear-Gaussian {\it reference} transitions, where the Gaussian structure permits sharp constants, refined exponential decay rates, and continuity estimates for Schrödinger bridges, Sinkhorn iterates, barycentric projections, conditional covariances, and proximal sampler semigroups.
We extend the Riccati matrix analysis of~\cite{adm-24}, originally developed for linear-Gaussian models to the general setting of log-concave target distributions. In the linear-Gaussian case, Theorem 3.1 and Corollary 3.5 in~\cite{adm-24} show that both Sinkhorn and Schrödinger bridges admit closed-form expressions in terms of Riccati equations (see also \eqref{Sigma-Gauss}, Example~\ref{examp-bary}, and Appendix~\ref{sec-a-stm}). The exponential rates of the Sinkhorn algorithm are likewise characterized by the fixed points of Riccati matrix equations (Section 5 in~\cite{adm-24}). These stability results and explicit formulas for Riccati fixed points \cite{adm-24,dh-23} continue to play an essential role in the log-concave setting.

Building on this framework, Section~\ref{sec-ent-maps} provides a new estimation of barycentric projections (entropic maps) of Schrödinger bridges in terms of the fixed point of Riccati maps (Theorem~\ref{theo-baryp}). For entropically regularized problems with quadratic costs, these results yield extended and refined versions of Caffarelli’s contraction theorem~\cite{caffarelli} and of the contraction bounds of Chewi and Pooladian~\cite{chewi}, now applicable to general cost functions of the form \eqref{def-W} and to non-commuting covariance structures. As shown in Example~\ref{examp-gauss-sh}, these estimates are sharp and reduce to exact matrix identities in the linear–Gaussian case.

In the context of the Sinkhorn algorithm, we further obtain non-asymptotic exponential decay estimates expressed in terms of the Riccati fixed point (Theorem~\ref{ineq-theo-impp-2-i-th} and Corollary~\ref{cor-impp}).

A major strength of our approach is its applicability to a broad class of models in machine learning and artificial intelligence. Section~\ref{sec-illu} illustrates the impact of our results in regularized entropic transport, the stability analysis of proximal samplers, diffusion-based generative models, and diffusion flow-matching methods.

The rest of the article is organized as follows.

Section~\ref{sec-statements} provides a brief summary of the main results.
Section~\ref{sec-ent-continuity} introduces the entropic continuity theorem, 
while Section~\ref{sinkcontract-sec} presents several contraction theorems for Sinkhorn bridges.
Section~\ref{sec-ent-maps} develops a conditional covariance representation for the gradient of entropic maps, including extended and refined versions of Caffarelli’s contraction theorem. 

Section~\ref{sec-illu} illustrates the implications of our results for regularized entropic transport, proximal samplers, diffusion-based generative models, and diffusion flow-matching methods.

Section~\ref{transport-sec} introduces the basic notions of divergences and transportation costs used throughout the article.
Section~\ref{sec-continuity-sb} establishes quantitative upper and lower bounds for Schrödinger bridges, including estimates for Markov transport maps.
Section~\ref{sinkhorn-sec} examines the stability properties of Sinkhorn bridges.
Section~\ref{sec-lingauss} presents contraction theorems for Sinkhorn semigroups and several estimates for entropic and Sinkhorn maps in the setting of linear-Gaussian reference transitions.

Finally, the appendix contains a brief review of Sinkhorn semigroups, various gradient and Hessian formulas for potential functions associated with Sinkhorn bridges, and the proofs of several technical results used throughout the article.

\section{Description of the models}

\subsection{Some basic notation}\label{basics-sec}

\subsubsection*{Integral operators}

Let $\Ma(\RR^d)$ be the set of locally bounded signed measures on $\RR^d$, and
let $\Ma_b(\RR^d)\subset \Ma(\RR^d)$ the subset of bounded measures. Let $\Ma_1(\RR^d)\subset \Ma_b(\RR^d)$ be the convex subset of probability measures on $\RR^d$. 
For a measure $\mu(dx)$, a function $f(x)$ and an integral operator $\Ka(x,dy)$, whenever they exist we denote by  $\mu(f)$,  $(\mu\Ka)(dy)$ and $\Ka(f)(x)$  the Lebesgue integral,  the measure and the function defined by
\begin{equation}\label{ref-prodm}
\mu(f):=\int_{\RR^d} f(x)~\mu(dx)\qquad
(\mu\Ka)(f):=\mu(\Ka(f))
\qquad\Ka(f)(x):=(\delta_x\Ka)(f).
\end{equation}
For indicator functions $f=1_{A}$ sometimes we write $\mu(A)$ instead of $\mu(1_A)$.
A  measure $\nu_1\in \Ma(\RR^d)$ is said to be absolutely continuous with respect to another measure $\nu_2\in \Ma(\RR^d)$ and we write $\nu_1\ll \nu_2$
as soon as $\nu_2(A)=0\Longrightarrow \nu_1(A)$ for any measurable subset $A\subset\RR^d$. Whenever $\nu_1\ll \nu_2$, we denote by ${d\nu_1}/{d\nu_2}$ the Radon-Nikodym derivative of $\nu_1$ w.r.t. $\nu_2$.
   
For indicator functions $f=1_{A}$ sometimes we write $\Ka(x,A)$ instead of $\Ka(1_A)(x)$.
For bounded integral operator  $\Ka_1(x,dy)$ from $\RR^d$ into $\RR^d$ and  $\Ka_2(y,dz)$ from $\RR^d$ into $\RR^d$,  we also denote by $(\Ka_1\Ka_2)$ the bounded integral operator from $\RR^d$ into $\RR^d$ defined by the integral composition formula
\begin{equation}\label{int-compo}
(\Ka_1\Ka_2)(x,dz):=\int\Ka_1(x,dy).
\Ka_2(y,dz)
\end{equation}

We also consider the product measures
\begin{equation}\label{prod-nott}
(\mu\times \Ka)(d(x,y)):=\mu(dx)\, \Ka(x,dy)
\quad\mbox{\rm and set}\quad
(\mu\times \Ka)^{\flat}(d(x,y)):=\mu(dy)\, \Ka(y,dx).
\end{equation}
The tensor product $\mu_1\otimes \mu_2$ of measures $\mu_1$ and $\mu_2$ is denoted by
$$
(\mu_1\otimes \mu_2)(d(x,y))=\mu_1(dx)\mu_2(dy)
\quad\mbox{\rm and}\quad \mu^{\otimes 2}=\mu\otimes\mu.
$$

Assume that $\delta_x\Ka\ll \mu\Ka$ for any $x\in\RR^d$. In this situation we have
\begin{equation}\label{dual-transition-r-d}
P=
\mu\times \Ka\Longrightarrow P^{\flat}=
(\mu\times \Ka)^{\flat}=(\mu\Ka)\times\Ka^{\ast}_{\mu}.
\end{equation}
with the dual (a.k.a. backward) Markov transition
\begin{equation}\label{dual-transition-r}
\Ka^{\ast}_{\mu}(y,dx):=\mu(dx)~\frac{d\delta_x\Ka}{d\mu\Ka}(y)
\Longleftrightarrow
\mu\left(f\Ka(g)\right)=(\mu\Ka)(~\Ka^{\ast}_{\mu}(f)~g)\Longrightarrow \mu \Ka\Ka^{\ast}_{\mu}=\mu.
\end{equation}

\subsubsection*{Differential operators}\label{def-hessian-grad}

We let $\nabla f(x)=\left[\partial_{x_i}f(x)\right]_{1\leq i\leq d}$ be the gradient column vector associated with some smooth function $f(x)$ from $\RR^d$ into $\RR$.
 Given some smooth function $h(x)$ from $\RR^d$ into $\RR^p$
we denote by $\nabla h=\left[\nabla h^1,\ldots,\nabla h^p\right]\in \RR^{n\times p}$ the gradient matrix associated with the column vector
 function $h=(h^i)_{1\leq i\leq p}$, for some $p\geq 1$. For instance, for a given matrix 
 $A\in\RR^{q\times p}$ we have $Ah(x)\in\RR^q$ and therefore
\begin{equation}\label{rule-ah}
 \nabla (A h(x))=(\nabla h(x))~A^{\prime}.
\end{equation}
Given a smooth  function $f(x,y)$ from $\RR^d\times\RR^d$ into $\RR$ and $i=1,2$
  we denote by $\nabla_i$ the gradient operator defined by
   $$
 \nabla_if=\left[
 \begin{array}{c}
 (\nabla_if)^1\\
 \vdots\\
  (\nabla_if)^d
 \end{array}
 \right] \quad\mbox{\rm with}\quad
  (\nabla_1f)^j:= \partial_{x^j}f \quad\mbox{\rm and}\quad
  (\nabla_2f)^j:= \partial_{y^j}f .
 $$
For smooth $p$-column valued functions $F(x,y)=(F^j(x,y))_{1\leq j\leq p}$ we also set
$$
  \nabla_i F=\left[ \nabla_i F^1,\ldots,
  \nabla_iF^p\right].
   $$ 
In this notation, the second order differential operators  $\nabla_i \nabla_j$ on function $f(x,y)$ from $\RR^d\times\RR^d$ into $\RR$  are defined by $(\nabla_i \nabla_j) (f):=  \nabla_i (\nabla_j(f))$ and we set $
  \nabla_i^2:=  \nabla_i \nabla_i$ when $i=j$.

\subsubsection*{Positive definite matrices}
Let $\Ga l_d$ be the general linear group of $(d\times d)$-invertible matrices, $\Sa^0_d$ the set 
 of (symmetric) positive semi-definite matrices and $\Sa^+_d\subset \Sa_d^0$ the subset of positive definite matrices. We denote by $\ell_{\tiny min}(v)$ and $\ell_{\tiny max}(v)$ the minimal and the maximal eigenvalues, respectively, of a symmetric matrix $v\in\RR^{d\times d}$ for some $d\geq 1$.
The spectral norm is defined by $\Vert v\Vert_2=\sqrt{\ell_{\tiny max}(v^{\prime}v)}$, with  $v^{\prime}$ the transpose. 

We sometimes use the L\" owner partial ordering notation $v_1\geq v_2$ to mean that a symmetric matrix $v_1-v_2$ is positive semi-definite (equivalently, $v_2 - v_1$ is negative semi-definite), and $v_1>v_2$ when $v_1-v_2$ is positive definite (equivalently, $v_2 - v_1$ is negative definite). Given $v\in\Sa_d^+$ we denote by $v^{1/2}$  the principal (unique) symmetric square root.

We also quote a local Lipschitz property of the square root function $v^{1/2}$ on (symmetric) definite positive matrices. For any $u,v\in\Sa_{r}^+$ we have the Ando-Hemmen inequality
\begin{equation}\label{square-root-key-estimate}
\Vert u-v\Vert \leq \left[\ell^{1/2}_{min}(u)+\ell^{1/2}_{min}(v)\right]^{-1}~\Vert u-v\Vert,
\end{equation}
for any unitary invariant matrix norm (e.g. the spectral, or Frobenius norm); see e.g.~\cite{higham} and~\cite{hemmen}.

The geometric mean $u~\sharp~ v$ of {two} positive definite matrices $u,v\in\Sa^+_d$ is defined by
\begin{equation}\label{sym-sharp}
u~\sharp~v=v~\sharp~u:=v^{1/2}~ \left(v^{-1/2}~u~v^{-1/2}\right)^{1/2}~v^{1/2}.
\end{equation}
For commuting matrices $uv=vu$ we recall that $u~\sharp~v=(uv)^{1/2}$.
 
With a slight abuse of notation, we denote by $I$ the $(d\times d)$-identity matrix as well as the identity integral operator.
We also denote by $0$ the null $(d\times d)$-matrix and the null $d$-dimensional vector, for any choice of the dimension $d\geq 1$. 

We usually represent points $x=(x^j)_{1\leq j\leq d}\in\RR^d$ by $d$-dimensional column vectors and $1 \times d$ matrices. In this notation, the Frobenius and the spectral norm $\Vert x\Vert_2=\sqrt{x^{\prime}x}$ coincides with the Euclidean norm. For $a,b\in\RR$ we set $a\wedge b:=\min{(a,b)}$ and  $a\vee b:=\max{(a,b)}$.

\subsection{Riccati matrix equations}
 We associate with some given $\varpi\in\Sa^+_d$  the Riccati map
\begin{equation}\label{ricc-maps-def}
s\in \Sa^0_d\mapsto \mbox{\rm Ricc}_{\varpi}(s):=(I+(\varpi+s)^{-1})^{-1}\in \Sa^+_d.
\end{equation}
As shown in~\cite{adm-24} the unique positive definite fixed point $ r_{\varpi}= \mbox{\rm Ricc}_{\varpi}(r_{\varpi})$ is given by 
\begin{equation}\label{def-fix-ricc-1}
(I+\varpi^{-1})^{-1}< r_{\varpi}:=-\frac{\varpi}{2}+\left(\varpi+\left(\frac{\varpi}{2}\right)^2\right)^{1/2}< I.
\end{equation}
Note the increasing property
$$
v_1\leq v_2\Longrightarrow \mbox{\rm Ricc}_{\varpi}(v_1)\leq \mbox{\rm Ricc}_{\varpi}(v_2).
$$
Let $r_n$ be the solution of the Riccati matrix difference equation
\begin{equation}\label{def-rn}
r_n=\mbox{\rm Ricc}_{\varpi}(r_{n-1}),
\end{equation}
starting from some initial condition $r_0\in \Sa^0_d$.
Note that for any $n\geq p\geq 1$ and $v\in\Sa^d_0$ we have the estimates
\begin{equation}\label{ricc-maps-incr}
 \mbox{\rm Ricc}_{\varpi}^p(0)\leq \mbox{\rm Ricc}_{\varpi}^n(0)\leq \mbox{\rm Ricc}_{\varpi}^n(v)
\leq \mbox{\rm Ricc}_{\varpi}^{n-1}(I)\leq  \mbox{\rm Ricc}_{\varpi}^{p-1}(I)\leq I.
\end{equation}
These matrix equations belong to the class of discrete algebraic Riccati equations (DARE), {and} no analytical solutions are available for general models.  
In our context, we recall that the unique positive definite fixed point $ r_{\varpi}= \mbox{\rm Ricc}_{\varpi}(r_{\varpi})$ is given in closed form by the formula (\ref{def-fix-ricc-1}).
 In addition, we have the formulae
\begin{equation}\label{prop-fp}
  r_{\varpi}+ r_{\varpi} \varpi^{-1} r_{\varpi}=I\quad \mbox{\rm and}\quad
 I<  r_{\varpi}^{-1}=I+u_{\varpi}^{-1}< I+\varpi^{-1}
 \quad \mbox{\rm with}\quad
 u_{\varpi}:=
\varpi+ r_{\varpi}.
\end{equation}
 There also exists some constant $c_{\varpi}$
such that for any $n\geq 1$ we have
\begin{equation}\label{cv-ricc-intro}
 \Vert r_{n}-r_{\varpi}\Vert_2\leq c_{\varpi}~\delta_{\varpi}^n~
\Vert r_0-r_{\varpi}\Vert_2\quad \mbox{\rm with}\quad
\delta_{\varpi}:=(1+\ell_{\tiny min}(u_{\varpi}))^{-2}.
\end{equation}
By Proposition A.6 in~\cite{adm-24} the constant $c_{\varpi}$ can be chosen so that
$$
c_{\varpi}\leq \frac{1}{\delta_{\varpi}}~\left(\frac{\Vert\varpi\Vert_2~\Vert\varpi^{-1}\Vert_2}{1+\Vert\varpi^{-1}\Vert_2}\right)^2~
\left((1+\Vert u_{\varpi}\Vert_2)(1+\Vert u_{\varpi}^{-1}\Vert_2)\right)^2.
$$

\begin{examp}
Lemma 4.3 in~\cite{dh-23} provides closed-form solutions of Riccati difference equations. For instance, for one dimensional models the solution of (\ref{def-rn}) is given by
$$
(r_{n}-r_{\varpi})=(r_{0}-r_{\varpi})~\frac{(\varpi+2r_{\varpi})~\delta_{\varpi}^n}{(r_0+\varpi+r_{\varpi})(1-\delta_{\varpi}^n)+(\varpi+2r_{\varpi})~\delta_{\varpi}^n},
$$
with the positive fixed point $r_{\varpi}$ defined in (\ref{def-fix-ricc-1}) and the exponential decay parameter $\delta_{\varpi}$ defined in (\ref{cv-ricc-intro}).
\end{examp}

\subsection{Potential functions}
 
 Equip $\RR^d$ with the Lebesgue measure $\lambda(dx)$ and
denote by  $\lambda_U$ and  $K_{W}$  Boltzmann-Gibbs measures and 
Markov transitions  of the following form
\begin{equation}
\lambda_U(dx):=e^{-U(x)}~\lambda(dx)\quad \mbox{\rm and}\quad
K_W(x,dy):=e^{-W(x,y)}~\lambda(dy),\label{ref-KW}
\end{equation}
for some potential functions $U(x)$ and $W(x,y)$ on $\RR^d$ and on $(\RR^d\times\RR^d)$. We also set 
\begin{equation}\label{ref-Wflat}
W^{\flat}(y,x):=W(x,y),
\end{equation}
 and denote by $  \Sigma_{W}$ the conditional covariance function
\begin{equation}\label{cov-mat-W}
  \Sigma_{W}(x):=\frac{1}{2}\int~
  K_W(x,dy)K_W(x,dz)~(y-z)~(y-z)^{\prime}.
\end{equation}
Consider the Schr\" odinger bridge $P_{\mu,\eta}$ with reference measure $P:=(\mu\times \Ka)$ associated with the probability measures and the reference Markov transition
\begin{equation}\label{def-triple}
(\mu,\eta,\Ka)=(\lambda_U,\lambda_V,K_W),
\end{equation}
for some potential functions $(U,V,W)$.
  The prototypical reference transition $K_W$ we have in mind is the
linear Gaussian transition  associated with the potential function
\begin{equation}\label{def-W}
W(x,y)=-\log{g_{\tau}(y-(\alpha+\beta x))},
\end{equation}
for a positive definite matrix $\tau>0$, some  $\alpha\in\RR^d$ and some invertible matrix $\beta\in\RR^{d\times d}$. In the above display, $g_{\tau}$ stands for the density of a centered Gaussian random variable with covariance matrix $\tau$. We underline that the linear Gaussian transition (\ref{def-W}) encapsulates all continuous time Gaussian models used in machine learning applications of Schr\" odinger bridges (see Section~\ref{sec-illu} as well as~\cite{adm-24} and references therein).

In this context, the Schr\" odinger bridges discussed in (\ref{ob-eq}) and (\ref{ref-back-bridge}) take the following product form
\begin{equation}\label{sbr-ref}
P_{\mu,\eta}=\mu~\times~K_{\WW_{\mu,\eta}}\quad \mbox{\rm and}\quad
P_{\eta,\mu}^{\flat}=\eta~\times~K_{\WW^{\flat}_{\eta,\mu}},
\end{equation}
for some potential functions  $\WW_{\mu,\eta}$ and $\WW^{\flat}_{\eta,\mu}$  (see Section~\ref{sec-sbm}, as well as Section~\ref{sec-pf} and (\ref{log-densities-fp})).  
For clarity of exposition, we write $\Sigma_{\mu,\eta}$ and $\Sigma^{\flat}_{\eta,\mu}$  instead of $\Sigma_{\WW_{\mu,\eta}}$ and $\Sigma_{\WW^{\flat}_{\eta,\mu}}$ the conditional covariance matrices (\ref{cov-mat-W}) associated with the transitions
$K_{\WW_{\mu,\eta}}
$ and $K_{\WW^{\flat}_{\eta,\mu}}$. 

The entropic maps $(\ZZ_{\mu,\eta},\ZZ^{\flat}_{\eta,\mu})$ are the random maps defined by
\begin{equation}\label{e-maps}
K_{\WW_{\mu,\eta}}(x,dy)=\PP(\ZZ_{\mu,\eta}(x)\in dy)\quad\mbox{and}\quad
K_{\WW^{\flat}_{\mu,\eta}}(x,dy)=\PP(\ZZ^{\flat}_{\eta,\mu}(x)\in dy).
\end{equation}

\subsection{Regularity conditions}\label{marg-reg-conditions}

This section introduces the structural and regularity conditions imposed on the marginal potentials  $(U,V)$
 throughout the paper. These assumptions are formulated in terms of curvature bounds on the Hessians of 
$(U,V)$, together with functional inequalities of transportation type. We first present general curvature conditions allowing for non-uniform convexity, in particular convexity at infinity, and relate them to logarithmic Sobolev and quadratic transportation inequalities. We then specialize to strongly convex models. Several illustrative examples are provided to clarify the scope of the assumptions and to connect them with classical models considered in the literature.

 \subsubsection*{Curvature conditions}
 Next regularity conditions depend on the quadratic transportation cost inequality $\TT_2(\rho)$ with constant $\rho>0$ (cf. Section~\ref{sec-tineq}) as well as on  the curvature of the  functions $U,V$.\\

\noindent 
 {\it $C_V^-(U):$ There exists positive definite matrix
$u_+,v_+>0$, $\overline{u}_+\in\Sa_d^0$ and some  $\delta\geq 0$ s.t.
\begin{equation}\label{ex-lg-cg-intro}
\nabla^2 U(x)\geq 
u_+^{-1}~1_{\Vert x\Vert_2 \geq  \delta}-\overline{u}_+ 
~1_{\Vert x\Vert_2< \delta}\quad\mbox{and $\lambda_V$ satisfies the $\TT_2(\Vert v_+\Vert)$-inequality.}
\end{equation}}

\noindent {\it $C_V(U):$    There exists positive definite matrices $u_+,v_+>0$ such that 
\begin{equation}\label{UV-compU2V}
\nabla^2 U\geq u^{-1}_+\quad\mbox{and $\lambda_V$ satisfies the $\TT_2(\Vert v_+\Vert)$-inequality.}
\end{equation}}
 Note that $C_V(U)\Longleftrightarrow C_V^-(U)$ with $\delta=0$.
Conditions $C_U(V)$ and $C_U^-(V)$ is defined as (\ref{UV-compU2V}) and (\ref{ex-lg-cg-intro}) switching the role of the potential functions. 

The Hessian condition in the l.h.s. of (\ref{ex-lg-cg-intro}) can be seen as a "uniform convexity property at infinity". This condition
(and thus the one in the l.h.s. of (\ref{UV-compU2V})) ensures that
$\lambda_U$ satisfies the  log-Sobolev inequality $LS(\rho)$ and thus the $\TT_2(\rho)$-inequality for some $\rho>0$ (cf. Section~\ref{sec-tineq}). 
Thus, $C_V^-(U)$ and $C_U^-(V)$ are met as soon as $U$ and $V$ are  "uniformly convex at infinity". 
To be specific, consider the condition:\\

\noindent 
 {\it $C_0^-(U,V):$ There exists positive definite matrix
$u_+,v_+>0$, $\overline{u}_+,\overline{v}_+\in\Sa_d^0$ and  $\delta_u,\delta_v\geq 0$ s.t.
\begin{equation}\label{ex-lg-cg-intro-uv-lls}
\begin{array}{rcl}
\nabla^2 U(x)&\geq &
u_+^{-1}~1_{\Vert x\Vert_2 \geq  \delta_u}-\overline{u}_+ 
~1_{\Vert x\Vert_2< \delta_u},\\
&&\\
\nabla^2 V(y)&\geq &
v_+^{-1}~1_{\Vert y\Vert_2 \geq  \delta_v}-\overline{v}_+ 
~1_{\Vert y\Vert_2< \delta_v}.
\end{array}
\end{equation}}

Whenever $C_0^-(U,V)$ is satisfied, the probability measure $\lambda_U$  satisfies the  log-Sobolev inequality $LS(\rho_u)$ and thus the $\TT_2(\rho_u)$ inequality; while  $\lambda_V$ satisfies the  log-Sobolev inequality $LS(\rho_v)$ and thus the $\TT_2(\rho_v)$ inequality, with some parameters \begin{equation}\label{rho-uv}
\rho_u=\rho(u_+,\overline{u}_+,\delta_u)\quad \mbox{\rm and}\quad 
\rho_v=\rho(v_+,\overline{v}_+,\delta_v).
\end{equation}
 This shows that
$$
C_0^-(U,V)\Longrightarrow (C_V^-(U)~\&~C_U^-(V)).
$$

 \subsubsection*{Strongly convex models}

The next regularity conditions correspond to strongly convex potential functions.\\

\noindent  {\it $C_0(U,V):$  There exists some positive definite matrices $u_+,v_+>0$ such that
     \begin{equation}\label{UV-d0}
\nabla^2 U\geq u^{-1}_+\quad\mbox{and}\quad
\nabla^2 V\geq v^{-1}_+.
\end{equation}}
\noindent  {\it $C_1(U,V):$ Condition $C_0(U,V)$ is satisfied for some positive definite matrices $u_+,v_+>0$. In addition, there exists positive semi-definite matrices  $u_{-},v_{-}\in\Sa_d^0\cup\{0\}$ such that
\begin{equation}\label{UV-d2}
\nabla^2 U\leq u^{-1}_-\quad\mbox{\rm and}\quad
\nabla^2 V\leq v^{-1}_-.
\end{equation} }

When $s=0$ we use the convention $s^{-1}=\infty\times I$, so that $\nabla^2 U\leq u^{-1}_-$ is for instance trivially met when $u_-=0$. When 
$u_-=0=v_-$  condition $C_1(U,V)$ coincides with $C_0(U,V)$.

Also observe that $$
C_1(U,V)\Longrightarrow
C_0(U,V)\Longleftrightarrow (C_U(V)~\&~ C_V(U))\Longrightarrow (C_V^-(U)~\&~C_U^-(V)).
$$

When both $u_-,v_-$ are positive definite matrices, condition (\ref{UV-d2}) is much stronger than the conventional convex condition (\ref{UV-d0}).

 \subsection{Some illustrations}

 Condition $C_1(U,V)$  with $u_-,v_->0$ has been used in several articles~\cite{chewi,conforti-ptrf,durmus}.
It encapsulates Gaussian potential functions
\begin{equation}\label{def-U-V}
U(x)=-\log{g_{u}(x-m)}
\quad \mbox{\rm and}\quad V(x):=-\log{g_{v}(x-\overline{m})},
\end{equation}
for some given parameters
$
(m,\overline{m})\in(\RR^d\times\RR^{d})$ and $ 
 (u,v)\in (\Sa^+_d\times \Sa^+_{d})
$. In this context (\ref{UV-d0}) as well as (\ref{UV-d2}) trivially hold with the covariances matrices
$$
(u_-,v_-)=(u,v)=(u_+,v_+).
$$

Consider the
linear-Gaussian reference transitions $K_W$ associated with some potential $W$ of the form (\ref{def-W}) for some parameters $(\alpha,\beta,\tau)$ and set $\cchi:=\tau^{-1}\beta$.
In this context, whenever condition $C_1(U,V)$ holds  for some $u_+,v_+>0$ and some positive semi-definite matrices $u_{-},v_{-}\in\Sa_d^+\cup\{0\}$ we set

\begin{equation}\label{def-over-varpi}
\begin{array}{rcl}
\varpi_0^{-1}&:=&v^{1/2}_-~\left(\cchi~u_+~\cchi^{\prime}\right)~v^{1/2}_-\\
\varpi_1^{-1}&:=&u_-^{1/2}~\left(\cchi^{\prime}~v_+~\cchi\right)~u_-^{1/2},
\end{array}
\quad\&\quad
\begin{array}{rcl}
\overline{\varpi}_0^{-1}&:=&v^{1/2}_+~(\cchi~u_-~\cchi^{\prime})~v^{1/2}_+\\
\overline{\varpi}_1^{-1}&:=&u_+^{1/2}~(\cchi^{\prime}~v_-~\cchi)~u_+^{1/2}.
\end{array}
\end{equation}
\begin{examp}\label{examp-lg-intro-2}
For 
linear-Gaussian models (\ref{def-W}) 
and (\ref{def-U-V}) we have
\begin{equation}\label{def-over-varpi-g}
\begin{array}{l}
(u,v)=(u_{\pm},v_{\pm})\\
\\
\Longrightarrow
(\overline{\varpi}_0,
\overline{\varpi}_1)=({\varpi}_0,{\varpi}_1)=\left(v^{-1/2}~\left(\cchi~u~\cchi^{\prime}\right)^{-1}~v^{-1/2}, u^{-1/2}~(\cchi^{\prime}~v~\cchi)^{-1}~u^{-1/2}
\right)
\end{array}\end{equation}
In this situation, applying Theorem 3.1 and Corollary 3.5  in~\cite{adm-24} the conditional covariances are constant functions defined for any $x\in\RR^d$ by
\begin{equation}\label{Sigma-Gauss}
\Sigma_{\mu,\eta}(x)=\sigma_{\mu,\eta}:=v^{1/2} ~r_{\overline{\varpi}_0}~v^{1/2}\quad
\mbox{\rm and}
\quad
\Sigma^{\flat}_{\eta,\mu}(x)=\sigma^{\flat}_{\eta,\mu}=u^{1/2} ~r_{\overline{\varpi}_1}~u^{1/2}
\end{equation}
In this context, the bridge potentials (\ref{sbr-ref}) are given in closed form by the Gaussian densities
\begin{eqnarray}
\WW_{\mu,\eta}(x,y)&=&-\log{g_{\,\sigma_{\mu,\eta}}\left(y-\left(\overline{m}+\sigma_{\mu,\eta}~\cchi~ (x-m)\right)\right)}\nonumber\\
\WW^{\flat}_{\eta,\mu}(x,y)&=&-\log{g_{\,\sigma^{\flat}_{\eta,\mu}}\left(y-\left(m+\sigma^{\flat}_{\eta,\mu}~\cchi^{\prime}~ (x-\overline{m})\right)\right)}.\label{ref-WW-Gauss}
\end{eqnarray}
\end{examp}

Based on our conventions, we underline the following situations:

$\bullet$ Condition $C_1(U,V)$ with  $u_-=0$ and $v_->0$ takes the form
\begin{equation}\label{convention-uv}
u_+^{-1}\leq \nabla^2 U\quad\mbox{\rm and}\quad
v_+^{-1}\leq\nabla^2 V\leq v^{-1}_-.
\end{equation}
In this case, we use the convention $\overline{\varpi}_0^{-1}=0=\varpi_1^{-1}$ and for any $s\in\Sa^0_d$ we set
\begin{equation}\label{convention-u}
\mbox{\rm Ricc}_{\overline{\varpi}_0}(s)=I=r_{\overline{\varpi}_0}
\quad
\mbox{and}\quad\mbox{\rm Ricc}_{{\varpi}_1}(s)=I=r_{\varpi_1}.
\end{equation}
Note that in this situation we have $u_+^{-1}>0$  so that $\varpi_0^{-1}>0$ and $\overline{\varpi}_1^{-1}>0$.

$\bullet$  Condition $C_1(U,V)$ with  $u_->0$ and $v_-=0$ takes the form
\begin{equation}\label{convention-vu}
u_+^{-1}\leq \nabla^2 U\leq u_-^{-1}\quad\mbox{\rm and}\quad
v_+^{-1}\leq\nabla^2 V.
\end{equation}
In this case, we use the convention $\overline{\varpi}_1^{-1}=0=\varpi_0^{-1}$ and for any $s\in\Sa^0_d$
 we set
\begin{equation}\label{convention-v}
\mbox{\rm Ricc}_{\overline{\varpi}_1}(s)=I= r_{\overline{\varpi}_1}
\quad
\mbox{and}\quad\mbox{\rm Ricc}_{{\varpi}_0}(s)=I=r_{\varpi_0}.
\end{equation}
 Note that in this situation we have $v_+>0$ we have $\varpi_1^{-1}>0$ and $\overline{\varpi}_0^{-1}>0$.

 \section{Statement of some main results}\label{sec-statements}
 \subsection{An entropic continuity theorem}\label{sec-ent-continuity}

One of the main contributions of this article is a new and flexible entropic continuity theorem that applies to broad classes of reference measures. Our goal is to transfer regularity properties of the Markov transport maps 
$$
\eta\mapsto\eta\Ka
$$
to the corresponding Schr\" odinger entropic projection map
\begin{equation}\label{ot-eq-etaK-i}
\eta\mapsto P_{\mu,\eta\Ka}:= \argmin_{Q\,\in\, \Pi(\mu,\eta\Ka)}\Ha(Q~|~P)\quad \mbox{\rm with the reference}\quad
P:=\mu\times \Ka=P_{\mu,\mu\Ka}.
\end{equation}

Understanding the regularity and stability of transport maps of the form $\eta \mapsto \eta\Ka$ is a central topic in applied probability, particularly in the study of the stability of stochastic processes; see, for example, \cite{dpa,bgl,dm-03,douc} and the references therein.

With these notations, our first main result establishes an entropic continuity principle, stated informally as follows.

\begin{theo}\label{th1-intro}
When the Markov kernel $\Ka$ in the reference measure (\ref{ot-eq-etaK-i}) is regular enough, there exists constants 
$a,b>0$ such that for any probability measures $\mu$ and $\eta$,
\begin{equation}\label{ref-Hep-D}
~\Da_{2}(\eta \Ka,\mu \Ka)^2\leq  a~
 \Ha(P_{\mu,\eta\Ka}~|~P_{\mu,\mu\Ka})\leq b
~\Da_{2}(\eta,\mu)^2.
\end{equation}
Moreover, if $\mu$ satisfies a quadratic transportation cost inequality, then there exists  $\varepsilon>0$ such that for any   probability measures $\eta$,
\begin{equation}\label{ref-Hep}
  \Ha(\eta\Ka~|~\mu\Ka)\leq
 \Ha(P_{\mu,\eta\Ka}~|~P_{\mu,\mu\Ka})
 \leq \varepsilon~ \Ha(\eta~|~\mu).
\end{equation}
\end{theo}

A precise version of this theorem is given in Section~\ref{sec-entrop-theo} (see Theorem~\ref{lem-Int-log-sob-talagrand}). For instance, (\ref{ref-Hep-D}) and (\ref{ref-Hep}) hold when the Hessian of the log density of $\Ka$ is uniformly bounded and its gradient is globally Lipschitz (cf. Condition  $A_2(\Ka)$ on page~\pageref{strong-convex-2}). Weaker conditions are discussed in Section~\ref{reg-cond-section-ab}.
The inequalities in \eqref{ref-Hep-D} provide $2$-Wasserstein stability bounds for Schrödinger bridges. When $b<1$, they yield a contraction property for the transport map $\eta\mapsto \eta\Ka$ in the $2$-Wasserstein distance. Similarly, when $\varepsilon<1$, the bounds in \eqref{ref-Hep} give a contraction estimate in relative entropy (see for instance (\ref{ref-contract-tt})).

\subsection{Sinkhorn contraction theorems}\label{sinkcontract-sec}

Without additional effort, Theorem~\ref{th1-intro} applied to the entropy formulae (\ref{bridge-form-v2}) and (\ref{bridge-form-2-V2}) already yields a variety of quantitative estimates.
To be more specific, the estimates (\ref{ref-Hep})  applied to (\ref{bridge-form-v2}) and (\ref{bridge-form-2-V2}) 
hold when the Hessians of the log density of $\Ka_{n}$ are uniformly bounded and their gradients is globally Lipschitz (cf. Section~\ref{general-sec}). 

A first key structural property of Sinkhorn log-densities is that {\it the oscillation of their gradients does not depend on the iteration index  $n\geq 0$} (cf. Lemma~\ref{nabla1-WN-lem}).

A second important feature is that their Hessians can be controlled in terms of the curvature of the potential functions $(U,V)$ (cf. Lemma~\ref{lem-appendix-Hess} and section~\ref{general-sec}).
Consequently, when these Hessians are uniformly bounded, applying \eqref{ref-Hep-D} yields the existence of constants
$a,b>0$ such that
\begin{equation}\label{first-was}
~\Da_{2}(\eta,\pi_{2n})^2\leq  a~
 \Ha(P_{\mu,\eta}~|~\Pa_{2n})\leq b
~\Da_{2}(\pi_{2n-1},\mu)^2.
\end{equation}
Moreover, if both $\mu$ and $\eta$ satisfy a quadratic transportation cost inequality, then  (\ref{ref-Hep}) applied to (\ref{bridge-form-v2}) and (\ref{bridge-form-2-V2})  guarantees the existence of  $\varepsilon>0$ such that
\begin{equation}\label{ref-HH-ep}
 \Ha(P_{\mu,\eta}~|~\Pa_{2n})
 \leq \varepsilon~ \Ha(\pi_{2n-1}~|~\mu)
\quad\mbox{\rm and}\quad
\Ha(P_{\mu,\eta}~|~\Pa_{2n+1})\leq 
\varepsilon~ \Ha(\pi_{2n}~|~\eta).
\end{equation}
Weaker conditions are discussed in Section~\ref{general-sec}.
Combining these inequalities with \eqref{ent-r1} yields the following quantitative convergence results.
\begin{theo}\label{th2-intro}
When the Markov transitions $\Ka_n$ are regular enough,  there 
exist constants
$a,b>0$ such that (\ref{first-was}) holds and for any $n\geq 0$,
$$
 \Da_{2}(\pi_{2n+1},\mu)^2\leq a~
 \Ha( P_{\mu,\eta}~|~\Pa_{2n+1})\leq b~\Da_2(\pi_{2n},\eta)^2.
$$
In addition, if
$\mu$ and $\eta$ satisfy a quadratic transportation cost inequality,
 there exists 
$\varepsilon>0$ such that for all $n\geq 0$,
\begin{equation}\label{ref-intro-ep}
 \Ha(P_{\mu,\eta}~|~\Pa_{n+2})\leq \left(1+  \varepsilon^{-1}
\right)^{-1}~ \Ha(P_{\mu,\eta}~|~\Pa_{n}).
\end{equation}
If moreover $C_0^-(U,V)$ holds with sufficiently large curvature parameters, then a sharper exponential rate holds:
     \begin{equation}\label{def-vareps-intro}
 \Ha(P_{\mu,\eta}~|~\Pa_{n+2})\leq
\left(1+\phi(\varepsilon)\right)^{-n}~  \Ha(P_{\mu,\eta}~|~\Pa_{0})
\quad\mbox{with}\quad
(1+\phi(\varepsilon))^{-2}<(1+\varepsilon^{-1})^{-1}.
\end{equation}
\end{theo}
A more precise statement and proof are given in Section~\ref{sinkhorn-sec}; see Lemma~\ref{theo-sinkhorn-1}, the contraction bounds \eqref{f1-c}–\eqref{f2-c}, and Theorem~\ref{theo-imp}.

\begin{theo}\label{theo-intro-lingauss}
Consider the
linear-Gaussian reference potential (\ref{def-W}). Then the estimates in (\ref{ref-intro-ep})  hold for even indices $n$ (respectively for odd indices $n$) whenever $C_U^-(V)$ (respectively $C_V^-(U)$) is satisfied. 
Moreover, if  $C_U(V)$ (respectively $C_V(U)$) holds for some $u_+,v_+>0$,
 the parameters $\varepsilon$ in (\ref{ref-intro-ep}) is given by
  \begin{equation}\label{def-vareps-intro-t}
\varepsilon:=\Vert \tau^{-1}\beta\Vert^2_2~\Vert u_+\Vert_2~\Vert v_+\Vert_2.
\end{equation} 
\end{theo}

A more detailed statement is provided in Section~\ref{sec-lingauss}; see (\ref{f1-c-uv}),(\ref{f1-c-uv-i}), Theorems~\ref{th-1-in} and~\ref{theo-thc-s}.

\begin{rmk}
Since
$$
   \Ha(P_{\mu,\eta}~|~\Pa_{2n+1})\leq    \Ha(P_{\mu,\eta}~|~\Pa_{2n})\leq  \Ha(P_{\mu,\eta}|~\Pa_{2n-1}),
   $$
    exponential decay for even or odd subsequences is equivalent.  
   In particular, for linear–Gaussian reference potentials, Theorem~\ref{theo-intro-lingauss} ensures exponential convergence of Sinkhorn bridges to the Schrödinger bridge as soon as either 
     $C_U^-(V)$  or $C_V^-(U)$ holds. For instance, when  $C_U(V)$ holds, then for all $n\geq 1$,
     \begin{equation}\label{ref-cv-i}
   \Ha(P_{\mu,\eta}~|~\Pa_{2n+1})\leq
\Ha(P_{\mu,\eta}~|~\Pa_{2n})\leq \left(1+ \varepsilon^{-1}
\right)^{-n}~ \Ha(P_{\mu,\eta}~|~\Pa_{0}).
\end{equation}

\end{rmk}

Combining Theorem~\ref{th1-intro} and Theorem~\ref{th2-intro} yields the following contraction estimates.. 

\begin{theo}\label{th3-intro}
Assume a linear–Gaussian reference potential \eqref{def-W} and that
  $C_0^-(U,V)$ holds.  Then, for any $n\geq 1$, 
the exponential decays
\begin{eqnarray*}
   \Ha(\eta~|~\pi_{2n})&\leq &\varepsilon
\left(1+\phi(\varepsilon)\right)^{-2(n-1)}~ \Ha(\eta~|~\pi_{0}),\\
\Da_{2}(\eta,\pi_{2n})
 &\leq &\varepsilon
\left(1+\phi(\varepsilon)\right)^{-(n-1)}~\Da_2(\eta,\pi_{0}),
 \end{eqnarray*}
  with $(\varepsilon,\phi(\varepsilon))$ as in (\ref{def-vareps-intro}).
  The same estimates hold upon replacing $(\eta,\pi_{2n})$ by $(\mu,\pi_{2n+1})$.
\end{theo}
A more precise version is given in Section~\ref{sec-lingauss} (Corollary~\ref{cor-ref-i}).

\begin{theo}\label{ineq-theo-impp-2-i-th}
Assume a linear–Gaussian reference potential and that  $C_1(U,V)$ holds. Then, there exist  $a>0$ and  $0\leq b<1$ such that
 for $n\geq p\geq 1$,
\begin{equation}\label{ineq-theo-impp-2-i}
\begin{array}{l}
 \displaystyle \Ha(P_{\mu,\eta}~|~\Pa_{2(n+1)})
 \leq  \left(1+   \varepsilon^{-1}~ \iota(\overline{\varpi})~\left(1+
a\, b^{\,p}\right)^{-1}\right)^{-(n-p)}  \Ha(P_{\mu,\eta}~|~\Pa_{2p-1}),
\end{array}\end{equation}
with  $\varepsilon$ as in (\ref{def-vareps-intro-t}) and the parameter $\iota(\overline{\varpi})$ defined by
$$
\iota(\overline{\varpi}):=
\frac{\Vert v_+\Vert_2}{\Vert v^{1/2}_+~ r_{\overline{\varpi}_0}~v^{1/2}_+\Vert_2}\vee
\frac{\Vert u_+\Vert_2}{\Vert u^{1/2}_+~  r_{\overline{\varpi}_1}~u^{1/2}_+\Vert_2}~\geq 1.
$$
\end{theo}

A more precise statement is provided in Section~\ref{scvx-sec} (see Corollary~\ref{cor-impp}).

When both matrices $u_-=0=v_-$ are null we have $\iota(\overline{\varpi})=1$. In this case,
 the estimates (\ref{ineq-theo-impp-2-i}) reduce to the estimate  (\ref{ref-cv-i}) for linear Gaussian models stated in Theorem~\ref{theo-intro-lingauss}.
 
 Consider the situation $\tau=t~\tau_0$ for some $\tau_0>0$ and some scalar $t>0$.
 As shown in Example~\ref{examp-ref-tau0}, whenever $u_->0$ or $v_->0$, for $t$ sufficiently small we have a sharper exponential convergence rate
 $$
( 1+\varepsilon^{-1}~\iota(\overline{\varpi}))^{-1}
< (1+\phi(\varepsilon))^{-2}<(1+\varepsilon^{-1})^{-1}.
 $$

\subsection{Entropic bridge maps}\label{sec-ent-maps}

This subsection investigates the regularity and contraction properties of {\it entropic bridge maps}, also known as barycentric projections of Schr\"odinger bridges. We first recall the interpretation of entropic bridges as regularized optimal transport problems and introduce notation emphasizing the role of the regularization parameter. We then derive a conditional covariance representation for the gradients of entropic maps, extending identities known in the linear--Gaussian case to general target measures. These covariance formulae provide quantitative estimates for entropic maps in terms of Riccati fixed points.

Schr\"odinger bridges and the Sinkhorn algorithm are widely used in practice to address high-dimensional optimal transport problems; see the pioneering work of Cuturi~\cite{cuturi} and the more recent contributions~\cite{chewi,genevay-cuturi,genevay-cuturi-2,pooladian,pooladian-2,seguy}.
More specifically, consider the Schr\"odinger bridge $P_{\mu,\eta}$ defined in \eqref{ob-eq} with target measures $(\mu,\eta)=(\lambda_U,\lambda_V)$ and reference measure $P:=\mu\times K_W$, where the potential $W$ is of the form \eqref{def-W} with $\tau=t\,\tau_0$ for some positive definite matrix $\tau_0>0$ and scalar $t>0$.

\begin{defi}
To emphasize the role of the regularization parameter $t$, we denote by
\begin{equation}\label{sbr-ref-tt}
P_{\mu,\eta}^t:=\mu\times K_{\WW^t_{\mu,\eta}},
\qquad
P_{\eta,\mu}^{t,\flat}:=\eta\times K_{\WW^{t,\flat}_{\eta,\mu}},
\end{equation}
the Schr\"odinger bridges \eqref{sbr-ref} and the associated transition potentials
$(\WW^t_{\mu,\eta},\WW^{t,\flat}_{\eta,\mu})$ corresponding to the choice $\tau=t\,\tau_0$.
We denote by $(\Sigma^t_{\mu,\eta},\Sigma^{t,\flat}_{\eta,\mu})$ the conditional covariance matrices \eqref{cov-mat-W}, and by $(\ZZ^t_{\mu,\eta},\ZZ^{t,\flat}_{\eta,\mu})$ the entropic maps \eqref{e-maps} associated with the transitions
$(K_{\WW^t_{\mu,\eta}},K_{\WW^{t,\flat}_{\eta,\mu}})$.
\end{defi}

Note that
 $$
 W(x,y)=\frac{1}{t}~c(x,y)+\frac{1}{2}\log{(\mbox{\rm det}(2t\pi\tau_0 I))},
$$
where
$$
c(x,y):=\frac{1}{2}~(y-(\alpha+\beta x))^{\prime}\tau_0^{-1}(y-(\alpha+\beta x)).
 $$
  In this setting, for any $Q\in \Pi(\mu,\eta)$ we have the entropic representation
 \begin{equation}\label{entoo}
\Ha(Q~|~P)=\frac{1}{t}~\int~c(x,y)~Q(d(x,y))+\Ha(Q~|~\mu\otimes\eta)-\eta(V)+\frac{1}{2}\log{(\mbox{\rm det}(2t\pi\tau_0 I))}.
\end{equation}

The parameter $t$ thus plays the role of an entropic regularization parameter (see Section~\ref{reg-transp-sec}), while $c(x,y)$ can be interpreted as the underlying transportation cost.

In this interpretation, as $t\to0$ the Schr\"odinger bridge $P^t_{\mu,\eta}$ converges, in a suitable sense (see for instance \eqref{MK} and the Gaussian case discussed in Example~\ref{examp-otgauss}), to the optimal transport bridge $P^0_{\mu,\eta}$ solving the Monge--Kantorovich problem
$$
P_{\mu,\eta}^0
:=
\argmin_{Q\in\Pi(\mu,\eta)}\int c(x,y)\,Q(d(x,y)).
$$
In this section, we establish new quantitative estimates for general cost functions of the form \eqref{def-W}, including refined versions of results by Chewi and Pooladian~\cite{chewi} and extensions of Caffarelli’s contraction theorem~\cite{caffarelli}.

\subsubsection*{Covariance formulae}

Consider the linear--Gaussian reference transitions $K_W$ associated with a potential $W$ of the form \eqref{def-W}, with parameters $(\alpha,\beta,\tau)$, and recall the notation $\cchi:=\tau^{-1}\beta$.

Following~\cite{pooladian,pooladian-2,seguy}, the barycentric projections (or {\it entropic maps}) of the Schr\"odinger bridges $K_{\WW_{\mu,\eta}}$ and $K_{\WW^{\flat}_{\eta,\mu}}$ are defined via the identity map $I(y):=y$ by the conditional expectations

\begin{equation}\label{bary-ref}
\begin{array}{rcl}
K_{\WW_{\mu,\eta}}(I)(x)&:=&\displaystyle\int~K_{\WW_{\mu,\eta}}(x,dy)~y=\EE(\ZZ_{\mu,\eta}(x)),
\\
K_{\WW^{\flat}_{\eta,\mu}}(I)(x)&:=&\displaystyle\int~K_{\WW^{\flat}_{\eta,\mu}}(x,dy)~y
=\EE(\ZZ^{\flat}_{\eta,\mu}(x)).
\end{array}
\end{equation}
In the above display, $(\ZZ_{\mu,\eta},\ZZ^{\flat}_{\eta,\mu})$   stands for entropic maps 
defined in (\ref{e-maps}). These conditional expectations provide approximations of optimal transport maps; see
\cite{dp-ent-26,ferradan,reich,seguy-15,seguy}.
 \begin{examp}\label{examp-bary}
Consider the linear--Gaussian model defined in \eqref{def-W} and \eqref{def-U-V}, with covariance matrices $u,v>0$ and parameters $(\alpha,\beta,\tau)$.
By \eqref{ref-WW-Gauss}, the entropic maps $(\ZZ_{\mu,\eta},\ZZ^{\flat}_{\eta,\mu})$ admit the representations
\begin{eqnarray}
\ZZ_{\mu,\eta}(x)&=&\overline{m}+\sigma_{\mu,\eta}~\cchi~ (x-m)+\sigma_{\mu,\eta}^{1/2}~G,
\nonumber\\
\ZZ^{\flat}_{\eta,\mu}(x)&=&m+\sigma_{\eta,\mu}^{\flat}~\cchi^{\prime}~ (x-\overline{m})+(\sigma^{\flat}_{\eta,\mu})^{1/2}~G.
\label{ref-WW-Gauss-eom}
\end{eqnarray}
where $G$ is a centered Gaussian random variable with unit covariance.

In terms of the Riccati fixed-point matrices $(r_{\overline{\varpi}_0},r_{\overline{\varpi}_1})$ defined in \eqref{ricc-maps-def}, we obtain
\begin{eqnarray}
\nabla K_{\WW_{\mu,\eta}}(I)(x)&=&\nabla \ZZ_{\mu,\eta}(x)=
\cchi^{\prime}~v^{1/2} ~r_{\overline{\varpi}_0}~v^{1/2},
\nonumber\\
\nabla K_{\WW^{\flat}_{\eta,\mu}}(I)(x)&=&\nabla \ZZ_{\mu,\eta}^{\flat}(x)=
\cchi~u^{1/2} ~r_{\overline{\varpi}_1}~u^{1/2}.\label{cov-grad-gauss}
\end{eqnarray}
\end{examp}

The following theorem shows that the covariance-based gradient representation \eqref{cov-grad-gauss}, which holds exactly in the Gaussian case, extends to general target measures.

\begin{theo}\label{theo-baryp}
We have the conditional covariance identities
\begin{equation}\label{conditional-cov}
\cchi^{\prime}~\Sigma_{\mu,\eta}~\cchi=\nabla K_{\WW_{\mu,\eta}}(I)~\cchi\quad\mbox{and}\quad
\cchi~\Sigma^{\flat}_{\eta,\mu}~\cchi^{\prime}=\nabla K_{\WW^{\flat}_{\eta,\mu}}(I)~\cchi^{\prime}.
\end{equation}
\begin{itemize}
\item 
Assume there exists positive definite matrices
$u_->0$ and $v_+>0$ such that
\begin{equation}\label{bary-ref-vu}
 \nabla^2 U\leq u^{-1}_-\quad\mbox{and}\quad
v_+^{-1}\leq \nabla^2 V.
\end{equation}
Then
\begin{equation}\label{ch-vu}
\nabla K_{\WW_{\mu,\eta}}(I)~\cchi\leq
\cchi^{\prime}~v_+^{1/2} ~r_{\overline{\varpi}_0}~v_+^{1/2}~\cchi~~ \mbox{and}~~
 \nabla K_{\WW^{\flat}_{\eta,\mu}}(I)~\cchi^{\prime}\geq \cchi~u_-^{1/2}~r_{\varpi_1}~u_-^{1/2}~\cchi^{\prime}.
\end{equation}
\item Assume there exists some positive definite matrices
$u_+>0$ and $v_->0$ such that
\begin{equation}\label{bary-ref-uv}
u_+^{-1}\leq \nabla^2 U\quad\mbox{and}\quad
 \nabla^2 V\leq v^{-1}_-.
\end{equation}
Then
\begin{equation}\label{ch-uv}
 \nabla K_{\WW^{\flat}_{\eta,\mu}}(I)~\cchi^{\prime}\leq
\cchi~u_+^{1/2}~r_{\overline{\varpi}_1}~u_+^{1/2}~\cchi^{\prime} 
~~ \mbox{and}~~
\nabla K_{\WW_{\mu,\eta}}(I)~\cchi\geq \cchi^{\prime}~ v_-^{1/2} ~r_{\varpi_0}~v_-^{1/2}~\cchi.
\end{equation}

\end{itemize}
\end{theo}

he proof of the covariance identities \eqref{conditional-cov} relies on elementary differential calculus and is therefore deferred to Appendix~\ref{sec-grad-hess-a}; see in particular \eqref{finex} and Example~\ref{ex-gauss-UVn}.

Condition \eqref{bary-ref-vu} includes, for instance, the case where $\nabla^2 U\ge -a_u I$ for some $a_u\ge0$, while \eqref{bary-ref-uv} covers the case $\nabla^2 V\ge -a_v I$ for some $a_v\ge0$.

Theorem~\ref{theo-baryp} is in fact a special case of the more general Theorem~\ref{last-est-th}, which applies to potential functions $(U,V)$ satisfying condition $C_1(U,V)$.
In the linear--Gaussian setting \eqref{def-W}–\eqref{def-U-V}, the entropic maps are affine, and the inequalities in Theorem~\ref{theo-baryp} reduce to equalities; see \eqref{cov-grad-gauss}.

\subsection{Entropic regularizations}

This subsection studies the small-$t$ regime of Schr\"odinger bridges and entropic maps, interpreting the parameter $t$ as an entropic regularization strength. We introduce a rescaled parametrization of the Riccati matrices adapted to the limit $t\to 0$ and derive asymptotic estimates for the gradients and conditional covariances of entropic bridge maps. These results quantify the convergence of entropic maps toward optimal transport maps and yield refined contraction bounds, extending Caffarelli-type estimates and recent results of Chewi and Pooladian to general cost functions and non-commuting covariance structures.

Set $\cchi_0:=\tau_0^{-1}\beta$.  
To emphasize the dependence on $t$ we write $({\varpi}_i(t),\overline{\varpi}_i(t))_{i=0,1}$ the  matrices defined in (\ref{def-over-varpi}) associated with the choice $\tau=t~\tau_0$.
With this notation, one readily checks that
\begin{equation}\label{def-over-varpi-check-p}
{\varpi}_i(t)=t^2~ \widecheck{\varpi}_i 
\quad \mbox{\rm and}\quad
  \overline{\varpi}_i=t^2~ \widehat{\varpi}_i,   
\end{equation}
where
\begin{equation}\label{def-over-varpi-check}
\begin{array}{rcl}
 \widecheck{\varpi}_0^{-1}&:=&v^{1/2}_-~\left(\cchi_0~u_+~\cchi_0^{\prime}\right)~v^{1/2}_-\\
 \widecheck{\varpi}_1^{-1}&:=&u_-^{1/2}~\left(\cchi_0^{\prime}~v_+~\cchi_0\right)~u_-^{1/2}
\end{array}
\quad\&\quad
\begin{array}{rcl}
\widehat{\varpi}_0^{-1}&:=&v^{1/2}_+~(\cchi_0~u_-~\cchi_0^{\prime})~v^{1/2}_+,\\
\widehat{\varpi}_1^{-1}&:=&u_+^{1/2}~(\cchi_0^{\prime}~v_-~\cchi_0)~u_+^{1/2}.
\end{array}
\end{equation}

With these definitions, the conditional covariance formulae of Theorem~\ref{theo-baryp} yield the following asymptotic estimates in the small-$t$ regime.

\begin{cor}\label{cor-entropic-map-t}
Assume that $\tau=t,\tau_0$ with $\tau_0>0$, and that condition \eqref{bary-ref-vu} holds for some $u_->0$ and $v_+>0$. Then, for every $x\in\RR^d$,
\begin{eqnarray*}
\limsup_{t\rightarrow 0}\ell_{\tiny max}\left(\nabla K_{\WW^t_{\mu,\eta}}(I)(x)~\cchi_0\right)&\leq& \ell_{\tiny min}\left(\cchi^{\prime}_0~\left((\cchi_0~u_-~\cchi_0^{\prime})^{-1}\sharp~v_+\right)~\cchi_0\right),
\\
\liminf_{t\rightarrow 0}\ell_{\tiny min}\left(\nabla K_{\WW^{t,\flat}_{\eta,\mu}}(I)(x)~\cchi_0^{\prime}\right)&\geq& \ell_{\tiny max}\left(\cchi_0~\left((\cchi_0^{\prime}~v_+~\cchi_0)^{-1}\sharp~u_-\right)~\cchi_0^{\prime}\right),
\end{eqnarray*}
where $u\sharp v$ denotes the geometric mean of two positive definite matrices (see \eqref{sym-sharp}).
\end{cor}
\proof
Since $\widehat{\varpi}_0^{-1}>0$ and $\widecheck{\varpi}_1^{-1}>0$, equation~\eqref{def-fix-ricc-1} implies
 $$
 \frac{r_{\overline{\varpi}_0}}{t}~\longrightarrow_{t\rightarrow 0}~ 
\widehat{\varpi}_0^{1/2}\quad\mbox{\rm and}\quad
 \frac{r_{{\varpi}_1}}{t}~\longrightarrow_{t\rightarrow 0}~ 
\widecheck{\varpi}_1^{1/2}.
 $$
Combining these limits with \eqref{ch-vu} yields
$$
\nabla K_{\WW^t_{\mu,\eta}}(I)(x)~\cchi_0\leq
\cchi^{\prime}_0~v_+^{1/2} ~\frac{r_{\overline{\varpi}_0}}{t}~v_+^{1/2}~\cchi_0,$$
and
$$
 \nabla K_{\WW^{t,\flat}_{\eta,\mu}}(I)(x)~\cchi^{\prime}_0\geq \cchi_0~u_-^{1/2}~\frac{r_{\varpi_1}}{t}~u_-^{1/2}~\cchi^{\prime}_0.
$$
The stated bounds follow by taking limits and applying elementary spectral inequalities.
\cqfd
\begin{rmk}\label{rmk-compcp}
Under condition \eqref{bary-ref-vu}, inequality \eqref{ch-vu} further implies
$$
\frac{\Sigma^t_{\mu,\eta}(x)}{t}~\leq~v_+^{1/2} ~\frac{r_{\overline{\varpi}_0}}{t}~v_+^{1/2}\longrightarrow_{t\rightarrow 0}~ v_+^{1/2} ~\widehat{\varpi}_0^{1/2}~v_+^{1/2}.
$$
Consequently the covariance matrix $\Sigma^t_{\mu,\eta}(x)$ of the Schr\" odinger bridge transition $K_{\WW^t_{\mu,\eta}}(x,dy)$ converges  towards the null matrix (uniformly w.r.t. the variable $x$). Whenever (\ref{ch-vu}) holds, we recover the fact that the optimal $2$-Wasserstein coupling in  (\ref{MK-tr}) is almost surely deterministic:
$$
T_{\mu,\eta}(x)=\ZZ^0_{\mu,\eta}(x):=\lim_{t\rightarrow 0}\EE(\ZZ^t_{\mu,\eta}(x))\quad\mbox{and}\quad
\nabla \ZZ^0_{\mu,\eta}(x)~\cchi_0\leq \cchi^{\prime}_0~\left((\cchi_0~u_-~\cchi_0^{\prime})^{-1}\sharp~v_+\right)~\cchi_0.
$$ 
In addition, when $\tau_0=I=\beta$ the above estimate resumes to 
$$
\nabla \ZZ^0_{\mu,\eta}(x)\leq u_-^{-1}~\sharp~ v_+.
$$
If, moreover $u_-$ and $v_+$ commute, this coincides with the estimate of Chewi and Pooladian~\cite{chewi}:
$$
\nabla \ZZ^0_{\mu,\eta}(x)\leq (u_-^{-1}~v_+)^{1/2}=u_-^{-1/2}~v_+^{1/2}.
$$
 \end{rmk}

Following word-for-word the proof of Corollary~\ref{cor-entropic-map-t} we check the following result.
\begin{cor}\label{cor-entropic-map-t-2}
Consider the situation $\tau=t~\tau_0$ for some positive definite matrix  $\tau_0>0$ and some scalar $t>0$.
Assume that (\ref{bary-ref-uv}) holds for some $u_+>0$ and $v_->0$. In this situation, for any $x\in\RR^d$ we have
\begin{eqnarray*}
\liminf_{t\rightarrow 0}\ell_{\tiny min}\left(\nabla K_{\WW^t_{\mu,\eta}}(I)~\cchi_0\right)&\geq&  \ell_{\tiny max}\left(\cchi^{\prime}_0~\left((\cchi_0~u_+~\cchi_0^{\prime})^{-1}~\sharp~v_-\right) ~\cchi_0\right),
\\
\limsup_{t\rightarrow 0}\ell_{\tiny max}\left(
 \nabla K_{\WW^{t,\flat}_{\eta,\mu}}(I)~\cchi_0^{\prime}\right)&\leq&
\ell_{\tiny min}\left(
\cchi_0~\left((\cchi_0^{\prime}~v_-~\cchi_0)^{-1}~\sharp~ u_+\right)~\cchi_0^{\prime} 
\right).
\end{eqnarray*}

\end{cor}

For the linear-Gaussian model discussed in (\ref{def-W}) and (\ref{def-U-V}) the entropic maps are linear and the limiting properties stated in Corollary~\ref{cor-entropic-map-t} and Corollary~\ref{cor-entropic-map-t-2} can be turned into the convergence of entropy maps.   
\begin{examp}\label{examp-gauss-sh}
Consider the linear-Gaussian model discussed in Example~\ref{examp-bary}.
Also assume that $\tau=t~\tau_0$ for some positive definite matrix  $\tau_0>0$ and some scalar $t>0$.
In this context, we have $(\overline{\varpi}_0(t),
\overline{\varpi}_1(t))=({\varpi}_0(t),{\varpi}_1(t))$ as well as 
$$
(\widecheck{\varpi}_0,\widecheck{\varpi}_1)=
(\widehat{\varpi}_0,\widehat{\varpi}_1)=\left(v^{-1/2}~\left(\cchi_0~u~\cchi_0^{\prime}\right)^{-1}~v^{-1/2},u^{-1/2}~\left(\cchi_0^{\prime}~v~\cchi_0\right)^{-1}~u^{-1/2}\right).
$$
Denote by $(\ZZ^t_{\mu,\eta},\ZZ^{t,\flat}_{\eta,\mu})$ the entropic optimal maps (\ref{ref-WW-Gauss-eom}) associated with the covariance matrices $ \tau=t~\tau_0$.
In this notation,  for any $x\in\RR^d$ we have
$$
\lim_{t\rightarrow 0}\nabla K_{\WW^t_{\mu,\eta}}(I)(x)=\lim_{t\rightarrow 0}\nabla \ZZ^t_{\mu,\eta}(x)
=\cchi_0^{\prime}~\left(\left(\cchi_0~u~\cchi_0^{\prime}\right)^{-1}~\sharp~ v\right),
$$
$$
\lim_{t\rightarrow 0}\nabla K_{\WW^{t,\flat}_{\eta,\mu}}(I)(x)=
\lim_{t\rightarrow 0}\nabla \ZZ^{t,\flat}_{\eta,\mu}(x)
=\cchi_0~
\left(\left(\cchi_0^{\prime}~v~\cchi_0\right)^{-1}~\sharp~u\right).
$$
Quantitative estimates can be obtained following the proof of Corollary 3.10 in~\cite{adm-24}.
\end{examp}

 \subsection*{Notes and references}\label{ref-sec-i}

Theorem~\ref{theo-baryp} and Corollary~\ref{cor-entropic-map-t} in
Section~\ref{sec-ent-maps} extend and refine Caffarelli’s contraction
theorem~\cite{caffarelli}, as well as results of Chewi and Pooladian, to
{\it general cost functions of the form~\eqref{def-W}} and to covariance
structures that are not necessarily commuting (see Theorems~1 and~4
in~\cite{chewi} and Remark~\ref{rmk-compcp} in the present article).
As shown in Example~\ref{examp-gauss-sh}, the resulting bounds are sharp
for broad classes of linear–Gaussian models.

The proof of Theorem~\ref{th1-intro} combines optimal entropic cost
couplings with Talagrand’s quadratic transportation inequality and the
$LS(\rho)$ logarithmic Sobolev inequality. Closely related arguments are
used in~\cite{chewi-phd} and~\cite{lee} to establish contraction
properties of proximal samplers. This line of reasoning provides the
main motivation for the contraction framework for Schr\"odinger and
Sinkhorn bridges developed in this work, which is based on
transportation–cost inequalities.

Related quantitative stability estimates of the form~\eqref{ref-intro-ep},
with different constants $\varepsilon>0$, are also derived
in~\cite{chiarini} for certain classes of entropic optimal transport
problems, under the assumption that the (unknown) Schr\"odinger bridge
transitions admit uniformly semi-concave log-densities. In that setting,
the reference measure is of the form $P := \mu \otimes \Ka$, where
$\Ka = K_W$ is associated with a potential
\begin{equation}\label{quadratic-case-ref}
  W(x,y) = c + d(x,y)^2,
\end{equation}
for some distance $d$ on $\RR^d$ and a normalizing constant $c$.

Note that the  linear-Gaussian models (\ref{def-W}) with $(\alpha,\beta)=(0,I)$ corresponds to the situation where the  distance is a quadratic form associated with some positive matrix $\tau^{-1}>0$. 
The case $\tau=t\tau_0$ for some positive definite matrix $\tau_0>0$ and some scalar $t>0$ can be interpreted as an entropic regularization of the Monge-Kantorovich problem (cf. Section~\ref{reg-transp-sec} in the present article). In this context the parameter $\varepsilon$ defined in
(\ref{def-vareps-intro-t}) reduces to
\begin{equation}\label{ref-imp-eps}
\varepsilon=\frac{1}{t^2}~\Vert \tau_0^{-1}\beta\Vert^2_2~\Vert u_+\Vert_2~\Vert v_+\Vert_2.
\end{equation}

\subsection{Some illustrations}\label{sec-illu}

\subsubsection{Regularized transport}\label{reg-transp-sec}
Consider the model $(\mu,\eta,\Ka)=(\lambda_U,\lambda_V,K_W)$ defined in (\ref{def-triple}) with the  linear Gaussian potential $W$ given by (\ref{ref-KW}) with $(\alpha,\beta)=(0,I)$ and $\tau=t~I$ for some parameter $t>0$. In this context, the parameter defined in (\ref{def-vareps-intro-t}) reduces to
$$
\varepsilon=\frac{1}{t^2}~\Vert u_+\Vert_2~\Vert v_+\Vert_2.
$$
In terms of the bridges (\ref{sbr-ref-tt}) associated with the parameter $\tau=t I$  by (\ref{entoo}) we have
\begin{equation}\label{ref-tcost}
P^t_{\mu,\eta}=
\argmin_{Q\,\in\, \Pi(\mu,\eta)}\left(\frac{1}{2t}\int~\Vert x-y\Vert_2^2~Q(d(x,y))+\Ha(Q~|~\mu\otimes\eta)\right).
\end{equation}
Following the discussion provided in section~\ref{sec-ent-maps}, the parameter $t$ can be seen as a regularization parameter of the Monge-Kantorovich problem. When $t\rightarrow 0$ the variational problem amounts to find an optimal coupling  that realizes the infimum
$$
P_{\mu,\eta}^0:=
\argmin_{Q\,\in\, \Pi(\mu,\eta)}\int~\Vert x-y\Vert_2^2~Q(d(x,y)).
$$
In terms of the $2$-Wasserstein distance, the optimal transport plan $T_{\mu,\eta}$ solving the Monge-Kantorovich problem is given by
\begin{equation}\label{MK}
T_{\mu,\eta}:=\argmin_{T:~\mu\circ T^{-1}=\eta}\Da_2(\mu,\mu\circ T^{-1})\quad\mbox{\rm with}\quad
(\mu\circ T^{-1})(f):=\mu( f\circ T).
\end{equation}
This shows that the optimal $2$-Wasserstein coupling is deterministic and given by
\begin{equation}\label{MK-tr}
\begin{array}{l}
P_{\mu,\eta}^0(d(x,y))=
\mu(dx)~ \delta_{T_{\mu,\eta}(x)}(dy),\\
\\
\mbox{\rm and}\quad
K_{\WW^0_{\mu,\eta}}(x,dy)= \delta_{\ZZ^0_{\mu,\eta}(x)}(dy)
\quad\mbox{\rm with}\quad \ZZ^0_{\mu,\eta}(x)=T_{\mu,\eta}(x).
\end{array}\end{equation}
\begin{examp}\label{examp-otgauss}
Consider the linear-Gaussian model discussed in (\ref{def-W}) and (\ref{def-U-V})  with some mean values $(m,\overline{m})$, some covariance matrices $u,v>0$ and the parameters $\alpha=0$ and $(\beta,\tau)$=(I,tI). In this context, using (\ref{ref-WW-Gauss-eom}) and following the proof of the estimates provided in Example~\ref{examp-gauss-sh} we check that
$$
\lim_{t\rightarrow 0}K_{\WW^t_{\mu,\eta}}(I)(x)=\ZZ^0_{\mu,\eta}(x)=
T_{\mu,\eta}(x)=\overline{m}+(u^{-1}~\sharp~v)~(x-m).
$$
For a more detailed discussion on the convergence of entropic maps general classes of linear Gaussian models we refer to~\cite{adm-24}, see also Theorem 2.1 in~\cite{gelbrich} for a closed form expression of $2$-Wasserstein distance and optimal transport plans.
\end{examp}
In the reverse angle when $t\rightarrow\infty$ the limiting 
variational problem (\ref{ref-tcost}) reduces to the tensor product measure
$$
\mu\otimes\eta=\argmin_{Q\,\in\, \Pi(\mu,\eta)}\Ha(Q~|~\mu\otimes\eta).
$$
In this context, using (\ref{ref-cv-i}) we check the uniform convergence
$$
\sup_{n\geq 1}\Ha(P_{\mu,\eta}~|~\Pa_{n})~\longrightarrow_{t\rightarrow\infty}~0.
$$
Whenever $C_0(U,V)$ holds, Theorem~\ref{th3-intro} also yields the uniform convergence
$$
\sup_{n\geq 1}\left(\Da_{2}(\eta,\pi_{2n})\vee \Da_{2}(\mu,\pi_{2n+1})\right)
~\longrightarrow_{t\rightarrow\infty}~0.
$$
\subsubsection{Proximal samplers}\label{prox-s}

When $n=0$ the distribution  $\Pa_0=P=\mu\times \Ka_0$ with $(\mu,\Ka_0)=(\lambda_U,K_W)$ is  the reference measure the Schr\" odinger bridge (\ref{ob-eq}) associated with the parameters (\ref{def-triple}). We further assume that $W$ is given by (\ref{def-W}) with $(\alpha,\beta)=(0,I)$ and $\tau=t I$. Recall that $\Vert z\Vert_2$ is the Euclidian norm of a state $z\in\RR^d$ and $(\Ka_0,\Ka_1)$ are the Sinkhorn transitions defined in (\ref{sinhorn-entropy-form-Sch}). In these settings the  two-blocks Gibbs sampler with the forward/backward transition $S_{1}:=\Ka_{0}\Ka_{1}$ and target product measure 
\begin{equation}\label{proximal-ex}
\Pa_0(d(x,y))\propto
\exp{\left(-U(x)-\frac{1}{2t}~\Vert x-y\Vert_2^2\right)}~\lambda(dx)\lambda(dy).
\end{equation} coincides with the proximal sampler introduced by Mou, Flammarion,  Wainwright and  Bartlett in~\cite{mou} and further developed in~\cite{chewi-chen,chewi-phd,jiaming,jiaojiao,lee}, see also the recent article~\cite{guan} for an extended version of the proximal sampler in Riemannian manifolds. These samplers 
can also be interpreted as sampling technique of $\lambda_U$ based on data augmentation or auxiliary variables~\cite{damlen1999gibbs,rendell2020global,vono2020asymptotically,vono2019split}.

In this context, the parameter $t$ is sometimes interpreted as the {\it step size} of the sampler and the backward sample is sometimes called a {\it restricted Gaussian oracle}
 (see for instance Section 4.2.2 in~\cite{chewi-phd}). The terminology "oracle" introduced in~\cite{lee,mou} comes from the fact that we generally cannot sample from the backward/dual transition $\Ka_1=K_{W^{\flat}_1}$. An oracle outputs samples from the backward transition "restricted" to some potential function $U$.  Sinkhorn iterates as well as the Gibbs sampler discussed above are  idealized samplers.   Metropolis Hasting within Gibbs type samplers can be used to approximate these samples, see for instance the adjusted proximal Metropolis type algorithm developed in~\cite{mou}.

\subsubsection{Diffusion generative models}\label{dgm-sec}

Forward-backward semigroups of the form discussed in Section~\ref{prox-s} also arise in the design and the analysis of  diffusion-based generative models~\cite{ho,sohl,song,song2020score,song2020score-arxiv}. As in (\ref{proximal-ex}), the forward step represented by $\Ka_0$ on some time horizon gradually transforms an unknown complex data distribution $\mu$ into a simpler distribution $\mu\Ka_0$ by injecting noise into data. The backward process represented by $\Ka_1$ reverses the noising process running backward from the final time horizon and by Gibbs' principle eventually generates samples from the initial data distribution. This backward Gibbs  transition is sometimes called the denoising transition step. 
Note that the backward transition $\Ka_{1}$ can alternatively be defined by the conjugate (a.k.a. dual) formula 
\begin{equation}
(\mu\Ka_0)\times \Ka_{1}=(\mu\times \Ka_{0})^{\flat}.\label{ref-intro-fb}
\end{equation}

We often start from a forward process represented by a discrete or continuous time stochastic flow $X_{s,t}(x)$
running from $s$ to $t\geq s$ and starting as $X_{s,s}(x)=x$ at time $t=s$. 
 At every time $t\geq 0$, the distributions of the internal random states of the forward process starting from $\mu$ are defined by the transport equations
 $$
 \nu_t=\mu P_{0,t}=\nu_sP_{s,t}\quad
\mbox{\rm with}\quad P_{s,t}(x,dy):=\PP(X_{s,t}(x)\in dy).
 $$
The backward flow $X_{t,s}(y)$ from $t$ to $s\leq t$ starting as $X_{t,t}(y)=y$ at time $s=t$
has a backward transition semigroup~\cite{anderson} defined by the conjugate formulae
$$
P_{t,s}(y,dx):=\PP(X_{t,s}(y)\in dx) \Longleftrightarrow
(\nu_s P_{s,t})\times P_{t,s}=(\nu_s\times P_{s,t})^{\flat}.
$$
Applying the above to $(s,t)=(0,t)$ we recover the conjugate formula (\ref{ref-intro-fb}) with
$$
\Ka_{0}(x,dy):=\PP(X_{0,t}(x)\in dy)\quad
\mbox{\rm with}\quad
\Ka_1(y,dx)=\PP(X_{t,0}(y)\in dx).
$$
For instance diffusion flows associated with some drift function $b_t(x)$ a Brownian motion $B_t$ and a diffusion matrix $\Sigma_t$ with appropriate dimensions  are given by
$$
dX_{s,t}(x)=b_t(X_{s,t}(x))~dt+\Sigma_t^{1/2}~dB_t.
$$
The null drift case $b_t=0$ and $\Sigma_t=I$ yields the transition discussed in (\ref{proximal-ex}); that is we have
\begin{equation}\label{r-HSG}
 K_0(x,dy)~\propto~
\exp{\left(-\frac{1}{2t}~\Vert x-y\Vert_2^2\right)}~\lambda(dy).
\end{equation}
A linear drift  $b_t(x)=Ax$ associated with some conformal matrix $A$ also yields the transition $\Ka_0=K_W$  discussed in (\ref{def-W}) with the potential function 
\begin{equation}\label{r-OU}
W(x,y)=-\log{g_{\tau}(y-\beta x)}\quad \mbox{\rm with}\quad
\beta=e^{tA}\quad \mbox{\rm and}\quad
\tau=\int_0^t e^{sA}~\Sigma~ e^{sA^{\prime}}~ds.
\end{equation}

Consider some positive  functions $(\beta_t,\tau_t)$ such that $(\beta_0,\tau_0)=(1,0)$ and for $t\in [0,1]$ set
\begin{equation}\label{ref-cnf}
b_t(x)=\partial_{t}\log{\beta_t}~x\quad\mbox{\rm and}\quad
\Sigma_t=\sigma_t~I\quad\mbox{\rm with}\quad \sigma_t=\beta^2_t~\partial_t\left(\frac{\tau_t}{\beta_t^2}\right).
\end{equation}
In this case the diffusion flow on the unit time interval $t\in [0,1]$ is given by
$$
X_{0,t}(x)=\beta_t~x+\beta_t~\int_0^t~
~\sqrt{\partial_s\left(\frac{\tau_s}{\beta_s^2}\right)}~dB_s\stackrel{\tiny law}{=}
\beta_t~x+\tau_t^{1/2}~B_1.
$$
This situation coincides with the linear-Gaussian transition  discussed in (\ref{def-W}) with the parameter $\alpha=0$ and $(\beta,\tau)$
replaced by $(\beta_t,\tau_t)$.  These diffusions are sometimes called diffusion flow matching~\cite{albergo,lipman,liu}, see also the more recent articles~\cite{bortoli,silveri} and references therein.

In all the  cases discussed above, whenever $\mu$ is Gaussian
the state distributions $\nu_t$ at every time horizon  as well as
 the backward transition  $\Ka_1$ are easy to sample and admits a closed form Gaussian formulation~\cite{adm-24}.

More generally, denoting $\Va_s$ the  potential function of the state distribution $\nu_s=\lambda_{\Va_s}$ when $\Sigma_s=I$ the time reversal process  is itself a Markov diffusion~\cite{anderson} defined backward in time by the diffusion
\begin{equation}\label{reverse}
-dX_{t,s}(x)=\left(-b_s(X_{t,s}(x))- \nabla \Va_s(X_{t,s}(x))\right) ds+~dB_s.
\end{equation}

In machine learning literature, the  two-blocks Gibbs sampler with transition $S_{1}:=\Ka_{0}\Ka_{1}$  is also called a {\it denoising diffusion model}~\cite{ho,bortoli-heng,shi,song2020score,song2020score-arxiv}.

 \begin{theo}\label{theo-S1n}
Let $S_1^n=S^{n-1}_1S_1$ be the semigroup of the   forward/backward Gibbs  sampler with transition $S_{1}:=\Ka_{0}\Ka_{1}$.
 Consider the
linear-Gaussian model (\ref{def-W}) associated with some parameters $(\alpha,\beta,\tau)$ and set
  \begin{equation}\label{def-delta-u-tau}
  \begin{array}{rcl}
   a_{u_+}(\beta,\tau)&:=&\Vert\tau^{-1}\beta\Vert_2^2~\Vert\tau\Vert_2~\Vert u_+\Vert_2,\\
   &&\\ 
 b_{u_+}(\beta,\tau)&:=&\Vert\tau^{-1}\beta\Vert_2^2~\Vert\tau\Vert_2~\Vert (u_+^{-1}+
     \beta^{\prime}\tau^{-1}\beta)^{-1}\Vert_2\leq a_{u_+}(\beta,\tau).
\end{array}   
 \end{equation}
Assume that $ b_{u_+}(\beta,\tau)<1$. In this situation,
for any $n\geq 0$ and any probability measure $\nu$ we have the exponential decay estimates
\begin{eqnarray*}
\Da_{2}(\nu S_{1}^n,\mu)\leq      b_{u_+}(\beta,\tau)^n~\Da_2(\nu ,\mu)
\quad\mbox{and}\quad
\Ha(\nu S_{1}^{n}~|~\mu )\leq a_{u_+}(\beta,\tau)~
 b_{u_+}(\beta,\tau)^{2n}~
\Ha(\nu~|~\mu).
\end{eqnarray*}
\end{theo}

A more precise statement is provided in Section~\ref{proxi-sec} (see Theorem~\ref{theo-proxi}). Note that condition $ b_{u_+}(\beta,\tau)<1$ is met for the Heat kernel (\ref{r-HSG})  as well as in linear-Gaussian model (\ref{r-OU})  when $A$ is an Hurwitz matrix as soon as the time horizon $t$ is sufficiently large.

Theorem~\ref{th3-intro} also ensures that Sinkhorn distributions $\pi_{2n+1}$ converge to $\mu$ as $n\rightarrow\infty$ with exponential rate $\varepsilon$ that depends on the potential functions $(U,V)$, while the one of forward/backward  Gibbs sampler discussed in Theorem~\ref{theo-S1n} {\it only depends on $U$}.
In generative Sinkhorn-type modeling, the perturbed data distribution  $\mu\Ka_0$ is approximating an easy-to-sample distribution $\eta=\lambda_V$ (often chosen as a Gaussian distribution) sometimes called the prior distribution~\cite{bortoli-heng}. As shown in (\ref{sinhorn-entropy-form-Sch}), the backward steps  $\pi_{2n+1}=\eta \Ka_{2n+1}$ of the generative model based on Sinkhorn transitions $\Ka_{2n+1}$ are initialized at the prior distribution.

In Example~\ref{examp-tau-t-22} we initiate a comparison of the rates obtained in Theorem~\ref{th2-intro} in the context of Sinkhorn semigroups with a prescribed target and the ones presented in 
Theorem~\ref{th3-intro} for proximal Gibbs-type samplers.
For instance, when $\beta=I$ and $\tau$ is diagonal, we show that
 $$
 \Vert u_+\Vert_2>\Vert v_+\Vert_2\Longleftrightarrow
(1+ \varepsilon^{-1})^{-1}<   b_{u_+}(\beta,\tau)^2.
 $$

Here again the backward sampler is an  idealized samplers since $\Va_s$ are generally unknown. To sample these backward denoising diffusions we often use
 the score matching trick introduced by~Hyv\"arinen and Dayan in~\cite{dayan} allows to estimate $\nabla \Va_s$ by some smooth score function $f_{s,\theta}$ indexed by some $\theta$ on some state space $\Theta$. In this context, 
 a neural network is often used to learn $\Va_s$.

\section{Transportation inequalities}\label{transport-sec}
This section collects the notions of entropy, Fisher information, and transportation costs that are used throughout the paper, together with the main functional inequalities linking them. We first introduce relative entropy, Fisher information, Wasserstein distances, and the entropic and Fisher transportation costs associated with a reference Markov transition. We then recall quadratic transportation inequalities and log-Sobolev inequalities, both in their classical and local forms, and discuss their relationships and equivalent formulations for Markov kernels. Several illustrative examples are provided, covering Gibbs measures and Gaussian transitions, with an emphasis on conditions ensuring log-Sobolev and transportation inequalities. Finally, we introduce the Schrödinger bridge framework and establish basic entropy estimates that will serve as key tools in the analysis of entropic optimal transport maps in later sections.
\subsection{Entropy and divergence criteria}
 This section introduce the basic notions of divergence and transportation costs 
 between pairs of probability measures used in the article.
\begin{itemize}
\item The relative entropy (a.k.a. Kullback-Leibler divergence and $I$-divergence) of  $\nu$ with respect to some $\mu\gg \nu$ is defined by the    formula
 $$
 \Ha(\nu~|~\mu):=\int~\log{\left(\frac{d\nu}{d\mu}(x)\right)}~\nu(dx).
$$
We also use the convention $  \Ha(\nu~|~\mu)=\infty$ when $\nu\not\ll \mu$.   
\item The (relative) Fisher information of $\nu$ with respect to some $\mu\gg \nu$ is defined by
$$
\Ja(\nu~|~\mu)= \mu\left(\frac{\Vert\nabla f\Vert_2^2}{f}\right)=\nu\left(\Vert\nabla \log{f}\Vert_2^2\right)
\quad \mbox{\rm with}\quad f=d\nu/d\mu.
$$
We also use the convention $  \Ja(\nu~|~\mu)=\infty$ when $\nu\not\ll \mu$.   
\item The $p$-th Wasserstein distance  associated with the Euclidian norm is defined for $p\geq 1$ by the formula
\begin{equation}\label{def-c2}
\Da_{p}(\nu,\mu)^p:= \inf_{P \in \Pi(\nu,\mu)}~P(e_p)~\quad
\mbox{\rm with}\quad e_p(x_1,x_2):=\Vert x_1-x_2\Vert_{2}^p.
\end{equation}

 \item The Kantorovich transport cost from $\nu$ to $\mu$ associated with the entropic transport cost
 is defined by
$$
\Ea_{\Ka}(\nu|\mu):=\inf_{P \in \Pi(\nu,\mu)}~P(h_{\Ka})\quad
\mbox{\rm with}\quad
h_{\Ka}(x_1,x_2):= \Ha(\delta_{x_1}\Ka~|~\delta_{x_2}\Ka).
$$
\item Given a Markov transitions $K$ such that $\delta_{x_1}\Ka\ll \delta_{x_2}\Ka$ for any $x_1,x_2\in\RR^d$ 
the Kantorovich transport cost from $\nu$ to $\mu$ associated with the Fisher information transport cost is defined by
\begin{equation}\label{def-fisher}
\Fa_{\Ka}(\nu|\mu):=\inf_{P \in \Pi(\nu,\mu)}~P(\jmath_{\Ka})
\quad
\mbox{\rm with}\quad
\jmath_{\Ka}(x_1,x_2):= \Ja(\delta_{x_1}\Ka~|~\delta_{x_2}\Ka).
\end{equation}

\end{itemize}

For any Markov transition $\Ka$ and any measures $\nu,\mu\in\Ma_1(\RR^d)$ we recall that
\begin{equation}\label{prop-Hent}
 \Ha(\nu\Ka~|~\mu\Ka)\leq  \Ha(\nu~|~\mu)\quad\mbox{\rm and}\quad
 \Ha(\mu\Ka~|~\nu)\leq \int~\mu(dx)~\Ha(\delta_{x}\Ka~|~\nu).
\end{equation}
For the convenience of the reader the detailed proofs of the above formulae are provided in the appendix on page~\pageref{prop-Hent-proof}.
\subsection{Quadratic and log-Sobolev inequalities}\label{sec-tineq}
 This section formally introduces the basic transportation inequalities used in the article and presents some properties and equivalent formulations. For a more thorough discussion on this subject we refer to~\cite{toulouse-team,bgl,cattiaux-14,cattiaux,bolbeault,ledoux,royer} and references therein.
\begin{defi}\label{def-log-sob}
We say that a measure $\mu\in\Ma_1(\RR^d)$ satisfies the quadratic transportation cost inequality $\TT_2(\rho)$ with constant $\rho>0$ if for any $\nu\in\Ma_1(\RR^d)$ we have
\begin{eqnarray}
\frac{1}{2}~\Da_{2}(\nu,\mu)^2&\leq& \rho~ \Ha(\nu~|~\mu).\label{K-talagrand}
\end{eqnarray}
We say that a Markov transition $\Ka$ satisfies the quadratic transportation cost inequality $\TT_2(\rho)$ when $\delta_x\Ka$ satisfies a quadratic transportation cost inequality $\TT_2(\rho)$ for any $x\in\RR^d$.
\end{defi}
\begin{defi}\label{def-log-sob-2}
We say that a measure $\mu\in\Ma_1(\RR^d)$ satisfies the log-Sobolev inequality $LS(\rho)$  if for any $\nu\in\Ma_1(\RR^d)$  we have
\begin{eqnarray}
 \Ha(\nu~|~\mu)
&\leq &\frac{\rho}{2}~
\Ja(\nu~|~\mu).\label{K-log-sob}
\end{eqnarray}
We say that a Markov transition $\Ka$ satisfies  the log-Sobolev inequality $LS(\rho)$ when $\delta_x\Ka$ satisfies a   log-Sobolev inequality $LS(\rho)$ for any $x\in\RR^d$.

We also say that $\Ka$ satisfies  the local log-Sobolev inequality $LS_{\tiny loc}(\rho)$ when we have for any $\mu,\eta\in \Ma_1(\RR^d)$
\begin{equation}\label{lem-H0-intro-def}
\Ea_{\Ka}(\mu~|~\eta)\leq \frac{\rho}{2}~\Fa_{\Ka}(\mu~|~\eta).
\end{equation}
\end{defi}

In (\ref{K-talagrand}) and (\ref{K-log-sob}) it is clearly enough to consider measures $\nu\ll\mu$.
Also note that $\Ka$ satisfies  the local log-Sobolev inequality $LS_{\tiny loc}(\rho)$ if and only if
for any $x_1,x_2\in\RR^d$ we have the inequalities
\begin{equation}\label{reflocls}
\Ha(\delta_{x_1}\Ka~|~\delta_{x_2}\Ka)\leq \frac{\rho}{2}~\Ja(\delta_{x_1}\Ka~|~\delta_{x_2}\Ka).
\end{equation}
This  condition is clearly met when $\Ka$ satisfies  the log-Sobolev inequality $LS(\rho)$. 

\begin{defi}\label{def-jj-2}
We say that $\Ka$ satisfies the 
Fisher-Lipschitz inequality  $\JJ_2(\kappa)$ with constant $\kappa>0$  if 
for any $x_1,x_2\in\RR^d$ we have the inequalities
\begin{equation}\label{lip-cond-H0-jmat}
\Ja(\delta_{x_1}\Ka~|~\delta_{x_2}\Ka)\leq 
  \kappa^2~\Vert x_1-x_2\Vert_2^2.
\end{equation}
\end{defi}

A theorem by Otto and Villani, Theorem 1 in~\cite{otto-villani}, ensures that the  log-Sobolev inequality $LS(\rho)$
implies the the quadratic transportation cost inequality $\TT_2(\rho)$.

Note that the strength of the inequalities $\TT_2(\rho)$ and $LS(\rho)$ defined in (\ref{K-talagrand}) and (\ref{K-log-sob}) are inversely proportional to the parameter $\rho$; that is we have  $\TT_2(\rho_1)\Longrightarrow \TT_2(\rho_2)$  and $LS(\rho_1)\Longrightarrow LS(\rho_2)$ 
as well as $LS_{\tiny loc}(\rho_1)\Longrightarrow LS_{\tiny loc}(\rho_2)$  for any $\rho_1\leq \rho_2$. 

Next lemma provides an equivalent formulation
of a quadratic transportation cost inequality for Markov integral operators.
\begin{lem}\label{lem-Int-log-sob}
The Markov transition $\Ka$ satisfies the $\TT_2(\rho)$ inequality with constant $\rho>0$ if and only if for any distribution $\mu\in\Ma_1(\RR^d)$ 
and any Markov transition $\La$  we  have 
\begin{equation}\label{lee-0}
\frac{1}{2}~\Da_{2}(\mu \La,\mu \Ka)^2\leq \rho ~\Ha(\mu\times \La~|~\mu\times \Ka).
 \end{equation} 
 \end{lem}
\proof
Assume $\Ka$ satisfies a quadratic transportation cost inequality $\TT_2(\rho)$.
Applying Talagrand's inequality (\ref{K-talagrand}) to $(\nu,\mu)=(\delta_x\La,\delta_x\Ka)$  we check that
$$
 \Ha(\delta_x\La~|~\delta_x\Ka)\geq \frac{1}{2 \rho}~\Da_{2}(\delta_x\La, \delta_x\Ka)^2.
$$
This implies that
\begin{eqnarray*}
\Ha(\mu\times \La~|~\mu\times \Ka)
&=&\int~\mu(dx)~  \Ha(\delta_x\La~|~\delta_x\Ka)\\
&\geq& \frac{1}{2\rho }~\int~\mu(dx)~\Da_{2}(\delta_x\La,\delta_x\Ka)^2\geq \frac{1}{2\rho }~\Da_{2}(\mu \La,\mu \Ka)^2.\end{eqnarray*}
This ends the proof of (\ref{lee-0}). In the reverse angle when (\ref{lee-0}) holds, choosing $\mu=\delta_x$ and $\La$ such that $\delta_x\La=\nu$ we have
$$
\frac{1}{2}~\Da_{2}(\nu,\delta_x \Ka)^2\leq \rho 
~\Ha(\mu\times \La~|~\mu\times \Ka)=  \rho~\Ha(\nu~|~\delta_x\Ka).
$$
This shows that $\Ka$ satisfies a quadratic transportation cost inequality $\TT_2(\rho)$.
This ends the proof of (\ref{lee-0}).

This ends the proof of the lemma.
\cqfd

\subsection{Some illustrations}

By a theorem of Bobkov, Theorem 1.3 in~\cite{bobkov}, a probability  measure $\mu=\lambda_U$ with  $\nabla^2U\geq 0$ satisfies the $LS(\rho)$-inequality for some $\rho>0$ if and only if there exists some constant $\delta>0$ such that
\begin{equation}\label{bobkov-condition}
\int~\exp{(\delta\Vert x\Vert_2^2)}~\mu(dx)<\infty.
\end{equation}
Consider a probability  measure $\mu=\lambda_U$ with
a smooth semiconvex potential, in the sense that $U(x)+\frac{a}{2}\Vert x\Vert_2^2$  is convex for some $a>0$ (so that $\nabla^2U\geq -a~I$). In this situation,  Aida~\cite{aida} and Wang~\cite{wang} (see also~\cite{ledoux-99,ledoux}) show that $\lambda_U$ satisfies the $LS(\rho)$-inequality for some $\rho>0$ as soon as there exists some constant $\delta>0$ such that 
\begin{equation}\label{aida-condition}
\int~\exp{(2(\delta+a)\Vert x\Vert_2^2)}~\mu(dx)<\infty.
\end{equation}

\begin{examp}\label{nabla2-logsob}
Consider Gibbs measures $\mu=\lambda_U$ associated with a smooth convex potential $U$ satisfying $\nabla^2 U(x)\geq u^{-1}_+$,
for some $u_+>0$.  This condition ensures that $\mu$ and  satisfies the  $\TT_2(\rho)$-inequality (\ref{K-talagrand}) as well as a  the $LS(\rho)$-inequality (\ref{K-log-sob}) with $\rho=\Vert u_+\Vert$; see for instance Corollary 9.3.2 and Corollary 5.7.2 in~\cite{bgl}.
\end{examp}
\begin{examp}\label{tensorization-examp}
Log-Sobolev inequalities are well-adapted for design dimension free estimates. For instance,
consider Gibbs measures of the form
$
\lambda_{U}=
\lambda_{U^1}\otimes
\lambda_{U^2}
$
for some potential functions $U^1$ and $U^2$ on $\RR^{d_1}$ and $\RR^{d_2}$ with
 $d=d_1+d_2$. Equivalently, we have $U(x)=U^{1}(x^1)+U^2(x^2)$ for any $x=(x^1,x^2)\in\RR^{d}=\RR^{d_1+d_2}$.
 By the tensorization principle if $\lambda_{U^i}$ satisfies the log-Sobolev inequality $LS(\rho^i)$ for $i=1,2$ then $\lambda_{U}$ satisfies the log-Sobolev inequality $LS(\rho^1\vee\rho^2)$. More general criteria for the log-Sobolev inequality  on product spaces based on the Hessian of conditional distributions are also provided in~\cite{otto}.
\end{examp}
\begin{examp}\label{nabla2-logsob-v2}
By the Holley-Stroock perturbation lemma~\cite{holley}, Gibbs measures of the form $\mu=\lambda_{U+\overline{U}}$ with $\mbox{\rm osc}(\overline{U}):=\sup_{x\in\RR^d}\overline{U}(x)-
\inf_{x\in\RR^d}\overline{U}(x)<\infty$
and $\nabla^2 U(x)\geq u^{-1}_+$ also satisfy the $LS(\rho)$-inequality (\ref{K-log-sob}) with $\rho=\Vert u_+\Vert~ \exp{(\mbox{\rm osc}(\overline{U}))}$ (see Corollary 1.7 in~\cite{ledoux}, as well as Theorem 1.1 in~\cite{cattiaux}). 

If $ \overline{U}$ is unbounded but satisfies
an exponential bound $\int~\lambda_U(dx)~e^{\alpha\Vert \nabla \overline{U}(x)\Vert}<\infty$ for $\alpha$ large enough, then $\lambda_{U+\overline{U}}$ also satisfies a log-Sobolev inequality (see Remark 21.5 in~\cite{villani-2}).
 For instance,  one-dimensional double-well potential functions $U(x)=x^4-a~x^2$, with $a>0$ are convex at infinity and thus satisfy 
the $LS(\rho)$-inequality for some $\rho>0$.  \end{examp}
Several criteria with explicit estimates of the log-Sobolev constants of Gibbs measures with non uniformly convex potentials are available in the literature under rather mild conditions, including the celebrated Bakry-Emery criteria on Riemannian manifolds~\cite{bakry-emery}. It is clearly out of the scope of the present article to review these estimates, see for instance~\cite{cattiaux-14,cattiaux,monmarche,otto,villani,villani-2} and references therein. 
Next we discuss two important and rather simple criteria covering a large class of models.

$\bullet$ Assume there exists some
$a>0$ and $b\geq 1$ such that for all $(x_1,x_2)$ we have
\begin{equation}\label{ex-lg-cg}
(x_1-x_2)^{\prime}(\nabla U(x_1)-\nabla U(x_2))\geq a~\Vert x_1-x_2\Vert^{2b}_2.
\end{equation}
In this situation, Theorem 9 in~\cite{cattiaux-14} provides explicit estimates of the
log-Sobolev constant of $\lambda_U$ in terms of the parameters $(a,b)$.

$\bullet$ Assume there exists some $a>0$, $b\geq 0$ and $\delta\geq 0$ such that 
\begin{equation}\label{ex-lg-cg-v2}
\nabla^2 U(x)\geq 
a~1_{\Vert x\Vert_2 \geq  \delta}~I-b 
~1_{\Vert x\Vert_2< \delta}~I.
\end{equation}
In this context there exists some $\overline{U}$ with bounded oscillations such that
$\nabla^2 (U+\overline{U})\geq (a/2)~I$. By the Holley-Stroock perturbation lemma~\cite{holley},
$\lambda_U$ satisfies the $LS(\rho)$-inequality (\ref{K-log-sob}) with some parameter $\rho=\rho(a,b,\delta)=2a^{-1}~ \exp{(\mbox{\rm osc}(\overline{U}))}$ that depends on $(a,b,\delta)$.
 For an explicit  construction of such function we refer to~\cite{monmarche-AHL}.
 We also refer to~\cite{monmarche} for a variety of log-Sobolev criteria when the curvature is not positive everywhere based on  Holley-Stroock and Aida-Shigekawa perturbation arguments~\cite{aida-shi,holley}.

A well known key feature of 
the log-Sobolev inequality is that the tensor product of a measure satisfying a $LS(\rho)$-inequality also satisfy the $LS(\rho)$-inequality with the same constant (thus independent of the dimension).
As underlined in~\cite{ledoux} in multivariate settings the perturbation estimates discussed in Example~\ref{nabla2-logsob-v2} behave poorly as functions of the dimension.

By (\ref{bobkov-condition}) the transition $K_W$ with $\nabla^2_2W\geq 0$ satisfies the $LS(\rho)$-inequality for some $\rho>0$  as soon as  there exists some constant $\delta>0$ such that
\begin{equation}\label{bobkov-condition}
\sup_{x\in\RR^d}\int~\exp{(\delta\Vert y\Vert^2)}~K_W(x,dy)<\infty.
\end{equation}

\begin{examp}\label{examp-logsob}
 Following the discussion given in Example~\ref{nabla2-logsob}, Markov transitions of the form $\Ka=K_W$ with $\nabla^2_2W(x,y)\geq \rho^{-1}~I$ also satisfy the $LS(\rho)$ inequality.
In this context, Markov transitions $K_{W+\overline{W}}$ with $\Vert \overline{W}\Vert<\infty$
 also satisfy the $LS(\overline{\rho})$ inequality with $\overline{\rho}=\rho~e^{2\Vert \overline{W}\Vert}$. More generally, Markov transitions $K_W$  associated with a potential function $y\mapsto W(x,y)$ convex at infinity uniformly w.r.t. the first variable also  satisfy 
the $LS(\rho)$-inequality for some $\rho>0$.  
\end{examp}

\begin{rmk}
Note that the proof of the estimate (\ref{lee-0}) is immediate when the distribution $\mu\Ka$ satisfies the  $\TT_2(\rho)$ inequality.
Indeed in this case we have
\begin{equation}\label{lee-0-v2}
\frac{1}{2}~\Da_{2}(\mu \La,\mu \Ka)^2\leq \rho ~\Ha(\mu\La~|~\mu\Ka)
\leq \rho ~\Ha(\mu\times \La~|~\mu\times \Ka).
 \end{equation} 
 Nevertheless,  the  $\TT_2(\rho)$ inequality is not always stable by Markov transport, in the sense that
 the distribution  $\mu\Ka$ may not satisfy the  $\TT_2(\rho)$ inequality
even if Markov transition $\Ka$ satisfy
 the  $\TT_2(\rho)$ inequality.
\end{rmk}

\begin{examp}\label{examp-bx}
Consider a Gaussian transition $\Ka=K_{W}$ of the form  
$$
\begin{array}{l}
W(x,y)=-\log{g_{\tau}(y-b(x))}\\
\\
\Longrightarrow
\nabla_2W(x,y)=\tau^{-1}(y-b(x))\quad \mbox{and}\quad
\nabla_2^2W(x,y)=\tau^{-1}\geq \rho^{-1}~I\quad \mbox{with}\quad\rho:=\Vert\tau\Vert,
\end{array}$$
for some drift function $b(x)$ such that
$$
(b(x_1)-b(x_2))^{\prime}\tau^{-2}
(b(x_1)-b(x_2))\leq \kappa^2~\Vert x_1-x_2\Vert_2^2
\quad \mbox{with}\quad
\kappa= \Vert \tau^{-1}\Vert_2~\mbox{\rm lip}(b),
$$
for some Lipschitz constant $\mbox{\rm lip}(b)$. In this case, we have
the Fisher information formula
\begin{eqnarray*}
\jmath_{\Ka}(x_1,x_2)
&=&\int~\Ka(x_1,dy)~\Vert\nabla_2W(x_1,y)-\nabla_2W(x_2,y)\Vert_2^2\\
&=&
(b(x_1)-b(x_2))^{\prime}\tau^{-2}
(b(x_1)-b(x_2))\leq \kappa^2~e_2(x_1,x_2).
\end{eqnarray*}
In addition, we have the Gaussian entropy formula
\begin{eqnarray*}
h_{\Ka}(x_1,x_2)&=& \Ha(\delta_{x_1}\Ka~|~\delta_{x_2}\Ka)=\int \Ka(x_1,dy)~\left(W(x_2,y)-W(x_1,y)\right)\\
& =&\frac{1}{2}
 \left(b(x_1)-b(x_2)\right)^{\prime}\tau^{-1}\left(b(x_1)-b(x_2)\right)\leq  \frac{\rho}{2}~\jmath_{\Ka}(x_1,x_2)\quad \mbox{with}\quad
 \rho:=\Vert \tau\Vert_2.
\end{eqnarray*}
The estimate $\nabla_2^2W(x,y)\geq \rho^{-1}~I$ also ensures that  $\Ka=K_W$ satisfies the log-Sobolev inequality $LS(\rho)$. Applying (\ref{K-log-sob}) we recover the fact that
$$   
h_{\Ka}(x_1,x_2)
\leq 
\frac{\rho}{2}~\jmath_{\Ka}(x_1,x_2).
$$
\end{examp}
\subsection{Schr\" odinger bridge maps}\label{sec-sbm}
The forward and backward disintegration of  the entropic optimal transport bridge $P_{\mu,\eta}$ with reference measure $P$ from $\mu$ to $\eta$ defined in (\ref{ob-eq}) are given by the product formulae
$$
P_{\mu,\eta}=\mu\times P_{\eta|\mu}.
$$
In the above display $P_{\eta|\mu}$ is the conditional distribution of the second coordinate of the bridge given the first and $P_{\eta|\mu}$ is the conditional distribution of the first coordinate given the second. 
To underline the role of the reference measure,
note that for any Markov transition $\Ta$ such that $\mu\Ta=\eta$ we have
$$
R:=(\mu\times \Ta)\in \Pi(\mu,\eta)
\quad\mbox{\rm and}\quad
Q:=\mu\times (\Ta\Ka)\quad\Longrightarrow\quad
Q_{\mu,\mu\Ta\Ka}=Q_{\mu,\eta\Ka}=Q.
$$
Equivalently $Q_{\mu,\eta\Ka}$ is the 
Schr\" odinger bridge   from
$\mu$ to $\eta\Ka$ with reference measure $Q$ but this coupling measure is generally not the Schr\" odinger bridge  $P_{\mu,\eta\Ka}$ {\it with reference measure $P=(\mu\times \Ka)$.}
Nevertheless we  have the estimate
\begin{equation}
\begin{array}{l}
P:=(\mu\times \Ka)\qquad
Q:=\mu\times (\Ta\Ka)\in\Pi(\mu,\eta\Ka)\\
\\
\Longrightarrow
 \Ha(P_{\mu,\eta\Ka}~|~P)\leq \Ha(Q~|~P)\leq  R(h_{\Ka}) 
 \quad\mbox{\rm with}\quad
R:=(\mu\times \Ta)^{\flat}\in \Pi(\eta,\mu).
 \end{array}\label{EKC-0}
\end{equation}
To prove this claim we apply (\ref{prop-Hent}) to $\nu=\delta_{x_2}\Ka$ and $\delta_{x_2}\Ta=\mu$ to check that
\begin{eqnarray}
\displaystyle\Ha(Q~|~P)&=&\int\mu(dx_2)~ \Ha((\delta_{x_2}\Ta)\Ka ~|~\delta_{x_2}\Ka)\\ 
&\leq& \int~\mu(dx_2)~\Ta(x_2,dx_1)~\nonumber
 \Ha(\delta_{x_1}\Ka~|~\delta_{x_2}\Ka)= R(h_{\Ka}).\label{ref-compc}
\end{eqnarray}

\section{Stability of  Schr\" odinger bridges}\label{sec-continuity-sb}

This section is devoted to quantitative stability properties of Schrödinger bridges with respect to perturbations of the target marginals and the reference transition. We first introduce a set of flexible regularity assumptions on the reference Markov transition, based on local log-Sobolev, Fisher–Lipschitz, and quadratic transportation inequalities. Under these conditions, we establish a general entropic continuity theorem providing sharp entropy and Wasserstein bounds for Schrödinger bridge maps. We then derive several corollaries yielding Lipschitz-type estimates for entropic and Wasserstein distances, including contraction properties of Markov transport maps. Finally, we present concrete sufficient conditions ensuring the validity of our assumptions, discuss their relevance for quadratic-cost and Gaussian reference models, and relate our results to existing stability estimates in the literature.

\subsection{Regularity conditions}\label{reg-cond-section-ab}
Consider the following regularity conditions:\\

\noindent {\it $A_0(\Ka)$:  The transition $\Ka$ satisfies  the local log-Sobolev inequality $LS_{\tiny loc}(\rho)$ and the Fisher-Lipschitz inequality  $\JJ_2(\kappa)$  for some $\rho,\kappa>0$. }\\

\noindent{\it $A_1(\Ka)$:  Condition $A_0(\Ka)$ is satisfied for some parameters $(\kappa,\rho)$. In addition, the Markov transition $\Ka$ satisfies a $\TT_2(\overline{\rho})$ inequality.  }\\

The rather mild condition $A_0(\Ka)$ is met for a large class of models.  For instance, as shown in Example ~\ref{examp-bx} condition  $A_0(\Ka)$ is also met for any Gaussian transitions associated with a Lipschitz drift function.
For Markov transitions $\Ka=K_W$, the $\JJ_2(\kappa)$ inequality resumes to an $\LL_2(\delta_{x_1}\Ka)$-Lipschitz property of the functions $\nabla_2W(x,y)$ in the sense that
$$
\jmath_{\Ka}(x_1,x_2)
=\int~\Ka(x_1,dy)~\Vert\nabla_2W(x_1,y)-\nabla_2W(x_2,y)\Vert^2\leq  \kappa^2~e_2(x_1,x_2).
$$
For linear Gaussian reference transitions (\ref{r-OU}) we already mention (as a  consequence of  (\ref{nabla1-WN-ex})), that  for any $n\geq 0$ all the Sinkhorn transitions $\Ka_n$ satisfy
 the Fisher-Lipschitz inequality  $\JJ(\kappa)$  with the parameter
 \begin{equation}\label{kappa-lingauss}
 \kappa=\Vert \tau_t^{-1}e^{tA}\Vert_2
 \quad\mbox{\rm and}\quad \tau_t:=\int_0^t e^{sA}~\Sigma~ e^{sA^{\prime}}~ds
 \end{equation} 
Log-Sobolev inequalities are  a well studied technique for the stability analysis of Markov processes. The class of measures satisfying these inequalities is quite large, see for instance~\cite{toulouse-team,bgl,cattiaux-14,cattiaux, bolbeault,ledoux,royer} 
and the pioneering articles by Gross~\cite{gross,gross-2}.
 
 When $\Ka=K_W$, the log-Sobolev inequality $LS_{\tiny loc}(\rho)$ can be checked
easily using the Hessian estimates discussed in Example~\ref{examp-logsob}, see also the curvature estimates discussed in Sections~\ref{sec-lsob-s-kernels} and~\ref{sec-lsob-s-curvature} in the context of Sinkhorn transitions. 

\subsection{An entropic continuity theorem}\label{sec-entrop-theo}

By (\ref{lem-H0-intro-def}), condition $A_0(\Ka)$ ensures that for any probability measures
$\eta,\mu$ we have
\begin{equation}\label{lem-H0-intro}
\begin{array}{l}
\displaystyle\Ea_{\Ka}(\eta|\mu)\leq \frac{\rho}{2}~\Fa_{\Ka}(\eta|\mu)\quad\mbox{and by (\ref{lip-cond-H0-jmat})}\quad
\sqrt{\Fa_{\Ka}(\eta|\mu)}\leq \kappa~\Da_{2}(\eta,\mu)\\
\\
\displaystyle\Longrightarrow\quad
\Ea_{\Ka}(\eta|\mu)\leq \frac{\rho}{2}~\Fa_{\Ka}(\eta|\mu)\leq \frac{\rho}{2}~
\kappa^2~\Da_{2}(\mu,\eta)^2.
\end{array}
\end{equation}
In addition, when $\mu$  satisfies the  $\TT_2(\upsilon)$ inequality  (\ref{K-talagrand})
we have
\begin{equation}\label{A0-T2}
\Ea_{\Ka}(\eta|\mu)\vee \Ea_{\Ka}(\mu|\eta)\leq~\varepsilon_{\kappa}(\rho,\upsilon)~ \Ha(\eta~|~\mu).
\end{equation}

In the above display formula,  $\varepsilon_{\kappa}(\rho_1,\rho_2)$  stands for the collection of parameters indexed by $\kappa,\rho_1,\rho_2>0$ and defined by
\begin{equation}\label{A0-T2-def-vareps}
\varepsilon_{\kappa}(\rho_1,\rho_2):=\kappa^2\rho_1~\rho_2.
\end{equation}

On the other hand, the entropic cost function $h_{\Ka}$ being continuous there exists some entropic optimal transport  $R \in \Pi(\eta,\mu)$. Choosing in (\ref{EKC-0}) the distribution 
 \begin{equation}\label{dom-prop}
 R:=\argmin_{Q \in \Pi(\eta,\mu)}~Q(h_{\Ka})\quad \mbox{\rm we check that}\quad
 \Ha(P_{\mu,\eta\Ka}~|~P)\leq R(h_{\Ka})=
\Ea_{\Ka}(\eta|\mu).
 \end{equation}

We summarize the above discussion with the following quantitative upper and lower  bounds for 
  Schr\" odinger bridges.

\begin{theo}\label{lem-Int-log-sob-talagrand}
The Schr\" odinger bridge map $\eta\mapsto P_{\mu,\eta\Ka}$
from 
$\mu$ to $\eta\Ka$ and reference measure
$P:=(\mu\times \Ka)$ defined in (\ref{ot-eq-etaK-i}) satisfies for any probability measure $\eta$ the entropic continuity property
 \begin{equation}  
 \Ha(\eta \Ka~|~\mu\Ka)\leq
 \Ha(P_{\mu,\eta\Ka}~|~P_{\mu,\mu\Ka})\leq \Ea_{\Ka}(\eta~|~\mu).\label{intg-log-sob}
 \end{equation} 
 \begin{itemize}
 \item When condition $A_0(\Ka)$ holds for some  $(\kappa,\rho)$ we have
\begin{equation}  
 \Ha(P_{\mu,\eta\Ka}~|~P_{\mu,\mu\Ka})\leq \frac{\rho}{2}~
\kappa^2~\Da_{2}(\mu,\eta)^2.
 \label{intg-log-sob-T2-inter}
 \end{equation} 
\item When condition $A_0(\Ka)$ holds for some  $(\kappa,\rho)$ and  $\mu$  satisfies the  $\TT_2(\upsilon)$ inequality, we  have the entropy estimate
\begin{equation}  
 \Ha(P_{\mu,\eta\Ka}~|~P_{\mu,\mu\Ka})
 \leq \varepsilon_{\kappa}(\rho,\upsilon)~ \Ha(\eta~|~\mu).
 \label{intg-log-sob-T2}
 \end{equation} 
\item  When condition $A_0(\Ka)$ holds for some  $(\kappa,\rho)$ and  $\eta$  satisfies the  $\TT_2(\upsilon)$ inequality, we  have the entropy estimate
\begin{equation}  
 \Ha(P_{\mu,\eta\Ka}~|~P_{\mu,\mu\Ka})
 \leq \varepsilon_{\kappa}(\rho,\upsilon)~ \Ha(\mu~|~\eta).
 \label{intg-log-sob-T2-sym}
 \end{equation} 
\item  When condition $A_1(\Ka)$ holds for some  $(\kappa,\rho,\overline{\rho})$   we also have the Wasserstein estimates
\begin{equation}  
~\Da_{2}(\eta \Ka,\mu \Ka)^2\leq 2\,\overline{\rho}~ 
 \Ha(P_{\mu,\eta\Ka}~|~P_{\mu,\mu\Ka})\leq \varepsilon_{\kappa}(\rho,\overline{\rho})
~\Da_{2}(\eta,\mu)^2.\label{intg-log-sob-2}
 \end{equation} 
 \end{itemize}
 \end{theo}

\proof
The estimates (\ref{intg-log-sob}) are direct consequences of the inequalities stated in (\ref{dom-prop}). The estimate (\ref{intg-log-sob-T2-inter}) is now a consequence of 
 (\ref{lem-H0-intro}) and (\ref{intg-log-sob}).
The estimates (\ref{intg-log-sob-T2}) and (\ref{intg-log-sob-T2-sym}) are checked using (\ref{intg-log-sob})  and (\ref{A0-T2}).
 We also check the l.h.s. of (\ref{intg-log-sob-2}) applying (\ref{lee-0}) to $\La=P_{\eta\Ka|\mu}$. The r.h.s. of (\ref{intg-log-sob-2}) is checked combining (\ref{intg-log-sob}) with the estimates (\ref{lem-H0-intro}).
 This ends the proof of the theorem.\cqfd
 
 We illustrate Theorem~\ref{lem-Int-log-sob-talagrand} by considering the reference measure
$P:=(\mu\times \Ka)$ associated with the linear-Gaussian transition $\Ka=K_W$ introduced in (\ref{r-OU}) over a fixed time horizon $t>0$.
 Combining Example~\ref{examp-logsob} with the curvature estimate (\ref{HessW0}) we verify that $\Ka$
satisfies  the $LS(\Vert \tau_t\Vert_2)$ inequality for every $t>0$, where $\tau_t$ is defined  in (\ref{kappa-lingauss}).
In this setting, using (\ref{kappa-lingauss}), estimate (\ref{intg-log-sob-T2-inter}) implies that for any probability measure $\eta$,
 $$
  \Ha(\eta \Ka~|~\mu\Ka)\leq
  \Ha(P_{\mu,\eta\Ka}~|~P_{\mu,\mu\Ka})\leq \frac{1}{2}~
\Vert \tau_t^{-1}e^{tA}\Vert_2^2~\Vert \tau_t\Vert_2~\Da_{2}(\mu,\eta)^2
 $$
Moreover, if $\mu$  satisfies the  $\TT_2(\upsilon)$ inequality then (\ref{intg-log-sob-T2}) yields the following entropy contraction estimate
\begin{equation}  \label{ref-contract-tt}
  \Ha(\eta \Ka~|~\mu\Ka)\leq \Ha(P_{\mu,\eta\Ka}~|~P_{\mu,\mu\Ka})
 \leq \upsilon~
\Vert \tau_t^{-1}e^{tA}\Vert_2^2~\Vert \tau_t\Vert_2~ \Ha(\eta~|~\mu).
 \end{equation} 
$\bullet$ In the case $A=0\Longrightarrow \tau_t=t~I$ discussed in Section~\ref{reg-transp-sec}
we have
$$
(\ref{ref-contract-tt})\quad \Longleftrightarrow\quad
  \Ha(\eta \Ka~|~\mu\Ka)\leq \Ha(P_{\mu,\eta\Ka}~|~P_{\mu,\mu\Ka})
 \leq \frac{\upsilon}{t}~ \Ha(\eta~|~\mu).
$$
$\bullet$ If $A$ is an Hurwitz matrix, then there exists some $a,b>0$ such that
$\Vert e^{tA}\Vert_2\leq a~e^{-tb}$. In addition, for any $t_0>0$ we have
$c_{t_0}:=\sup_{t\geq t_0}{(\Vert \tau_t^{-1}\Vert_2^2 \vee\Vert \tau_t\Vert_2)}<\infty$. In this scenario for any $t\geq t_0$ we have
$$
(\ref{ref-contract-tt})\quad \Longrightarrow\quad
  \Ha(\eta \Ka~|~\mu\Ka)\leq \Ha(P_{\mu,\eta\Ka}~|~P_{\mu,\mu\Ka})
 \leq \upsilon~c_{t_0}~ a^2~e^{-2tb}~ \Ha(\eta~|~\mu).
$$

 \subsection{Reference flexibility}

Next corollary provides a way to estimate the regularity of  Schr\" odinger bridges 
$P_{\mu,\eta}$ for {\it some possibly unknown reference measure $P$} in terms of  the Schr\" odinger bridges
of the form $$
Q_{\mu,\eta}=\argmin_{R\,\in\, \Pi(\mu,\eta)}\Ha(R~|~Q),
$$ when 
  the reference measure $Q=(\mu\times \Ka)$  is associated with a Markov transition $\Ka$ such  that $\eta=\pi\Ka$ for some $\pi$. This result highlights the flexibility of our approach with respect to the choice of reference transition.

\begin{cor}\label{cor-key}
Consider the Schr\" odinger bridge $P_{\mu,\eta}$ defined in (\ref{ob-eq}) for some reference measure $P$ and a couple of probability measures $(\mu,\eta)$. 
Assume there exists reference measure of the form $Q=(\mu\times\Ka)$ with some Markov transition $\Ka$ and some probability measure $\pi$ such that
$$
P_{\mu,\pi\Ka}=Q_{\mu,\pi\Ka}
\quad \mbox{and}\quad\pi\Ka=\eta.
$$
In this situation, we have
\begin{equation}  
 \Ha(\pi\Ka~|~\mu\Ka)\leq \Ha(P_{\mu,\pi\Ka}~|~Q_{\mu,\mu\Ka})\leq  \Ea_{\Ka}(\pi|\mu).\label{intg-log-sob-2-ex}
 \end{equation}  
  \begin{itemize}
\item
 When  $A_0(\Ka)$ holds for some parameters $(\kappa,\rho)$ 
and  the probability measure $\mu$ satisfies the $\TT_2(\upsilon)$  inequality 
we have
\begin{equation}  
 \Ha(\pi\Ka~|~\mu\Ka)\leq \Ha(P_{\mu,\pi\Ka}~|~Q_{\mu,\mu\Ka})\leq \varepsilon_{\kappa}(\rho,\upsilon)~ \Ha(\pi~|~\mu).
\label{intg-log-sob-2-ex-cor}
 \end{equation}   
\item When $A_1(\Ka)$ holds for some parameters $(\kappa,\rho,\overline{\rho})$   we also have
\begin{equation}  
\Da_{2}(\pi \Ka,\mu \Ka)^2\leq {2\overline{\rho}~
 \Ha( P_{\mu,\pi\Ka}~|~Q_{\mu,\mu\Ka})}\leq \varepsilon_{\kappa}(\rho,\overline{\rho})
~\Da_2(\pi,\mu)^2.\label{intg-log-sob-2-cor}
 \end{equation} 
 \end{itemize}
\end{cor}
\proof
We have
$$
\begin{array}{l}
P_{\mu,\pi\Ka}
=Q_{\mu,\pi\Ka}\quad \mbox{\rm with the reference}\quad
Q=\mu\times \Ka=Q_{\mu,\mu\Ka}
\\
\\
\Longrightarrow
 \Ha(P_{\mu,\pi\Ka}~|~Q_{\mu,\mu\Ka})= \Ha(Q_{\mu,\pi\Ka}~|~Q_{\mu,\mu\Ka}).
\end{array}$$
The end of the proof of Corollary~\ref{cor-key} is now a direct application of Theorem~\ref{lem-Int-log-sob-talagrand} to the reference $Q$ (combined with (\ref{lem-H0-intro})). This ends the proof of the corollary.
\cqfd

\subsection{Some Lipschitz estimates}\label{lip-sec-A2}

This section is devoted to quantitative regularity properties of Markov transitions and their associated bridge and transport maps. Building on the general framework developed in Theorem~\ref{lem-Int-log-sob-talagrand}, we show that a broad class of Markov kernels satisfying suitable local functional inequalities induces Lipschitz-type controls at the level of entropy and quadratic transportation costs.

More precisely, the estimates established in Theorem~\ref{lem-Int-log-sob-talagrand} and its corollaries are remarkably flexible: they apply to general bridge maps associated with a reference Markov transition $\Ka$ satisfying a local logarithmic Sobolev inequality $LS_{\mathrm{loc}}(\rho)$. This setting naturally encompasses many kernels arising in stochastic dynamics, including transitions defined by strongly convex potentials and their perturbations. Within this framework, the interplay between logarithmic Sobolev inequalities, Talagrand-type transportation inequalities, and Fisher information bounds yields explicit contraction and stability estimates under the action of $\Ka$.
As a first consequence, we obtain Lipschitz continuity properties of the mapping $\mu \mapsto \mu \Ka$ with respect to both entropy and Wasserstein distances. These estimates can be interpreted as entropic and metric regularization effects induced by the Markov transition. In particular, they provide quantitative bounds on how entropy and quadratic transportation costs evolve under the action of $\Ka$, and clarify the precise role played by curvature, convexity, and Lipschitz-type assumptions on the transition kernel.
 We conclude the section by presenting a set of easily verifiable sufficient conditions, formulated in terms of the gradient and Hessian of the underlying potential, which ensure that the abstract assumptions required for these Lipschitz estimates are satisfied.

 \begin{cor}\label{cor-1-intro}
 Consider a Markov transition $\Ka$ satisfying $A_0(\Ka)$ for some parameters $(\rho,\kappa)$. In this situation, the following assertions are satisfied:
 
 \begin{itemize} 
 \item  For any  $\mu,\eta\in\Ma_1(\RR^d)$ we have the entropic transportation cost inequality
 \begin{equation}  
2~ \Ha(\eta \Ka~|~\mu\Ka)\leq~ \varepsilon_{\kappa}(\rho,1)~\Da_{2}(\mu,\eta)^2.\label{intg-log-sob-v2}
 \end{equation} 
 \item When $\mu$  satisfies a $\TT_2(\upsilon)$  inequality,   for any  $\eta\in\Ma_1(\RR^d)$ we have 
   \begin{equation}\label{lip-ent}
\Ha(\eta \Ka~|~\mu \Ka)~\leq~  \left(1\wedge\varepsilon_{\kappa}(\rho,\upsilon)\right)~ \Ha(\eta~|~\mu).
  \end{equation}
  \item  When $\eta$  satisfies the  $\TT_2(\upsilon)$ inequality, we  have the entropy estimate
\begin{equation}  
\Ha(\eta \Ka~|~\mu \Ka)
 \leq \left(1\wedge\varepsilon_{\kappa}(\rho,\upsilon)\right)~ \Ha(\mu~|~\eta).
 \label{intg-log-sob-T2-sym-2}
 \end{equation} 

  \item   When $A_1(\Ka)$ holds for some  $(\kappa,\rho,\overline{\rho})$, for any  $\mu,\eta\in\Ma_1(\RR^d)$ we have 
\begin{equation}  \label{lipD2}
\Da_{2}(\eta \Ka,\mu \Ka)^2
\leq   \varepsilon_{\kappa}(\rho,\overline{\rho})~\Da_2(\eta,\mu)^2.
 \end{equation} 

 \end{itemize}
 \end{cor}
 The estimate (\ref{intg-log-sob-v2})  is checked using  (\ref{lem-H0-intro}) and (\ref{intg-log-sob}).
 The estimates (\ref{lip-ent}), (\ref{intg-log-sob-T2-sym-2}) and (\ref{lipD2}) are direct consequence of
 (\ref{intg-log-sob-T2}) and (\ref{intg-log-sob-T2-sym}) and (\ref{intg-log-sob-2}).

Note that the entropic transportation cost inequality
 (\ref{intg-log-sob-v2}) "goes in the reverse direction" of the quadratic transportation cost inequality  (\ref{K-talagrand}).

We end this section with an easily checked sufficient condition in terms of 
the gradient  $\nabla_2W(x,y)$  and the Hessian  $\nabla_2^2W(x,y)$  of the function $y\mapsto W(x,y)$. \\
\\
{\it $A_2(\Ka)$:  The Markov transition has the form $\Ka=K_W$ for some potential function $W(x,y)$. In addition, there exists some $\rho>0$ such that for any $x,y\in \RR^d$ we have
\begin{equation}\label{strong-convex-2}
\nabla^2_2W\geq \rho^{-1}~I.
\end{equation} 
Moreover, there exists some $\kappa>0$ such that for any  $x_1,x_2\in\RR^d$ and $y\in\RR^d$ we have 
\begin{equation}\label{lip-cond-H0}
 \Vert \nabla_2W(x_1,y)- \nabla_2W(x_2,y)\Vert_2\leq \kappa~\Vert x_1-x_2\Vert_2.
\end{equation} }

 The convexity condition (\ref{strong-convex-2}) ensures that the Markov transition $K_W$ satisfies the   $\TT_2(\rho)$-inequality as well as the $LS(\rho)$-logarithmic
Sobolev inequalities. The Lipschitz condition (\ref{lip-cond-H0}) ensures that $K_W$ satisfies  the Fisher-Lipschitz inequality  $\JJ_2(\kappa)$.

We conclude that
\begin{equation}\label{H01}
 A_2(\Ka)\Longrightarrow \left( A_1(\Ka)\quad \mbox{\rm with}\quad \overline{\rho}=\rho\right)\Longrightarrow 
  A_0(\Ka).
\end{equation}
Following the perturbation argument discussed in Example~\ref{examp-logsob} 
weaker sufficient conditions of  $A_1(\Ka)$ can be designed replacing
the Hessian estimate condition (\ref{strong-convex-2}) by an asymptotic convex condition.

Note that we can similarly consider in (\ref{EKC-0}) and (\ref{dom-prop})  a Wasserstein optimal coupling. In this case, we loose the entropy upper bound in (\ref{intg-log-sob-T2}).  Nevertheless, whenever $A_1(\Ka)$ holds for some  parameters $(\kappa,\rho,\overline{\rho})$ we obtain the Wasserstein estimate (\ref{intg-log-sob-2}) combining the log-Sobolev inequality (\ref{reflocls}) with Fisher information cost estimate (\ref{lip-cond-H0-jmat}).

 The proof of (\ref{lem-H0-intro}) and (\ref{A0-T2}) as well as of (\ref{intg-log-sob-2}) and (\ref{lip-ent}) combine optimal entropic cost couplings, Talagrand's quadratic transportation  inequality and the $LS(\rho)$-logarithmic
Sobolev inequality. Closely related arguments based on Wasserstein couplings are employed in~\cite{lee} and~\cite{chewi-phd} to establish contraction properties of proximal samplers with target product measures of the form (\ref{proximal-ex}); see for instance,  the proof of Lemma 2 in~\cite{lee} and the proof of Theorem 4.3.2 in~\cite{chewi-phd}.

 Finally, note that condition 
$A_0(\Ka)$ ensures that both the entropic and Fisher information cost functions associated with the transition 
$\Ka$ are continuous and dominated by the quadratic cost. As shown in \eqref{dom-prop}, the proof of \eqref{intg-log-sob} combines this property with the existence of either an entropic or a Wasserstein optimal coupling. Such gluing arguments are standard in the optimal transport literature; see, for example,~\cite{letrouit}, Lemma~7.6 in~\cite{villani}, Lemma~1 in~\cite{cuturi}, and the proof of Theorem~1.1 in~\cite{chiarini}.

\section{Stability of Sinkhorn bridges}\label{sinkhorn-sec}
 
 \subsection{Some regularity properties}\label{sinkhorn-sec-reg}
This section is concerned with the stability properties of Sinkhorn bridges (\ref{sinhorn-entropy-form}) associated with a general  reference measure. Our regularity conditions are expressed in terms of the conditions $A_i(\Ka)$ discussed in Section~\ref{sec-entrop-theo} and Section~\ref{lip-sec-A2}:\\

 {\it $A_0(\Sa):$  For any $n\geq 1$,  $A_0(\Ka_{n})$ is satisfied
for some  $(\kappa_{n},\rho_{n})$.
 }\\

 {\it $A_1(\Sa):$  For any $n\geq 1$,  $A_1(\Ka_{n})$ is satisfied
for some  $(\kappa_{n},\rho_{n},\overline{\rho}_n)$.
 }\\
 
  {\it $A_2(\Sa):$  For any $n\geq 1$,  $A_2(\Ka_{n})$  is satisfied for some  $(\kappa_{n},\rho_{n})$.
 }\\

By (\ref{H01}) we have
$$
 A_2(\Sa)\Longrightarrow \left( A_1(\Sa)\quad \mbox{\rm with}\quad \overline{\rho}_n=\rho_n\right)\Longrightarrow 
  A_0(\Sa)
  $$
  
  We also consider the conditions:\\
  
    {\it $A_U(\Sa):$ For any $n\geq 1$, $A_0(\Ka_{2n})$ holds for some  $(\kappa_{2n},\rho_{2n})$ and 
 $\lambda_U$ satisfies the $\TT_2(\Vert u_+\Vert_2)$ inequality for some $u_+\in\Sa_d^+$.}\\
 
 {\it $A_V(\Sa):$ For any $n\geq 1$, $A_0(\Ka_{2n-1})$ holds for some  $(\kappa_{2n-1},\rho_{2n-1})$ and 
 $\lambda_V$ satisfies the $\TT_2(\Vert v_+\Vert_2)$  inequality for some $v_+\in\Sa_d^+$.}\\

\begin{defi}
When any of the conditions $A_U(\Sa)$, $A_V(\Sa)$ or $A_i(\Sa)$  is met for some parameters  $(\kappa_{n},\rho_{n})$ and $\overline{\rho}_n$, for any $n\geq 1$ we set
\begin{equation}\label{def-varepsi-SW}
\begin{array}{rcl}
\epsilon_{2n}&:=&\varepsilon_{\kappa_{2n}}(\rho_{2n},\Vert u_+\Vert_2)\\\\
 \epsilon_{2n-1}&:=&\varepsilon_{\kappa_{2n-1}}(\rho_{2n-1},\Vert v_+\Vert_2)
 \quad \mbox{and}\quad
\overline{\epsilon}_{n}:=\varepsilon_{\kappa_{n}}(\rho_{n},\overline{\rho}_n).
 \end{array}
\end{equation}
with the parameters $\varepsilon_{\kappa}(\rho_1,\rho_2)$ indexed by $\kappa,\rho_1,\rho_2\geq 0$ defined in (\ref{A0-T2-def-vareps}).
\end{defi}

We shall also use the following condition:\\
 
 {\it $A_0^{\star}(\Sa):$  Conditions $A_U(\Sa)$ and $A_V(\Sa)$ hold
for some parameters $(\rho_{n},\kappa_{n})$ such that
\begin{equation}\label{epsilon-star}
\epsilon_{\star}:=\sup_{n\geq 1}\left(\varepsilon_{\kappa_{2n}}(\rho_{2n},\Vert u_+\Vert_2)\vee
\varepsilon_{\kappa_{2n-1}}(\rho_{2n-1},\Vert v_+\Vert_2)\right)
<\infty.
\end{equation}
 }

  We illustrate these rather abstract conditions in 
   the context of reference transitions in Section~\ref{general-sec} and "uniformly convex at infinity" potential functions $(U,V)$ satisfying (\ref{ex-lg-cg-intro-uv-lls}).

   A more detailed discussion  in the context of a linear Gaussian reference transition $K_W$ associated with some potential $W$ of the form (\ref{def-W})   is provided in Section~\ref{sec-lingauss}. 
In this situation, the above conditions can be easily checked in terms of the marginal regularity conditions discussed in Section~\ref{marg-reg-conditions}. As we we shall see in  (\ref{v-cinfty})  (\ref{u-cinfty}) and (\ref{c0a0}), 
$$
\begin{array}{rcccl}
C_U(V)&\Longrightarrow&
C_U^-(V)&\Longrightarrow&
A_U(\Sa)\quad\mbox{\rm with $\sup_{n\geq 1}(\kappa_{2n}\vee \rho_{2n}\vee \epsilon_{2n})<\infty$},\\
\\
C_V(U)&\Longrightarrow&
C_V^-(U)&\Longrightarrow&
A_V(\Sa)\quad\mbox{\rm with $\sup_{n\geq 1}(\kappa_{2n-1}\vee\rho_{2n-1}\vee  \epsilon_{2n-1})<\infty$},
\\
\\
 C_0(U,V)&\Longrightarrow& C_0^-(U,V)&\Longrightarrow&(
A_U(\Sa), A_V(\Sa)~\mbox{\rm and}~A_0^{\star}(\Sa) ),
\end{array}
$$
as well as 
$$
C_0(U,V)\Longrightarrow(A_2(\Sa)~\mbox{\rm and}~A_0^{\star}(\Sa))
$$
As expected from the log-Sobolev criteria discussed in Section~\ref{sec-tineq} (cf. Examples~\ref{nabla2-logsob},~\ref{nabla2-logsob-v2} and~\ref{examp-logsob}), these regularity conditions will be checked using a variety of gradient and Hessian formulae relating Sinkhorn potentials to the potential functions $(U,V)$ of the target measures.

We emphasize that all regularity assumptions of the form $A_U(\Sa)$, $A_V(\Sa)$, or $A_i(\Sa)$ are formulated solely in terms of corresponding regularity conditions $A_i(\Ka_n)$ imposed on the family of Sinkhorn transitions $\Ka_n$. Owing to the monotonicity properties \eqref{e12} and \eqref{ent-r1}, all contraction-type results established in this article remain valid even if the conditions $A_i(\Ka_n)$ are satisfied only asymptotically, namely for all sufficiently large indices $n \geq n_0$, for some $n_0 \geq 0$.

\subsection{Sinkhorn potential functions}\label{sec-pf}
Next we recall some basic well known properties of Sinkhorn bridges.
For the convenience of the reader, detailed proofs in our framework are provided in the appendix~\ref{review-sinkhorn} on page~\pageref{review-sinkhorn}.

Consider marginal measures $(\mu,\eta)=(\lambda_U,\lambda_V)$ and a reference transition $\Ka=K_W$ associated with some potential functions $(U,V,W)$.
Sinkhorn bridges are expressed  in terms of a sequence of potential functions $(U_n,V_n)$ defined by the  integral recursions
\begin{eqnarray}
U_{2n}~&:=&U+\log{K_W(e^{-V_{2n}})}=U_{2n+1},\nonumber\\
V_{2n+1}&:=&V+\log{K_{W^{\flat}}(e^{-U_{2n+1}})}=V_{2(n+1)}.
\label{prop-schp}
\end{eqnarray}
with the initial condition $V_0=0$ and $W^{\flat}$ as in (\ref{ref-Wflat}).

With these notations at hand, 
the Markov transitions in (\ref{sinhorn-entropy-form-Sch}) are given by the formulae
$$
\Ka_{2n}:=K_{W_{2n}}
\quad
\mbox{\rm and}\quad
\Ka_{2n+1}:=K_{W^{\flat}_{2n+1}},
$$
with  the potential functions $W_n$ defined by
\begin{eqnarray}
  U(x)+W_{2n}(x,y)&:=&U_{2n}(x)+W(x,y)+V_{2n}(y),\nonumber
\\
  V(y)+W_{2n+1}^{\flat}(y,x)&:=&U_{2n+1}(x)+W(x,y)+V_{2n+1}(y).\label{log-densities}
\end{eqnarray}
Note that $W_0=W$ so that $\Pa_0=P$ coincides with the reference measure (\ref{def-triple}). 

In terms of the potential functions $(U_n,V_n)$ defined in (\ref{prop-schp}) the distributions (\ref{transport-g}) have the Boltzmann-Gibbs formulations
\begin{equation}\label{pi-gibbs}
\pi_{2n}=\lambda_{V^{\pi}_{2n}}\quad \mbox{\rm and}\quad
\pi_{2n+1}=\lambda_{U^{\pi}_{2n+1}},
\end{equation}
with the potential functions
\begin{equation}\label{pi-gibbs-p}
V^{\pi}_{2n}:=V+(V_{2n}-V_{2n+1})
\quad \mbox{\rm and}\quad
U^{\pi}_{2n+1}:=U+(U_{2n+1}-U_{2(n+1)}).
\end{equation}

\subsection{General reference transitions}\label{general-sec}

Consider the following condition:\\

{\it
$(\Wa):$ The potential functions $W$ and $W^{\flat}$
satisfy  (\ref{lip-cond-H0}) for some $\kappa$. In addition, there exists positive definite matrices $\widehat{v}>0$ and $\widehat{u}>0$ such that the following holds 
$$
\nabla_1^2W(x,y_1)-\nabla_1^2W(x,y_2)\leq \widehat{u}^{-1}
\quad \mbox{and}\quad
\nabla_1^2W^{\flat}(y,x_1)-\nabla_1^2W^{\flat}(y,x_2)\leq \widehat{v}^{-1}.
$$}

\begin{rmk}
Assume $\widecheck{W}$ and $\overline{W}$ satisfy condition $(\Wa)$ for some parameters
$(\widecheck{\kappa},\widecheck{u},\widecheck{v})$ and $(\overline{\kappa},\overline{u},\overline{v})$. In this case, the function $W:=\widehat{W}+\overline{W}$
satisfy condition $(\Wa)$  with the parameters
$$
\kappa:=\widehat{\kappa}+\overline{\kappa}
\quad\mbox{and}\quad
(\widehat{u}^{-1},\widehat{v}^{-1}):=
(\widecheck{u}^{-1}+\overline{u}^{-1},\widecheck{v}^{-1}+\overline{v}^{-1}).
$$

\end{rmk}
A key feature of Sinkhorn potential functions $W_n$ is the following gradient formula.
\begin{lem}\label{nabla1-WN-lem}
For any $n\geq 0$ we have
\begin{eqnarray}
  \nabla_2W_{2n}(x_1,y)-  \nabla_2W_{2n}(x_2,y)&=&
   \nabla_2W(x_1,y)- \nabla_2W(x_2,y)\nonumber
\\
  \nabla_2W^{\flat}_{2n+1}(y_1,x)-  \nabla_2W^{\flat}_{2n+1}(y_2,x)&=&
   \nabla_2W^{\flat}(y_1,x)- \nabla_2W^{\flat}(y_2,x).
\label{nabla1-WN}
\end{eqnarray}
\end{lem}
A proof of this property is provided in Appendix~\ref{sec-grad-hess-a} (see Lemma~\ref{lem-appendix-Hess}).
This shows that
$W_{2n}$ and $W^{\flat}_{2n+1}$ satisfy (\ref{lip-cond-H0}) with $\kappa_n=\kappa$  as soon as  $W$ and $W^{\flat}$
satisfy  (\ref{lip-cond-H0}) for some $\kappa$. We conclude that all kernels $\Ka_n$ satisfy the 
Fisher-Lipschitz inequality  $\JJ_2(\kappa)$.

$\bullet$ Assume $C_V^-(U)$ is satisfied 
for some $0<u_+<\widehat{u}$ and $\overline{u}_+\in\Sa_d^0$ and $\delta_u\geq 0$. 
Using the Hessian formula presented in Lemma~\ref{lem-appendix-Hess} there exists some $\overline{u}\in \Sa_d^0 $ for any $n\geq 1$ we have the uniform estimates
\begin{eqnarray}
  \nabla^2_2W^{\flat}_{2n-1}(y,x)&\geq&\widetilde{u}^{-1}~1_{\Vert x\Vert_2 \geq  \delta_u}-\overline{u} 
~1_{\Vert x\Vert_2< \delta_u}\quad \mbox{\rm with}\quad \widetilde{u}^{-1}:=u_+^{-1}- \widehat{u}^{-1}>0.\label{loc-hess-ref-od}
\end{eqnarray}
This implies that $\Ka_{2n-1}$ satisfy the log-Sobolev inequality $LS(\widetilde{\rho}_u)$ and thus the $\TT_2(\widetilde{\rho}_u)$  inequality  for some parameter  $\widetilde{\rho}_u>0$ that depends on $(\widetilde{u},\overline{u},\delta_u)$. Recalling that $\lambda_V$ satisfies the $\TT_2(\Vert v_+\Vert)$-inequality we conclude that $A_V(\Sa)$ is satisfied with
\begin{equation}\label{c-vu-ref-eps}
\kappa_{2n-1}=\kappa
\quad \mbox{\rm and}\quad
\rho_{2n-1}=\widetilde{\rho}_u\Longrightarrow
 \epsilon_{2n-1}=\varepsilon_{\kappa}(\widetilde{\rho}_u,\Vert v_+\Vert_2).
\end{equation}

$\bullet$ Assume $C_U^-(V)$ is satisfied 
for some $0<v_+<\widehat{v}$ and $\overline{v}_+\in\Sa_d^0$ and $\delta_v\geq 0$. 
Using the Hessian formula presented in Lemma~\ref{lem-appendix-Hess} there exists some $\overline{v}\in \Sa_d^0 $ for any $n\geq 1$ we have the uniform estimates
\begin{eqnarray}
  \nabla^2_2W^{\flat}_{2n}(x,y)&\geq&\widetilde{v}^{-1}~1_{\Vert x\Vert_2 \geq  \delta_v}-\overline{u} 
~1_{\Vert x\Vert_2< \delta_v}\quad \mbox{\rm with}\quad \widetilde{v}^{-1}:=v_+^{-1}- \widehat{v}^{-1}>0.\label{loc-hess-ref-ev}
\end{eqnarray}
This implies that $\Ka_{2n}$ satisfy the log-Sobolev inequality $LS(\widetilde{\rho}_v)$ and thus the $\TT_2(\widetilde{\rho}_v)$  inequality  for some parameter  $\widetilde{\rho}_v>0$ that depends on $(\widetilde{v},\overline{v},\delta_v)$. Recalling that $\lambda_U$ satisfies the $\TT_2(\Vert u_+\Vert)$-inequality we conclude that $A_U(\Sa)$ is satisfied with
\begin{equation}\label{c-uv-ref-eps}
\kappa_{2n}=\kappa
\quad \mbox{\rm and}\quad
\rho_{2n}=\widetilde{\rho}_v
\Longrightarrow
 \epsilon_{2n}=\varepsilon_{\kappa}(\widetilde{\rho}_v,\Vert u_+\Vert_2).
\end{equation}

$\bullet$ Assume $C_0^-(U,V)$ is satisfied 
for some $0<u_+<\widehat{u}$ and  $0<v_+<\widehat{v}$. 
In this case, we recall that the probability measure $\lambda_U$  satisfies the $LS(\rho_u)$ and thus the $\TT_2(\rho_u)$ inequality; while  $\lambda_V$ satisfies the $LS(\rho_v)$ and thus the $\TT_2(\rho_v)$ inequality, with some parameters as in (\ref{rho-uv}).
In this situation conditions  $A_U(\Sa)$,$A_V(\Sa)$ and $\Aa_1(\Sa)$ hold with $\kappa_{2n}=\kappa_{2n-1}=\kappa$ and the 
 the time homogeneous parameters
\begin{equation}\label{c-uvu-ref-eps}
\begin{array}{l}
\rho_{2n}=\widetilde{\rho}_v=\overline{\rho}_{2n}
  \quad \mbox{and}\quad
  \rho_{2n-1}=\widetilde{\rho}_u=\overline{\rho}_{2n-1}\\
  \\
  \Longrightarrow
  (\epsilon_{2n-1},\epsilon_{2n})=(\varepsilon_{\kappa}(\widetilde{\rho}_u,\rho_v),\varepsilon_{\kappa}(\widetilde{\rho}_v,\rho_u))
   \quad \mbox{and}\quad
         (\overline{\epsilon}_{2n-1}, \overline{\epsilon}_{2n})=((\kappa \widetilde{\rho}_u)^2,
         (\kappa \widetilde{\rho}_v)^2).
   \end{array}
\end{equation}
Thus, condition $A_0^{\star}(\Sa)$ holds with
\begin{equation}\label{c-uvu-ref-eps-star}
\epsilon_{\star}:=\varepsilon_{\kappa}(\widetilde{\rho}_v,\rho_u)\vee \varepsilon_{\kappa}(\widetilde{\rho}_u,\rho_v).
\end{equation}

By (\ref{H01}) condition $A_2(\Sa)$ as well as  $A_U(\Sa)$ and $A_V(\Sa)$  hold with the parameters $\rho_n$ as soon as we have
\begin{equation}\label{Hess-gen}
\nabla^2_2W_{2n}\geq \rho_{2n}^{-1}~I\quad \mbox{and}\quad
\nabla^2_2W^{\flat}_{2n-1}\geq \rho_{2n-1}^{-1}~I.
\end{equation}
$\bullet$ When  $C_0(U,V)$ holds for some positive definite matrices $v_+< \widehat{v}$ and $u_+< \widehat{u}$ 
condition (\ref{Hess-gen}) also holds with the time homogeneous parameters
$$
\rho_{2n}=\Vert \widetilde{v}\Vert_2
  \quad \mbox{and}\quad
  \rho_{2n-1}=\Vert \widetilde{u}\Vert_2.
$$
In this case $A_0^{\star}(\Sa)$ is clearly met 
\begin{equation}\label{c-uvu-dif-ref-eps}
\epsilon_{\star}:=\varepsilon_{\kappa}(\Vert \widetilde{v}\Vert_2,\Vert u_+\Vert_2)\vee \varepsilon_{\kappa}(\Vert \widetilde{u}\Vert_2,\Vert v_+\Vert_2).
\end{equation}

\subsection{Contraction theorems}\label{contract-theo-sec}

This section is devoted to quantitative contraction estimates for Sinkhorn bridges
$(\Pa_n)_{n\geq 0}$. Building on the functional inequalities established in the previous sections, we show that suitable regularity properties of the underlying Sinkhorn transitions yield explicit entropy and Wasserstein contraction rates along the iteration.

Throughout this section, the parameters $\epsilon_n$ and $\overline{\epsilon}_n$ denote the collections defined in~\eqref{def-varepsi-SW}. These quantities encode the strength of the local logarithmic Sobolev and transport inequalities satisfied by the corresponding Sinkhorn transitions and play a central role in the contraction estimates derived below.

\subsubsection*{Some basic lemmas}

The next lemma is a consequence of (\ref{bridge-form-v2}) and the entropy estimate (\ref{intg-log-sob-2-ex-cor}), and the discussion provided in Section~\ref{general-sec}.
\begin{lem}\label{lem-1-odd}
Assume condition $A_V(\Sa)$ hold for some parameters $(\rho_{2n+1},\kappa_{2n+1})$.
 In this situation,
we have
\begin{equation}\label{bridge-form-v2-inside}
  \Ha(\mu~|~\pi_{2n+1})\leq\Ha(P_{\mu,\eta}~|~\Pa_{2n+1})\leq   \epsilon_{2n+1}~ \Ha(\pi_{2n}~|~\eta).
 \end{equation}
When $(\Wa)$ and $C_V^-(U)$ are satisfied for some $\kappa>0$ and $0<u_+<\widehat{u}$, 
the entropy estimate (\ref{bridge-form-v2-inside}) holds with the parameter $ \epsilon_{2n+1}=\varepsilon_{\kappa}(\widetilde{\rho}_u,\Vert v_+\Vert_2)$ and $
(\widetilde{\rho}_u,v_+)$ as in (\ref{c-vu-ref-eps}).
 \end{lem}
In the same vein, the next lemma is a consequence of (\ref{bridge-form-2-V2}) and the  estimate (\ref{intg-log-sob-2-ex-cor}).
 \begin{lem}\label{lem-1-even}
Assume condition $A_U(\Sa)$  holds for some parameters $(\rho_{2n},\kappa_{2n})$. In this situation,
 we have
 \begin{equation}\label{bridge-form-2-V2-inside}
  \Ha(\eta~|~\pi_{2n})\leq \Ha(P_{\mu,\eta}~|~\Pa_{2n})\leq 
\epsilon_{2n}~ \Ha(\pi_{2n-1}~|~\mu).
 \end{equation} 
 When $(\Wa)$ and $C_U^-(V)$ are satisfied for some $\kappa>0$ and  $0<v_+<\widehat{v}$, 
the entropy estimate (\ref{bridge-form-2-V2-inside}) holds with the parameter $ \epsilon_{2n}=\varepsilon_{\kappa}(\widetilde{\rho}_v,\Vert u_+\Vert_2)$ and $
(\widetilde{\rho}_v,u_+)$ as in (\ref{c-uv-ref-eps}).
 \end{lem}
Combining the entropy formulae (\ref{bridge-form-v2})  and
(\ref{bridge-form-2-V2}) with the $2$-Wasserstein estimates (\ref{intg-log-sob-2-cor}) we check the following lemma.
\begin{lem} \label{theo-sinkhorn-1}

Whenever $A_1(\Sa)$  is satisfied, we have the Wasserstein estimates
\begin{eqnarray}
\Da_{2}(\pi_{2n},\eta)^2&\leq&{2\overline{\rho}_{2n} ~
 \Ha( P_{\mu,\eta}~|~\Pa_{2n})}\leq \overline{\epsilon}_{2n}~\Da_2(\pi_{2n-1},\mu)^2.\label{d2-2n}\\
 &&\nonumber\\
 \Da_{2}(\pi_{2n+1},\mu)^2&\leq&{2\overline{\rho}_{2n+1}~
 \Ha( P_{\mu,\eta}~|~\Pa_{2n+1})}\leq  \overline{\epsilon}_{2n+1}~\Da_2(\pi_{2n},\eta)^2.\label{d2-2n1}
 \end{eqnarray}
  When $(\Wa)$ and
 $C_0^-(U,V)$ are satisfied for some $\kappa>0$, $0<u_+<\widehat{u}$ and  $0<v_+<\widehat{v}$, the estimates (\ref{d2-2n})
and (\ref{d2-2n1}) hold with   $(\overline{\epsilon}_{2n+1}, \overline{\epsilon}_{2n})=((\kappa \widetilde{\rho}_u)^2,
         (\kappa \widetilde{\rho}_v)^2)$ and the parameters $( \widetilde{\rho}_u, \widetilde{\rho}_v)$  as in (\ref{c-uvu-ref-eps}).
 \end{lem}
 
 \subsubsection*{Contraction inequalities}
The above lemmas can also be used to derive a variety of contraction estimates.
For instance, using (\ref{e12}) and (\ref{ent-r1}) we have the estimate
\begin{eqnarray*}
 \Ha(P_{\mu,\eta}~|~\Pa_{2n})&=&      \Ha(P_{\mu,\eta}~|~\Pa_{2(n-1)})- \left(  \Ha(\eta~|~\pi_{2(n-1)})+\Ha(\mu~|~\pi_{2n-1})\right)\\
 &\leq &  \Ha(P_{\mu,\eta}~|~\Pa_{2(n-1)})- \left(     \Ha(\pi_{2n-1}~|~\mu)+ \Ha(\pi_{2n}~|~\eta)\right).
\end{eqnarray*}
The entropy inequality (\ref{bridge-form-2-V2-inside}) now yields the rather crude contraction estimates.
\begin{theo}\label{Theo-f1-c}
Under the assumption of Lemma~\ref{lem-1-even},
 for any $n\geq 1$ we have 
\begin{equation}\label{f1-c}
 \Ha(P_{\mu,\eta}~|~\Pa_{2n})\leq \left(1+  \varepsilon^{-1}_{2n}
\right)^{-1}~ \Ha(P_{\mu,\eta}~|~\Pa_{2(n-1)}).
\end{equation}
 When $(\Wa)$ and $C_U^-(V)$ are satisfied for some $\kappa>0$ and  $0<v_+<\widehat{v}$, 
the entropy estimate (\ref{bridge-form-2-V2-inside}) holds with the parameter $ \epsilon_{2n}=\varepsilon_{\kappa}(\widetilde{\rho}_v,\Vert u_+\Vert_2)$ and $
(\widetilde{\rho}_v,u_+)$ as in (\ref{c-uv-ref-eps}).
\end{theo}

In the same vein, using (\ref{e12}) and (\ref{ent-r1}) we have the estimates
\begin{eqnarray*}
   \Ha(P_{\mu,\eta}~|~\Pa_{2n+1})&=& 
 \Ha(P_{\mu,\eta}|~\Pa_{2n-1})-\left(\Ha(\mu~|~\pi_{2n-1})+  \Ha(\eta~|~\pi_{2n})\right)\\
 &\leq & \Ha(P_{\mu,\eta}|~\Pa_{2n-1})-\left( \Ha(\pi_{2n}~|~\eta)+    \Ha(\pi_{2n+1}~|~\mu)\right).
\end{eqnarray*}
In this case, the entropy inequality (\ref{bridge-form-v2-inside}) yields the  following estimates
\begin{theo}\label{Theo-f2-c}
Under the assumption of Lemma~\ref{lem-1-odd},
 for any $n\geq 1$ we have
\begin{equation}\label{f2-c}
   \Ha(P_{\mu,\eta}~|~\Pa_{2n+1})\leq  \left(1+  \varepsilon^{-1}_{2n+1}\right)^{-1}  \Ha(P_{\mu,\eta}~|~\Pa_{2n-1}).
\end{equation}
When $(\Wa)$ and $C_V^-(U)$ are satisfied for some $\kappa>0$ and $0<u_+<\widehat{u}$, 
the entropy estimate (\ref{bridge-form-v2-inside}) holds with the parameter $ \epsilon_{2n+1}=\varepsilon_{\kappa}(\widetilde{\rho}_u,\Vert v_+\Vert_2)$ and $
(\widetilde{\rho}_u,v_+)$ as in (\ref{c-vu-ref-eps}).
\end{theo}

More refined estimates are presented in Section~\ref{sec-lingauss} in the context of 
linear Gaussian reference transitions.

\begin{theo}\label{theo-imp}
Assume that $A_0^{\star}(\Sa)$ is satisfied for some parameter $\epsilon_{\star}$. In this situation, for any  $n\geq 2 $ and $p\geq 0$   we have
\begin{equation}\label{theo-impp-a}
 \Ha(P_{\mu,\eta}~|~\Pa_{n+p})\leq
\left(1+\phi(\epsilon_{\star})\right)^{-(n-2)}~  \Ha(P_{\mu,\eta}~|~\Pa_{p}).
\end{equation}
with the function
\begin{equation}\label{def-phi-eps}
\phi(\epsilon_{\star}):=\frac{\epsilon_{\star}^{-1} }{\sqrt{1/4+ \epsilon_{\star}^{-1} }+1/2}.
\end{equation}
When $C_0^-(U,V)$ is satisfied 
for some $0<u_+<\widehat{u}$ and  $0<v_+<\widehat{v}$,  the estimates (\ref{theo-impp-a}) hold with $\epsilon_{\star}$ as in (\ref{c-uvu-ref-eps-star}).
\end{theo}
The proof of the above theorem is provided in the appendix on page~\pageref{theo-imp-proof}.

 Under condition  $A_0^{\star}(\Sa)$ the estimates (\ref{f1-c}) and (\ref{f2-c}) ensures that for any $n\geq 1$ and $p\geq 0$ we have
 $$
 \Ha(P_{\mu,\eta}~|~\Pa_{2n+p})\leq
\left(1+\epsilon^{-1}_{\star}\right)^{-n}~  \Ha(P_{\mu,\eta}~|~\Pa_{p}).
 $$
Replacing $n$ by $(2n)$ in  (\ref{theo-impp-a}) up to some constant  the exponential decays are  given by
$$
(1+\phi(\epsilon_{\star}))^{-2n}=\left(1+\epsilon^{-1}_{\star}+\phi(\epsilon_{\star})\right)^{-n}<\left(1+\epsilon^{-1}_{\star}\right)^{-n}.
$$
To check this claim, observe that
$$
\begin{array}{l}
\displaystyle 1+\phi(\epsilon_{\star})=\frac{1+\sqrt{1+4\epsilon^{-1}_{\star}}}{2}\Longleftrightarrow 1+2\phi(\epsilon_{\star})=\sqrt{1+4\epsilon^{-1}_{\star}}
\\
\\
\displaystyle\Longrightarrow
(1+\phi(\epsilon_{\star}))^2=\frac{1}{2}+\epsilon^{-1}_{\star}+\frac{1}{2}~\sqrt{1+4\epsilon^{-1}_{\star}}=
\frac{1}{2}+\epsilon^{-1}_{\star}+\frac{1}{2}~\left(1+2\phi(\epsilon_{\star})\right)\end{array}
$$
This yields the estimate
\begin{equation}\label{mrf}
(1+\phi(\epsilon_{\star}))^2=1+\epsilon^{-1}_{\star}+\phi(\epsilon_{\star})>
1+\epsilon^{-1}_{\star}.
\end{equation}

\section{Some illustrations}\label{sec-lingauss}
This section illustrates the general contraction and regularity results obtained in the previous sections in the concrete setting of linear Gaussian reference transitions
$K_W$ associated with some potential $W$ of the form (\ref{def-W}) for some parameters $(\alpha,\beta,\tau)$ and recall that $\cchi:=\tau^{-1}\beta$.

 Owing to their explicit structure, these models allow for precise computations of gradients, Hessians, and curvature quantities, making it possible to verify the assumptions introduced earlier and to obtain explicit bounds for the associated Sinkhorn transitions and Schrödinger bridges. In particular, we show how log-Sobolev, transport, and contraction inequalities can be derived directly from the parameters of the linear Gaussian transitions.

\subsection{Gaussian reference transitions}\label{sec-lingauss}
Using the gradient formulae (\ref{nabla1-WN}) (see also Example~\ref{l-exam-Vpi}) the linear-Gaussian structure of the reference potential is transferred to the one of Sinkhorn potentials; in the sense that   for any $n\geq 0$ we have
\begin{eqnarray}
  \nabla_2W_{2n}(x_1,y)-  \nabla_2W_{2n}(x_2,y)&=&
\cchi ~(x_2-x_1),\nonumber
\\
  \nabla_2W^{\flat}_{2n+1}(y_1,x)-  \nabla_2W^{\flat}_{2n+1}(y_2,x)&=&
 \cchi^{\prime}~(y_2-y_1). \label{nabla1-WN-ex}
\end{eqnarray}
This implies that
\begin{equation}\label{eq-2n}
\kappa_-^2~\Vert x_1-x_2\Vert_2^2\leq 
\jmath_{K_{W_{2n}}}(x_1,x_2)
=(x_2-x_1)^{\prime}\cchi\cchi^{\prime}(x_2-x_1)\leq 
\kappa^2~\Vert x_1-x_2\Vert_2^2,
\end{equation}
with the parameters
\begin{equation}\label{eq-2n-def-kappa}
\kappa_-=\sqrt{\ell_{\tiny min}(\cchi\cchi^{\prime})}
\leq \kappa:=\sqrt{\ell_{\tiny max}(\cchi\cchi^{\prime})}=\Vert\cchi\Vert_2.
\end{equation}
Recalling that $\cchi\cchi^{\prime}$ and $\cchi^{\prime}\cchi$ have the same eigenvalues,
we also check that
\begin{equation}\label{eq-2n1}
\kappa_-^2~\Vert y_1-y_2\Vert_2^2\leq
\jmath_{K_{W_{2n+1}^{\flat}}}
(y_1,y_2)\leq 
\kappa^2~\Vert y_1-y_2\Vert_2^2.
\end{equation}  
This yields for any 
$ \nu_1,\nu_2\in\Ma_1(\RR^d)$ and $n\geq 0$ the formulae
\begin{equation}\label{def-kappa}
 \kappa_-^2~\Da_2(\nu_1,\nu_2)^2  \leq
\Fa_{\Ka_{n}}(\nu_1|\nu_2)\leq \kappa^2~\Da_2(\nu_1,\nu_2)^2.  
\end{equation}
Without further work, this shows  that Sinkhorn transitions $\Ka_n$ satisfy 
Fisher-Lipschitz inequality  (\ref{lip-cond-H0-jmat}) and $(W_{2n},W^{\flat}_{2n+1})$ satisfies  (\ref{lip-cond-H0}) with a common parameter 
$$
\kappa_{2n}=\kappa=\kappa_{2n+1}.
$$
Using the Hessian formula presented in
Example~\ref{ref-hessian-W-f-examp} we have
$$
\nabla_1^2W^{\flat}(y,x_1)-\nabla_1^2W^{\flat}(y,x_2)=0
\quad \mbox{and}\quad
\nabla_1^2W(x,y_1)-\nabla_1^2W(x,y_2)=0.
$$
Thus condition $(\Wa)$ is met for {\it any} positive definite matrices $\widehat{v}>0$ and $\widehat{u}>0$.

\subsection{Log-Sobolev inequalities}\label{sec-lsob-s-kernels}
Recalling that $W_0=W$ we check that
$$
 (y_1-y_2)^{\prime} ( \nabla_2W_{0}(x,y_1)-  \nabla_2W_{0}(x,y_2))= (y_1-y_2)^{\prime}\tau^{-1}(y_1-y_2).
$$
In addition, for any $n\geq 1$ we also have the rather crude uniform estimates
 \begin{equation}\label{ref-transf}
 \begin{array}{l}
 (x_1-x_2)^{\prime}  \left(\nabla_2W^{\flat}_{2n-1}(y,x_1)-    \nabla_2W^{\flat}_{2n-1}(y,x_2)\right)\geq 
 (x_1-x_2)^{\prime} \left(\nabla U(x_1)-\nabla U(x_2)\right),\\
 \\
   (y_1-y_2)^{\prime}\left(   \nabla_2W_{2n}(x,y_1)-  \nabla_2W_{2n}(x,y_2)\right)\geq 
 (y_1-y_2)^{\prime}\left(\nabla V(y_1)-\nabla V(y_2) \right).
\end{array} 
\end{equation}
The proof of the above estimates follows elementary differential calculus, thus they are provided in Appendix~\ref{sec-grad-hess-a} (see Example~\ref{l-exam-DeltaUV}).

The curvature estimates (\ref{ref-transf})  can be used to transfer the convexity properties of the potential functions $(U,V)$ to the ones of Sinkhorn transitions. 

\begin{examp}
By (\ref{ref-transf})  when the potential function $U$ satisfies (\ref{ex-lg-cg}), all Sinkhorn transitions $\Ka_{2n-1}$ with odd indices satisfies  the log-Sobolev inequality. 
\end{examp}
\begin{examp}
Assume that the potential functions $(U,V)$ satisfies (\ref{ex-lg-cg}) for some possibly different parameters $(a_u,b_u)$ and $(a_v,b_v)$. In this situation, by (\ref{ref-transf}) for any $n\geq 1$ Sinkhorn transitions $\Ka_{2n-1}$ and $\Ka_{2n}$  also satisfy the log-Sobolev inequality $LS(\rho(a_u,b_u))$ and the
$LS(\rho(a_v,b_v))$  inequality  for some parameters $\rho(a_u,b_u)>0$ and $\rho(a_v,b_v)>0$. 

In this context, conditions $A_0(\Ka_n)$ and thus $\Aa_0(\Sa)$ are satisfied with $(\rho_0,\kappa_0)=(\Vert\tau\Vert_2,\kappa)$,  $(\rho_{2n},\kappa_{2n})=(\rho(a_v,b_v),\kappa)$ and 
$(\rho_{2n-1},\kappa_{2n-1})=(\rho(a_v,b_v),\kappa)$, for any $n\geq 1$.
\end{examp}

More refined curvature estimates can be obtained combining the covariance estimates stated in Proposition~\ref{prop-BL-CR} with the estimates (\ref{nbU}) and (\ref{nbV}). For instance, when (\ref{bary-ref-uv}) holds for some $v_->0$ and $u_+>0$ combining (\ref{nbU}) with (\ref{cov-ev}) for any $n\geq 0$ we have
 \begin{equation}\label{ref-transf-imp}
 \begin{array}{l}
 (x_1-x_2)^{\prime}  \left(\nabla_2W^{\flat}_{2n+1}(y,x_1)-    \nabla_2W^{\flat}_{2n+1}(y,x_2)\right)\\
 \\
\displaystyle\geq 
 (x_1-x_2)^{\prime} \left(\nabla U(x_1)-\nabla U(x_2)\right) + (x_1-x_2)^{\prime}~\cchi^{\prime}~v^{1/2}_-~\overline{\varsigma}_{2n}~v^{1/2}_-~\cchi~(x_1-x_2).
\end{array} 
\end{equation}
In the above display,  $\overline{\varsigma}_{2n}$ stands for the flow of Riccati matrices
$$
\overline{\varsigma}_{2n}=\mbox{\rm Ricc}_{\varpi_0}\left(\overline{\varsigma}_{2(n-1)}\right) \geq \mbox{\rm Ricc}_{\varpi_0}\left(0\right)=(I+\varpi_0^{-1})^{-1}=
(I+v^{1/2}_-~\left(\cchi~u_+~\cchi^{\prime}\right)~v^{1/2}_-)^{-1},
$$
with the initial condition $\overline{\varsigma}_{0}=v^{-1/2}_-~\tau~
v^{-1/2}_-$ (and  $\varpi_0$ defined in (\ref{def-over-varpi})). In this situation, we have
the rather crude uniform estimate
\begin{equation}\label{sc-reff}
 \begin{array}{l}
 (x_1-x_2)^{\prime}  \left(\nabla_2W^{\flat}_{2n+1}(y,x_1)-    \nabla_2W^{\flat}_{2n+1}(y,x_2)\right)\\
 \\
\displaystyle\geq 
 (x_1-x_2)^{\prime} \left[\left(\nabla U(x_1)-\nabla U(x_2)\right) + ~\cchi^{\prime}~(v^{-1}_-+\left(\cchi~u_+~\cchi^{\prime}\right))^{-1}~~\cchi\right]~(x_1-x_2).
\end{array} 
\end{equation}
More refined curvature estimates can be derived when $n\geq p$ using the inequalities (\ref{ricc-maps-incr}).

Last but not least, the estimates (\ref{def-kappa}) yield the equivalence between the 
local log-Sobolev inequality and the Lipschitz continuity of the 
the entropic transport cost associated with Sinkhorn transitions. More precisely, we have
$$
\begin{array}{l}
\mbox{\rm $\Ka_{n}$ satisfies  the  $LS_{\tiny loc}(\rho_{n})$ inequality}\quad
(\mbox{\rm with}\quad \rho_{n}=c_{n}/\kappa_-^{2})\\
\\
\displaystyle\Longleftrightarrow\quad \forall \nu_1,\nu_2\in \Ma_1(\RR^d)\qquad
\Ea_{\Ka_{n}}(\nu_1~|~\nu_2)\leq  \frac{c_{n}}{2}~\Da_2(\nu_1,\nu_2)^2\quad
(\mbox{\rm with}\quad  c_{n}=\kappa^2~\rho_{n}).
\end{array}$$
The last assertion is equivalent to the fact that for any $z_1,z_2\in\RR^d$ we have
$$
\Ha(\delta_{z_1}\Ka_n~|~\delta_{z_2}\Ka_n)\leq \frac{c_n}{2}~\Vert z_1-z_2\Vert_2^2.
$$

\subsection{Curvature inequalities}\label{sec-lsob-s-curvature}

\subsubsection*{Hessian-based sufficient conditions}
The conditions $A_2(\Ka_{n})$ 
are characterized in terms of lower bounds on the Hessians $\nabla_2^2W_{2n}$ and $\nabla^2_2 W^{\flat}_{2n+1}$.
Computing these matrices relies on standard, though somewhat lengthy, differential calculus for Sinkhorn transition potentials. More precisely, the Hessians of the log-densities \eqref{log-densities} can be expressed in terms of the Schrödinger bridge potentials \eqref{prop-schp}; see, for instance, \cite{chewi,chiarini,durmus}.
A complete derivation, together with explicit expressions in terms of conditional covariance matrices, is provided in the appendix (page~\pageref{cov-hess-sec}).

Concrete illustrations of these formulae in the linear-Gaussian setting can be found in
Examples~\ref{ex-gauss-UVn} and~\ref{l-exam-Vpi}.
In particular, Example~\ref{l-exam-Vpi} shows that
\begin{equation}\label{HessW0}
 \nabla^2_2W_0(x,y)=\tau^{-1},
\qquad
  \nabla_2^2W^{\flat}_{1}(y,x)
    =\nabla^2 U(x)+\cchi^{\prime}~\tau~\cchi.
\end{equation}
Moreover, for all $n\geq 1$,  one has the crude but robust bounds
\begin{equation}\label{HessWn}
  \nabla_2^2W_{2n}(x,y)\geq \nabla^2 V(y),\qquad
   \nabla_2^2W^{\flat}_{2n+1}(y,x)
    \geq\nabla^2 U(x),
\end{equation}
as well as their bridge counterparts
\begin{equation}\label{HessWnbridge}
  \nabla_2^2\WW_{\mu,\eta}(x,y)\geq \nabla^2 V(y),\qquad
   \nabla_2^2\WW^{\flat}_{\eta,\mu}(y,x)
    \geq\nabla^2 U(x).
\end{equation}

\subsubsection*{Implications for log-Sobolev regularity}

$\bullet$ Assume $C_U^-(V)$ holds for some positive definite matrix $v_+>0$ some $\overline{v}_+\in\Sa_d^0$ and $\delta\geq 0$. Then $V$ satisfies (\ref{ex-lg-cg-v2}) with
$a_v=\Vert v_+\Vert_2^{-1}>0$, $b_v=\Vert \overline{v}_+\Vert_2$.
Using (\ref{HessWn}), we obtain, for all $x,y\in\RR^d$ 
$$
  \nabla_2^2W_{2n}(x,y)\geq a_v~1_{\Vert y\Vert \geq  \delta}~I-b_v 
~1_{\Vert y\Vert< \delta}~I.
$$
By (\ref{HessWnbridge}) the same bound holds when $W_{2n}$ is replaced by $\WW_{\mu,\eta}$.
Consequently, the measure $\lambda_V$ as well as the Markov kernels $\Ka_{2n}$ and $K_{\WW_{\mu,\eta}}$  satisfy a log-Sobolev inequality $LS(\rho(a_v,b_v,\delta))$  for some   $\rho(a_v,b_v,\delta)>0$ depending only on $(a_v,b_v,\delta)$. 
In summary, 
\begin{equation}\label{v-cinfty}
C_U^-(V)\quad\Longrightarrow\quad
A_U(\Sa),
\end{equation} 
with the parameters   $\rho_{2n}=\rho(a_v,b_v,\delta)$ for all $n\geq 1$, $\rho_0=\Vert\tau\Vert_2$ and $\kappa_{2n}=\kappa$ as in (\ref{eq-2n}).

$\bullet$ Assume instead that $C_V^-(U)$ holds  for some positive definite matrix
$u_+>0$, some $\overline{u}_+\in\Sa_d^0$ and  $\delta\geq 0$. 
Then $U$ satisfies (\ref{ex-lg-cg-v2}) with
$a_u=\Vert u_+\Vert_2^{-1}>0$ and $b_u=\Vert \overline{u}_+\Vert_2$.
Arguing as above  and using (\ref{HessWn}) we deduce that
  $\lambda_U$, the kernels $\Ka_{2n+1}$ and 
  $K_{\WW^{\flat}_{\eta,\mu}}$   satisfy a log-Sobolev inequality $LS(\rho(a_u,b_u,\delta))$  for some   $\rho(a_u,b_u,\delta)>0$ depending only on $(a_u,b_u,\delta)$. 
Hence,
\begin{equation}\label{u-cinfty}
C_V^-(U)\quad\Longrightarrow\quad
A_V(\Sa).
\end{equation} 
 with parameters $\rho_{2n+1}=\rho(a_u,b_u,\delta)$ and $\kappa_{2n+1}=\kappa$ as in (\ref{eq-2n}).

\subsubsection*{Examples and remarks}

\begin{examp}
Assume (\ref{bary-ref-uv})
holds for some $v_->0$ and $u_+>0$. 
Then the curvature bound $  \nabla_2^2W^{\flat}_{2n+1}
 \geq u_+^{-1}$   implies that Sinkhorn transitions $\Ka_{W^{\flat}_{2n+1}}$  satisfy a $LS(\Vert u_+\Vert_2)$-inequality, even if $V$ is not convex. The stronger uniform condition $ \nabla^2 V\leq v^{-1}_-$ in (\ref{bary-ref-uv}) is only needed to derive sharper curvature estimates; see
 Remark~\ref{rmk-0v}.
 
\end{examp}

\begin{examp}
Assume again (\ref{bary-ref-uv})
for some $v_->0$ and $u_+>0$. This includes the case where $V$ is semiconvex, i.e.,  $ \nabla^2 V\geq - a~ I$ for some $a\geq 0$.  If, in addition the Aida-type condition (\ref{aida-condition}) holds for some $\delta>0$, then $\lambda_V$ satisfies a $\TT_2(\Vert v_+\Vert)$-inequality for some $v_+>0$ and  condition $C_V(U)$ is fulfilled. In this setting, the assumptions of Lemma~\ref{lem-1-odd} hold with $(\rho_{2n+1},\kappa_{2n+1})=(\Vert u_+\Vert_2,\kappa)$ and $\kappa$ as in (\ref{def-kappa}). 
\end{examp}

\begin{examp}
Assume (\ref{bary-ref-vu}) for some $u_->0$ and $v_+>0$. Then
the curvature bound $  \nabla_2^2W_{2n}
 \geq v_+^{-1}$   ensures that $\Ka_{W_{2n}}$  satisfies the $LS(\Vert v_+\Vert_2)$-inequality, even when $U$ is not  convex. As before, the uniform bound $ \nabla^2 U\leq u^{-1}_-$ in (\ref{bary-ref-vu})  is only required for
sharper estimates; see Remark~\ref{rmk-0v}).
 
\end{examp}

\begin{examp}
Finally, assume (\ref{bary-ref-vu})
holds for some $u_->0$ and $v_+>0$, corresponding to  a semiconvex potential $U$, i.e.,   $ \nabla^2 U\geq - a~ I$ for some $a\geq 0$. If (\ref{aida-condition}) holds for some $\delta>0$ then $\lambda_U$ satisfies a $\TT_2(\Vert u_+\Vert_2)$-inequality for some $u_+>0$ and  condition $C_U(V)$ holds. In this case, the assumptions of Lemma~\ref{lem-1-even} are satisfied with $(\rho_{2n},\kappa_{2n})=(\Vert v_+\Vert_2,\kappa)$ and $\kappa$ as in (\ref{def-kappa}). 
\end{examp}

\subsection{Some entropy inequalities}
  The contraction theorems presented in Section~\ref{contract-theo-sec}
  combined with the examples  above yields a variety of entropy estimates.
  In what follows $\varepsilon_{a}(b,c)$ stands for the collection of parameters  indexed by $a,b,c\geq 0$ defined in (\ref{A0-T2-def-vareps}).
  \subsubsection*{Convex at infinity potentials}
\begin{itemize}
\item Assume that $\mu=\lambda_U$ satisfies the $\TT_2(\Vert u_+\Vert_2)$  inequality for some $u_+\in\Sa_d^+$:
\begin{itemize}
\item When $\nabla^2 V\geq v_+^{-1}$, for any $n\geq 0$ conditions $A_0(\Ka_{2n})$ are met with $\rho_0=\Vert \tau_+\Vert_2$, and 
$
\rho_{2n}=\Vert v_+\Vert_2
$ for any $n\geq 1$. In addition $\eta=\lambda_V$ satisfies the $\TT_2(\Vert v_+\Vert_2)$ inequality. 
Conditions  $C_U(V)$ and $A_U(\Sa)$ are also satisfied with $\kappa$ as in (\ref{def-kappa}). In this case, by (\ref{bridge-form-2-V2-inside}) for any $n\geq 1$ we have the entropy estimate
\begin{equation}\label{f1-c-bis}
  \Ha(P_{\mu,\eta}~|~\Pa_{2n})\leq \varepsilon~
  \Ha(\pi_{2n-1}~|~\mu),
\end{equation}
with $\varepsilon$   defined in (\ref{def-vareps-intro-t}).
\item When the potential function $V$ satisfies (\ref{ex-lg-cg-v2}) for some $a_v>0$ and $b_v,\delta\geq 0$. By (\ref{v-cinfty}), conditions  $C_U^-(V)$ and $A_U(\Sa)$ are met with   $\rho_{2n}=\rho(a_v,b_v,\delta)$ for any $n\geq 1$, (note that $\rho_0=\Vert\tau\Vert_2$) and $\kappa_{2n}=\kappa$ as in (\ref{eq-2n}).
In this case, by (\ref{bridge-form-2-V2-inside}) and (\ref{f1-c}) for any $n\geq 1$ we have the entropy estimate
\begin{equation}\label{f1-c-uv}
\begin{array}{rcl}
  \Ha(P_{\mu,\eta}~|~\Pa_{2n})&\leq& \varepsilon_{u,v}~
  \Ha(\pi_{2n-1}~|~\mu)\\
  &&\\
   \Ha(P_{\mu,\eta}~|~\Pa_{2n})&\leq& \left(1+  \varepsilon^{-1}_{u,v}
\right)^{-1}~ \Ha(P_{\mu,\eta}~|~\Pa_{2(n-1)}), 
\end{array}
\end{equation}
with the parameter
\begin{equation}\label{epsi-uv-ref-g}
  \varepsilon_{u,v}:=\varepsilon_{\kappa}(\rho(a_v,b_v,\delta),\Vert u_+\Vert_2)
  \quad \mbox{\rm with}\quad
   \kappa=\Vert\tau^{-1}\beta\Vert_2.
\end{equation}
\end{itemize}
\item  Assume that $\eta=\lambda_V$ satisfies the $\TT_2(\Vert v_+\Vert_2)$  inequality for some $v_+\in\Sa_d^+$:
\begin{itemize}
\item When $\nabla^2 U\geq u_+^{-1}$, for any $n\geq 0$ condition $A_0(\Ka_{2n+1})$ are met with
$
\rho_{2n+1}=\Vert u_+\Vert
$ and $\mu$ satisfies the $\TT_2(\Vert u_+\Vert_2)$ inequality. 
Conditions  $C_V(U)$ and $A_V(\Sa)$ are also satisfied with $\kappa$ as in (\ref{def-kappa}). In this case, by (\ref{bridge-form-v2-inside}) we have the entropy estimate
\begin{equation}\label{f2-c-bis}
  \Ha(P_{\mu,\eta}~|~\Pa_{2n+1})\leq ~\varepsilon~
~ \Ha(\pi_{2n}~|~\eta).
\end{equation}
 \item When the potential function $U$ satisfies (\ref{ex-lg-cg-v2}) for some $a_u>0$ and $b_u,\delta\geq 0$. By (\ref{u-cinfty}), conditions  $C_V^-(U)$ and $A_V(\Sa)$  are met with   $\rho_{2n+1}=\rho(a_u,b_u,\delta)$ for any $n\geq 0$ and $\kappa_{2n}=\kappa$ as in (\ref{eq-2n}).
In this case, by (\ref{bridge-form-v2-inside}) and (\ref{f2-c}) for any $n\geq 1$ we have the entropy estimate
\begin{equation}\label{f1-c-uv-i}
\begin{array}{rcl}
  \Ha(P_{\mu,\eta}~|~\Pa_{2n+1})&\leq& ~\varepsilon_{v,u}~
~ \Ha(\pi_{2n}~|~\eta) \\
  &&\\
   \Ha(P_{\mu,\eta}~|~\Pa_{2n+1})&\leq & \left(1+  \varepsilon^{-1}_{v,u}\right)^{-1}  \Ha(P_{\mu,\eta}~|~\Pa_{2n-1}),
\end{array}
\end{equation}
with the parameter
\begin{equation}\label{epsi-vu-ref-g}
  \varepsilon_{v,u}:=\varepsilon_{\kappa}(\rho(a_u,b_u,\delta),\Vert v_+\Vert_2)
    \quad \mbox{\rm with}\quad
   \kappa=\Vert\tau^{-1}\beta\Vert_2.
\end{equation}

\end{itemize}
\item Assume $C_0^-(U,V)$ is satisfied for some $u_+,v_+\in\Sa_d^+$ and $\overline{u}_+,\overline{v}_+\in\Sa_d^0$ and $\delta_u,\delta_v\geq 0$. 
In this situation, for any  $n\geq 2 $ and $p\geq 0$   we have
\begin{equation}\label{theo-impp-a-kg}
 \Ha(P_{\mu,\eta}~|~\Pa_{n+p})\leq
\left(1+\phi(\epsilon_{\star})\right)^{-(n-2)}~  \Ha(P_{\mu,\eta}~|~\Pa_{p}),
\end{equation}
with $\phi$ as in (\ref{def-phi-eps}). In this case condition $A_0^{\star}(\Sa)$ holds with
\begin{equation}\label{c-uvu-ref-eps-star-g}
\epsilon_{\star}:=\varepsilon_{\kappa}(\rho_u,\rho_v)\quad\mbox{\rm with $(\rho_u,\rho_v)$ as in (\ref{rho-uv})}    \quad \mbox{\rm and}\quad
   \kappa=\Vert\tau^{-1}\beta\Vert_2.
\end{equation}

\end{itemize}
  \subsubsection*{Strongly convex potentials}

 We further assume that condition $C_0(U,V)$ is satisfied for some  $u_+,v_+>0$. 
By (\ref{HessW0}), conditions  $A_1(\Ka_0)$ and $A_2(\Ka_0)$ hold with
\begin{equation}\label{defrho0}
\overline{\rho}_0=\rho_0:=\Vert\tau\Vert_2
  \quad \mbox{and $\kappa$ as in (\ref{def-kappa})}.
\end{equation}

Working a little harder, as shown in Example~\ref{l-exam-Vpi}, the probability measure $$\pi_0=\mu \Ka_{0}=\lambda_{V^{\pi}_0}$$ is associated with a potential function $V^{\pi}_0$ with the Hessian matrix  
\begin{equation}\label{pi0-t2}
  \nabla^2V^{\pi}_{0}\geq \frac{1}{\rho^{\pi}_0}~I\quad \mbox{\rm with the parameter}\quad
 \rho^{\pi}_0:= \Vert \tau+\beta~ u_+~\beta^{\prime}\Vert_2\geq \rho_0.
\end{equation}
This property ensures that $\pi_0$ satisfies the $\TT_2(\rho^{\pi}_0)$ inequality as well as the $LS(\rho^{\pi}_0)$ inequality.
\begin{rmk}
Choosing $(\tau,\beta)$ such that
$$
u_+^{-1}> \beta^{\prime}\tau^{-1}v_+\tau^{-1}\beta\quad \mbox{and}\quad
v_+^{-1}>\tau^{-1}\beta~u_+
 \beta^{\prime}\tau^{-1},
$$
there exists some $\rho_1^{\pi},\rho_2^{\pi}>0$ such that for any $n\geq 1$ we have
$$
  \nabla^2U^{\pi}_{2n-1}\geq \frac{1}{\rho^{\pi}_{1}}~I
\quad \mbox{and}\quad
  \nabla^2V^{\pi}_{2n}\geq \frac{1}{\rho^{\pi}_{2}}~I.
$$
These assertions are checked using the Hessian formulae provided in Appendix~\ref{sec-a-stm}
(see Remark~\ref{rmk-TLS-sd}).
In this context, Sinkhorn distributions $\pi_{2n-1}$ satisfy the  $LS(\rho^{\pi}_1)$ and $\TT_2(\rho^{\pi}_1)$ inequalities and $\pi_{2n}$  satisfies the  $LS(\rho^{\pi}_2)$ and $\TT_2(\rho^{\pi}_2)$ inequalities.
\end{rmk}
We also have
\begin{equation}\label{c0a0}
C_0(U,V)\Longrightarrow
A_0^{\star}(\Sa)\quad \mbox{\rm and}\quad 
A_2(\Sa).
\end{equation}
More precisely, using (\ref{HessW0}) we also check that conditions $A_1(\Ka_1)$  and $A_2(\Ka_1)$ are met with the parameters
\begin{equation}\label{defrho1}
\overline{\rho}_1=\rho_1:=1/\ell_{\tiny min}\left(u_+^{-1}+
     \beta^{\prime}\tau^{-1}\beta
\right)=\Vert (u_+^{-1}+
     \beta^{\prime}\tau^{-1}\beta)^{-1}\Vert_2\leq \Vert u_+\Vert_2.
\end{equation}
By (\ref{HessWn}), for any $n\geq 1$ conditions  $A_2(\Ka_{2n+1})$ and respectively  $A_2(\Ka_{2n})$ are also satisfied with $\rho=\Vert u_+\Vert_2$ and respectively  $\rho=\Vert v_+\Vert_2$.
The above discussion  also shows that condition $A_0^{\star}(\Sa)$ is satisfied. 
 In this context, the parameters  $ \epsilon_n$ and $\epsilon_{\star}$  introduced in (\ref{def-varepsi-SW})  and (\ref{epsilon-star}) resume  to 
 $$
  \epsilon_0=\kappa^2~\Vert\tau\Vert_2~\Vert u_+\Vert_2\quad \mbox{\rm and}\quad
   \epsilon_1=\kappa^2~\Vert (u_+^{-1}+
     \beta^{\prime}\tau^{-1}\beta)^{-1}\Vert_2~\Vert v_+\Vert_2.
 $$
In addition, in terms of the parameter $\varepsilon$ defined in (\ref{def-vareps-intro-t}) for any $n\geq 1$ we have
\begin{equation}\label{defepsilon12}
  \epsilon_{n+1}=\epsilon_2=\epsilon_{\star}=\varepsilon\geq 
   \epsilon_1.
\end{equation}
In summary, we have
$$
C_U(V)\Longrightarrow A_U(\Sa),
$$
 with  $\rho_0=\Vert \tau\Vert_2$ and $\kappa_0=\kappa$  as in (\ref{def-kappa}) and $(\rho_{2n},\kappa_{2n})=(\Vert v_+\Vert_2,\kappa)$,  for any $n\geq 1$.
 
We also  have
$$
C_V(U)\Longrightarrow A_V(\Sa),
$$
with 
$\rho_{1}$ as in  (\ref{defrho1}) and $\kappa_1=\kappa$ as in  (\ref{def-kappa}) 
and $(\rho_{2n+1},\kappa_{2n+1})=(\Vert u_+\Vert_2,\kappa)$ for any $n\geq 1$.

\begin{examp}\label{examp-tau-t}
Consider the case $\beta=I$ and $\tau=tI$ for some $t>0$ and diagonal matrices $u_+,v_+$. In this situation we have
$$
\rho_0=t
\qquad
\kappa:=1/t\qquad
 \rho_1:=\frac{t}{t+\Vert u_+\Vert_2}~\Vert u_+\Vert_2\quad \mbox{and}\quad
  \rho^{\pi}_0:= t+\Vert u_+\Vert_2.
$$
We also have
 $$
 \epsilon_1=\frac{1}{t}~\frac{\Vert u_+\Vert_2\Vert v_+\Vert_2}{t+\Vert u_+\Vert_2}~\leq  \epsilon_2=\varepsilon=\frac{1}{t^2}~\Vert u_+\Vert_2\Vert v_+\Vert_2\longrightarrow_{t\rightarrow\infty}~0.
 $$
\end{examp} 

\section{Some continuity theorems}

This section specializes the general contraction framework developed earlier to linear-Gaussian reference transitions and derives a collection of quantitative continuity results for Schrödinger bridges, Sinkhorn iterates, and associated entropic maps.
Exploiting the explicit structure of Gaussian kernels, we obtain sharp entropy and Wasserstein contraction inequalities, refined exponential decay rates, and continuity estimates for barycentric projections, conditional covariances, and proximal sampler semigroups.

We focus on linear-Gaussian reference transitionss $K_W$ associated with a potential $W$ of the form (\ref{def-W}) parametrized by $(\alpha,\beta,\tau)$ and recall the notation $\cchi:=\tau^{-1}\beta$.

Throughout this section, the parameters $(\varepsilon,  \varepsilon_{u,v},  \varepsilon_{v,u},
\epsilon_{\star})$ denote the collections defined in (\ref{def-vareps-intro-t}), (\ref{epsi-uv-ref-g}) (\ref{epsi-vu-ref-g}) and (\ref{c-uvu-ref-eps-star-g}), respectively, while 
$\phi$ refers to the function defined in (\ref{def-phi-eps}).

\subsection{Contraction inequalities}\label{sct}

Using the entropy estimates (\ref{f1-c-bis}) and (\ref{f2-c-bis}) and following word-for-word the proof of (\ref{f1-c}) and (\ref{f2-c}) we check the following theorem. 
\begin{theo}\label{th-1-in}
\begin{itemize}
\item Assume $C^-_U(V)$. In this situation, for any $n\geq 1$ and $p\geq 0$ we have
\begin{equation}\label{f1-c-t}
 \Ha(P_{\mu,\eta}~|~\Pa_{2(n+p)})\leq \left(1+ \varepsilon_{u,v}^{-1}
\right)^{-n}~ \Ha(P_{\mu,\eta}~|~\Pa_{2p}). 
\end{equation}
When $C_U(V)$ is satisfied, the above estimate holds replacing  $\varepsilon_{u,v}$ by $\varepsilon$.
\item Assume $C_V^-(U)$. In this situation, for any $n\geq 1$ and $p\geq 0$ we have
\begin{equation}\label{f2-c-t}
   \Ha(P_{\mu,\eta}~|~\Pa_{2(n+p)+1})\leq  \left(1+\varepsilon_{v,u}^{-1}
\right)^{-n}~\Ha(P_{\mu,\eta}~|~\Pa_{2p+1}).
\end{equation}
When $C_V(U)$ is satisfied, the above estimate holds replacing  $\varepsilon_{v,u}$ by $\varepsilon$.
\end{itemize}

\end{theo}

Theorem~\ref{th-1-in} can also be obtained as a direct corollary of Theorems~\ref{Theo-f1-c} and~\ref{Theo-f2-c} when applied to linear--Gaussian reference transitions.

Next theorem is a direct consequence of (\ref{defepsilon12}) and Theorem~\ref{theo-imp}. 
 \begin{theo}\label{theo-thc-s}
 Assume $C_0^-(U,V)$. 
 In this situation, for any  $n\geq 2 $ and $p\geq 0$   we have
\begin{equation}\label{expo-decays-imp}
 \Ha(P_{\mu,\eta}~|~\Pa_{n+p})\leq
\left(1+\phi(\varepsilon_{\star})\right)^{-(n-2)}~  \Ha(P_{\mu,\eta}~|~\Pa_{p}).
\end{equation}
When $C_0(U,V)$ is satisfied, the above estimate holds replacing  $\varepsilon_{\star}$ by $\varepsilon$.

 \end{theo}

\begin{cor}\label{cor-ref-i}
Under the assumptions of Theorem~\ref{theo-thc-s},
for any $p\geq 0$ and $n\geq 1$ we have the exponential decays
 \begin{eqnarray}
   \Ha(\eta~|~\pi_{2(n+p)})&\leq &\varepsilon_{\star}
\left(1+\phi(\varepsilon_{\star})\right)^{-2(n-1)}~ \Ha(\eta~|~\pi_{2p}),\label{thc-cor-s}
\\
\Da_{2}(\eta,\pi_{2(n+p)})
& \leq&\varepsilon_{\star}
\left(1+\phi(\varepsilon_{\star})\right)^{-(n-1)}~ \Da_2(\eta,\pi_{2p}).\label{thc-cor-s-2}
 \end{eqnarray}
  The same estimate holds replacing $(\eta,\pi_{2(n+p)},\pi_{2p})$ by $(\mu,\pi_{2(n+p)+1},\pi_{2p+1})$. Moreover, when $C_0(U,V)$ is satisfied, the above estimates holds replacing  $\varepsilon_{\star}$ by $\varepsilon$.
\end{cor}
 \proof
 We check (\ref{thc-cor-s}) using (\ref{expo-decays-imp}), (\ref{bridge-form-2-V2-inside}) and (\ref{e12}).
 More precisely (\ref{expo-decays-imp}) and (\ref{bridge-form-2-V2-inside}) yield for any $p\geq 1$ the estimate
 $$
    \Ha(\eta~|~\pi_{2(n+p)})\leq \Ha(P_{\mu,\eta}~|~\Pa_{2(n+p)})\leq 
\left(1+\phi(\varepsilon_{\star})\right)^{-2(n-1)}~\varepsilon_{\star}~
  \Ha(\pi_{2p-1}~|~\mu).
 $$
The end of the proof of (\ref{thc-cor-s}) is now a direct consequence of (\ref{e12}).

 Using (\ref{d2-2n}) and (\ref{expo-decays-imp})
 we also have
 $$
\Da_{2}(\pi_{2(n+(p+1))},\eta)^2
 \leq 2\Vert v_+\Vert_2 ~
\left(1+\phi(\varepsilon_{\star})\right)^{-2(n-1)}~ \Ha(P_{\mu,\eta}~|~\Pa_{2(p+1)}). 
 $$
By (\ref{d2-2n}) and (\ref{d2-2n1}) we also have 
 $$
  2\Vert v_+\Vert_2 ~
 \Ha( P_{\mu,\eta}~|~\Pa_{2(p+1)})\leq \varepsilon_{\star}~\Da_2(\pi_{2p+1},\mu)^2\leq
\varepsilon_{\star}^2~\Da_2(\pi_{2p},\eta)^2.
  $$
We conclude that
$$
\Da_{2}(\pi_{2(n+(p+1))},\eta)^2
 \leq 
\varepsilon_{\star}^2~
\left(1+\phi(\varepsilon_{\star})\right)^{-2(n-1)}~ ~\Da_2(\pi_{2p},\eta)^2.
$$
This ends the proof of (\ref{thc-cor-s-2}). The proof of the estimates for odd indices follows word-for-word the same lines of arguments, thus it is skipped.
This ends the proof of the corollary.
\cqfd

\subsection{Proximal sampler semigroup}\label{proxi-sec}

Recall that $S_1=\Ka_{0}\Ka_{1}$ is the Markov transition of the simple alternating forward/backward Gibbs sampler (a.k.a. proximal sampler) with target measure
$$
\Pa_0(d(x,y))=\mu(dx)\Ka_{0}(x,dy)= \exp{\left(-U(x)-W(x,y)\right)}~\lambda(dx)\lambda(dy).
$$
In what follows $\varepsilon_{a}(b,c)$ stands for the collection of parameters  indexed by $a,b,c\geq 0$ defined in (\ref{A0-T2-def-vareps}).
By (\ref{pi0-t2}), 
 the marginal $\pi_0$ of $\Pa_0$ with respect to the second coordinate  satisfies the  $\TT_2(\rho^{\pi}_0)$ with $\rho^{\pi}_0$ defined  in (\ref{pi0-t2}). Also recall that $A_1(\Ka_0)$ is met with the parameter $\overline{\rho}_0=\rho_0=\Vert\tau\Vert_2$ defined in (\ref{defrho0}), and
 $A_1(\Ka_1)$ is met with $\overline{\rho}_1=\rho_1$ defined in (\ref{defrho1}). 
 
 Observe that
 $$
  \epsilon_0=\varepsilon_{\kappa}(\rho_0,u_+)=a_{u_+}(\beta,\tau)
 \quad \mbox{\rm and}\quad
 \varepsilon_{\kappa}(\rho_0,\rho_1)= b_{u_+}(\beta,\tau),
 $$
  with  the parameters $ \epsilon_0$, $\kappa$  and $(a_{u_+}(\beta,\tau), b_{u_+}(\beta,\tau))$ defined in (\ref{def-varepsi-SW}) (\ref{def-kappa})  and Theorem~\ref{theo-S1n}.
Also note that
 $$
\varepsilon_{\kappa}(\rho_0,\rho_1)\leq     \epsilon_0\wedge    \epsilon^{\pi}_0 
\quad\mbox{\rm with}\quad 
   \epsilon^{\pi}_0:=\varepsilon_{\kappa}(\rho_0^{\pi},\rho_1).
  $$
 \begin{examp}\label{examp-tau-t-2}
In the example discussed in (\ref{examp-tau-t}) we have
 $$
   \epsilon^{\pi}_0=\frac{\Vert u_+\Vert_2}{t}
=  \epsilon_0 \longrightarrow_{t\rightarrow\infty}~0
\quad\mbox{
as well as}\quad
 b_{u_+}(\beta,\tau)=\frac{\Vert u_+\Vert_2 }{t+\Vert u_+\Vert_2}<1.
 $$

 \end{examp}
\begin{theo}\label{theo-proxi}
For any $\nu\in\Ma_1(\RR^d)$ we have the estimates
\begin{eqnarray}
\Da_{2}(\nu S_{1},\mu)&\leq & b_{u_+}(\beta,\tau)~\Da_2(\nu ,\mu),\label{d2-prox}
\\
\Ha(\nu S_{1}~|~\mu)
&\leq&~  \left(1\wedge
   \epsilon^{\pi}_0
\right) \left(1\wedge   \epsilon_0 \right)~~\Ha(\nu~|~\mu).\label{ent-prox}
\end{eqnarray}
In addition, for any $n\geq 0$ we have the entropy estimates
\begin{equation}\label{ent-prox-vv2}
\Ha(\nu S_{1}^{n+1}~|~\mu )~\leq ~ a_{u_+}(\beta,\tau) ~
 b_{u_+}(\beta,\tau)^{2n}~
\Ha(\nu~|~\mu).
\end{equation}
\end{theo}

\proof

Recalling that $(\overline{\rho}_0,\overline{\rho}_1)=(\rho_0,\rho_1)$ and using twice (\ref{lipD2})  we check that
\begin{eqnarray*}
\Da_{2}(\nu S_{1},\mu)
&\leq& \sqrt{\varepsilon_{\kappa}(\rho_1,\rho_1)}~\Da_2(\nu \Ka_{0},\mu \Ka_{0})\leq    \varepsilon_{\kappa}(\rho_0,\rho_1)
~\Da_2(\nu ,\mu).
\end{eqnarray*}
The last assertion comes from the fact that
$$
\sqrt{\varepsilon_{\kappa}(\rho_1,\rho_1)}\sqrt{\varepsilon_{\kappa}(\rho_0,\rho_0)}
=(\kappa~\rho_0)\times(\kappa~\rho_1)=\kappa^2~\rho_0~\rho_1=\varepsilon_{\kappa}(\rho_0,\rho_1).
$$
This ends the proof of (\ref{d2-prox}).
Consider a  probability distribution $\nu_2$  satisfying a $\TT_2(\upsilon)$ inequality.
By (\ref{lip-ent}) for any probability distribution $\nu_1$ we check that
$$
\Ha(\nu_1 \Ka_{1}~|~\nu_2 \Ka_{1})\leq~\left(1\wedge
\varepsilon_{\kappa}(\nu,\rho_1)
\right)~ \Ha(\nu_1~|~\nu_2).
$$

Choosing $$\nu_2=\mu \Ka_{0}=\pi_{0}\quad\mbox{\rm and}\quad \nu_1=\nu \Ka_{0},$$
and recalling that  $\pi_{0}$ satisfy a
$\TT_2(\rho^{\pi}_0)$-inequality
 this implies that
$$
\Ha(\nu S_{1}~|~\mu)\leq~ \left(1\wedge \varepsilon_{\kappa}(\rho^{\pi}_0,\rho_1)
\right)~ \Ha(\nu \Ka_{0}~|~\mu \Ka_{0}).
$$
Recalling that $\mu=\lambda_U$ satisfies a $\TT_2(\Vert u_+\Vert_2)$-inequality, applying (\ref{lip-ent}) and using (\ref{def-kappa}) we have
$$
\Ha(\nu \Ka_{0}~|~\mu \Ka_{0})\leq~ \left(1\wedge
\varepsilon_{\kappa}(\rho_{0},\Vert u_+\Vert_2)
\right)~ \Ha(\nu~|~\mu).
$$
We conclude that
$$
\begin{array}{l}
\Ha(\nu S_{1}~|~\mu)
\leq~  \left(1\wedge
\varepsilon_{\kappa}(\rho^{\pi}_0,\rho_1)
\right) \left(1\wedge \varepsilon_{\kappa}(\rho_{0},\Vert u_+\Vert_2)\right)~~\Ha(\nu~|~\mu).
\end{array}$$
Using (\ref{d2-prox}) we also check that
$$
\Da_{2}(\nu S_{1}^n,\mu)^2\leq   2~  b_{u_+}(\beta,\tau)^{2n}~\frac{1}{2}~\Da_2(\nu ,\mu)^2\leq 
2 ~  b_{u_+}(\beta,\tau)^{2n}~
\Vert u_+\Vert_2~
\Ha(\nu~|~\mu).
$$
The r.h.s. estimate comes from the fact that $\mu=\lambda_U$ satisfies a $\TT_2(\Vert u_+\Vert_2)$-inequality. On the other hand, using  (\ref{intg-log-sob-v2}) we check that
\begin{eqnarray*}
\Ha(\nu S_{1}^{n+1}~|~\mu )&=&\Ha((\nu S_{1}^{n}\Ka_{0})
\Ka_{1}~|~(\mu \Ka_{0})
\Ka_{1})\\
&\leq&\Ha(\nu S_{1}^{n}\Ka_{0}
~|~\mu \Ka_{0})\leq
 \frac{1}{2}~  \varepsilon_{\kappa}(\rho_{0},1)~\Da_{2}(\nu S_{1}^{n},\mu)^2.
\end{eqnarray*}
We conclude that
$$
\Ha(\nu S_{1}^{n+1}~|~\mu )~\leq ~ \varepsilon_{\kappa}(\rho_{0},\Vert u_+\Vert_2)
~ (\kappa^2~\rho_0~\rho_1)^{2n}~
\Ha(\nu~|~\mu).
$$
This ends the proof of the proposition.
\cqfd

 \begin{examp}\label{examp-tau-t-22}
Consider the model discussed in (\ref{examp-tau-t}) and Example~\ref{examp-tau-t-2}.
In this context, to compare the exponential converge decays of Sinkhorn iterates (\ref{thc-cor-s}) and the ones of proximal sampler (\ref{ent-prox-vv2}), observe that
 $$
(1+ \varepsilon^{-1})^{-1}=\frac{\Vert u_+\Vert_2\Vert v_+\Vert_2}{t^2+\Vert u_+\Vert_2\Vert v_+\Vert_2}<
 b_{u_+}(\beta,\tau)^2
  \Longleftrightarrow \frac{t}{2}~\left(\Vert v_+\Vert_2^{-1}-\Vert u_+\Vert_2^{-1}\right)>0.
$$
with the parameter  $\varepsilon$ defined in (\ref{def-vareps-intro-t}).
 We conclude that
 $$
 \Vert u_+\Vert_2>\Vert v_+\Vert_2\Longleftrightarrow
(1+ \varepsilon^{-1})^{-1}<  b_{u_+}(\beta,\tau)^2.
 $$

 \end{examp}

\subsection{Log-Lyapunov estimates}\label{scvx-sec}
By (\ref{ent-r1}) the sequences 
$ \Ha(P_{\mu,\eta}~|~\Pa_{2n})$ and $   \Ha(P_{\mu,\eta}~|~\Pa_{2n+1})$ are decreasing
and we have
$$
 \Ha(P_{\mu,\eta}~|~\Pa_{2(n+1)})\leq \Ha(P_{\mu,\eta}|~\Pa_{2n+1})\leq     \Ha(P_{\mu,\eta}~|~\Pa_{2n}).
$$
By Theorem~\ref{theo-thc-s} they converge exponentially fast to $0$ as $n\rightarrow\infty$. The aim of this section is to estimate the largest Lyapunov exponent
\begin{eqnarray*}
\Lambda_{\mu,\eta}&:=&\limsup_{n\rightarrow\infty}\frac{1}{n}\log{\Ha(P_{\mu,\eta}~|~\Pa_{n})}\leq
-\log{
(1+\phi(\epsilon_{\star}))
}. 
\end{eqnarray*}

Consider the
linear-Gaussian reference transitions $K_W$ associated with some potential $W$ of the form (\ref{def-W}) for some parameters $(\alpha,\beta,\tau)$.

Without further mention, in this section we assume that the potential functions $(U,V)$
satisfy  $C_1(U,V)$ for some positive definite matrices $u_+,v_+>0$ and some positive semi-definite matrices  $u_{-},v_{-}\in\Sa_d^+\cup\{0\}$.
In what follows $\varepsilon$ also stands for the parameter  defined in (\ref{def-vareps-intro-t}).

Consider the Riccati matrix difference equations
\begin{equation}\label{ricc-interm}
\overline{\tau}_{2n}=\mbox{\rm Ricc}_{\overline{\varpi}_0}\left(\overline{\tau}_{2(n-1)}\right)  \quad\mbox{and}\quad
 \overline{\tau}_{2n+1} =\mbox{\rm Ricc}_{\overline{\varpi}_1}\left(\overline{\tau}_{2n-1}\right),
\end{equation}
with the positive definite matrices $(\overline{\varpi}_0,\overline{\varpi}_1)$ defined in (\ref{def-over-varpi})
and the initial conditions 
$$
\overline{\tau}_{0}=v^{-1/2}_+~\tau~v^{-1/2}_+
\quad\mbox{\rm and}\quad
\overline{\tau}_{1} =u_+^{-1/2}
\left(~u_+^{-1}+~\cchi^{\prime}~\tau~\cchi \right)^{-1}~u_+^{-1/2}.
  $$
  Next lemma is a consequence of the monotone properties of Riccati maps.
      \begin{lem}\label{lem-HW}
For any $n\geq p\geq 1$ we have the estimates
$$
 \overline{\tau}_{2(n+1)}\leq \widehat{\tau}_{2p}:=\mbox{\rm Ricc}_{\overline{\varpi}_0}^{p}(I)
\quad \mbox{and}\quad
\overline{\tau}_{2n+1}\leq  \widehat{\tau}_{2p+1}:=\mbox{\rm Ricc}_{\overline{\varpi}_1}^p\left( I\right).
$$
  \end{lem}
  The detailed proof of the above lemma is provided in Appendix~\ref{sec-a-stm} on page~\pageref{lem-HW-proof}.
  Consider 
  the rescaled matrices
  $$
  {\tau}_{2n}:=v^{1/2}_+~\overline{\tau}_{2n}~v^{1/2}_+\quad\mbox{and}\quad
  {\tau}_{2n+1}:=u_+^{1/2}~\overline{\tau}_{2n+1}~u_+^{1/2}.
  $$
With this notation at hand we have the following Hessian estimates.
\begin{prop}\label{prop-pre-Lyap}
For any  $n\geq 0$ we have
  $$
   \nabla_2^2W_{2n}\geq {\tau}_{2n}^{-1}\quad \mbox{and}\quad
     \nabla_2^2W^{\flat}_{2n+1}\geq {\tau}_{2n+1}^{-1}.
     $$
  \end{prop}
  The proof of the above proposition is based on successive applications of Brascamp-Lieb and Cramer Rao inequalities. The detailed proof is provided in Appendix~\ref{sec-a-stm}, see for instance Lemma~\ref{prop-BL-CR} and formulae (\ref{tmde}). 
  
Recalling that $r_{\overline{\varpi}_0}\leq I$ and  $r_{\overline{\varpi}_1}\leq I$ we have
  $$
   r_{\overline{\varpi}_0}\leq
  \widehat{\tau}_{2(p+1)}\leq \widehat{\tau}_{2p}\leq I\quad \mbox{\rm and}\quad
    r_{\overline{\varpi}_1}\leq 
    \widehat{\tau}_{2p+1}\leq \widehat{\tau}_{2p-1}\leq I.
 $$
We also have the monotone properties
\begin{equation}\label{def-iota-v}
1\leq \xi_{2p}:=\frac{\Vert v_+\Vert_2}{\Vert v^{1/2}_+~\widehat{\tau}_{2p}~v^{1/2}_+\Vert_2}\leq \xi_{2(p+1)}\leq \iota_{0}(\overline{\varpi}_0):=
\frac{\Vert v_+\Vert_2}{\Vert v^{1/2}_+~ r_{\overline{\varpi}_0}~v^{1/2}_+\Vert_2},
\end{equation}
as well as
\begin{equation}\label{def-iota-u}
1\leq \xi_{2p-1}:=\frac{\Vert u_+\Vert_2}{\Vert u^{1/2}_+~\widehat{\tau}_{2p-1}~u^{1/2}_+\Vert_2}\leq  \xi_{2p+1}
\leq \iota_{1}(\overline{\varpi}_1):=\frac{\Vert u_+\Vert_2}{\Vert u^{1/2}_+~  r_{\overline{\varpi}_1}~u^{1/2}_+\Vert_2}.
\end{equation}

Next theorem is a consequence of Proposition~\ref{prop-pre-Lyap} combined with the estimates (\ref{f1-c}) and (\ref{f2-c}).

  \begin{theo}\label{theo-impp}
  For $n\geq p\geq 1$ we have the exponential estimates
\begin{equation}\label{ineq-theo-impp}
\begin{array}{l}
 \displaystyle \Ha(P_{\mu,\eta}~|~\Pa_{2(n+1)})
 \leq  \left(1+ (\xi_{2p}\vee \xi_{2p+1})~\varepsilon^{-1}\right)^{-(n-p)}  \Ha(P_{\mu,\eta}~|~\Pa_{2p-1}).
\end{array}\end{equation}
with  $\varepsilon$  as in (\ref{defepsilon12}) and  $(\xi_{2p},\xi_{2p+1})$ defined in (\ref{def-iota-v}) and (\ref{def-iota-u}).
  \end{theo}
    The proof of the above theorem is provided in Appendix~\ref{sec-tech-proofs} on page~\pageref{theo-impp-proof}.
    Using (\ref{cv-ricc-intro}) we check the following corollary.
\begin{cor}\label{cor-impp}
  For $n\geq p\geq 1$ we have the exponential estimates
\begin{equation}\label{ineq-theo-impp-2}
\begin{array}{l}
 \displaystyle \Ha(P_{\mu,\eta}~|~\Pa_{2(n+1)})\\
 \\
 \leq  \left(1+   \varepsilon^{-1}~ \iota(\overline{\varpi})~\left(1+
c_{\overline{\varpi}}~\delta_{\overline{\varpi}}^p\right)^{-1}\right)^{-(n-p)}  \Ha(P_{\mu,\eta}~|~\Pa_{2p-1}).
\end{array}\end{equation}
with the parameters
$$
\iota(\overline{\varpi}):=\vee_{i\in\{ 0,1\}}
\iota_{i}(\overline{\varpi}_i)\qquad
c_{\overline{\varpi}}:=\vee_{i\in\{ 0,1\}} c_{\overline{\varpi}_i}
\quad\mbox{and}\quad
\delta_{\overline{\varpi}}:=\vee_{i\in\{ 0,1\}}\delta_{\overline{\varpi}_i}.
$$
In the above display $( c_{\overline{\varpi}_i},\delta_{\overline{\varpi}_i},\iota_{i}(\overline{\varpi}_i))$ with $i\in\{0,1\}$ stands for the parameters defined in
(\ref{cv-ricc-intro}), (\ref{def-iota-v}), 
(\ref{def-iota-u}).
  \end{cor}
    The proof of the above corollary is provided in Appendix~\ref{sec-tech-proofs} on page~\pageref{cor-impp-proof}.
\begin{rmk}
  \begin{itemize}
\item When $u_-=0$ and $v_->0$ for any $n\geq 1$ we use the convention
 $
 \overline{\varpi}_0^{-1}=0=\delta_{\overline{\varpi}_0}=c_{\overline{\varpi}_0}
$ and 
$ \overline{\tau}_{2n}= \widehat{\tau}_{2n}=I=r_{\overline{\varpi}_0}
 $.
 In this case we have $\overline{\varpi}_1>0$ and $\widehat{\tau}_{2p+1}<I$ and therefore
 $$
\xi_{2p}=1=\iota_{0}(\overline{\varpi}_0)\quad \mbox{and}\quad  1<\xi_{2p+1}< \iota_{1}(\overline{\varpi}_1).
 $$
 \item When $v_-=0$ and $u_->0$ for any $n\geq 1$ we use the convention
 $
\overline{\varpi}_1^{-1}=0=\delta_{\overline{\varpi}_1}=c_{\overline{\varpi}_1}$ and $
 \overline{\tau}_{2n+1}=   \widehat{\tau}_{2n+1}=I=r_{\overline{\varpi}_1}
 $.  In this case we have $\overline{\varpi}_0>0$ and and $\widehat{\tau}_{2p}<I$ therefore
 $$
1<\xi_{2p}<\iota_{0}(\overline{\varpi}_0)\quad \mbox{and}\quad  \xi_{2p+1}=1=\iota_{1}(\overline{\varpi}_1).
 $$
 \end{itemize}
 \end{rmk}

Letting $n\rightarrow\infty$ and then $p\rightarrow\infty$ in (\ref{ineq-theo-impp-2}) we check the following corollary. 
\begin{cor}
We have
 the log-Lyapunov estimate
\begin{equation}\label{log-lyap}
\Lambda_{\mu,\eta}
\leq -\frac{1}{2}
\log{\left(1+\iota(\overline{\varpi})~\varepsilon^{-1}\right)},
\end{equation}
with $\varepsilon$ as in  (\ref{defepsilon12}) and $\iota(\overline{\varpi})$ defined in (\ref{ineq-theo-impp-2}).
 \end{cor}

Recalling that the square root is concave we have $(u+v)^{1/2}\geq (u^{1/2}+v^{1/2})/\sqrt{2}$
for any positive definite matrices $u,v>0$. This yields the estimate
\begin{equation}\label{inq-fpr}
\frac{1}{\sqrt{2}}~\left( \varpi^{1/2}-\frac{\varpi}{2(1+\sqrt{2})}\right) \leq r_{\varpi}:=\left(\varpi+\left(\frac{\varpi}{2}\right)^2\right)^{1/2}-\frac{\varpi}{2}.
\end{equation}
Using (\ref{square-root-key-estimate}) we check that
\begin{equation}\label{inq-sqrtfp}
\ell_{\tiny max}(r_{\varpi}) \leq \frac{\ell_{\tiny max}(\varpi)}{(\lambda^{1/2}_{min}(\varpi)\vee (\lambda_{min}(\varpi)/2))+
\lambda_{min}(\varpi)/2}~\leq \frac{\ell_{\tiny max}(\varpi)}{\lambda^{1/2}_{min}(\varpi)}.
\end{equation}
\begin{examp}\label{examp-ref-tau0}
Consider the situation $\tau=t~\tau_0$ for some $\tau_0>0$ and some scalar $t>0$ and set $\cchi_0:=\tau_0^{-1}\beta$. In this case we have
$$
\varepsilon:=\frac{1}{t^2}~\widehat{\varepsilon}\quad \mbox{with the parameter}\quad 
\widehat{\varepsilon}:=\Vert\cchi_0\Vert^2~\Vert u_+\Vert_2~\Vert v_+\Vert_2
$$
We also have
$
  \overline{\varpi}_i=t^2~ \widehat{\varpi}_i
$
with $i\in\{0,1\}$ and the matrices $\widehat{\varpi}_i$ defined in (\ref{def-over-varpi-check}).
  Using (\ref{inq-sqrtfp}), for any $i\in\{0,1\}$ we check that
  $
\ell_{\tiny max}(r_{  \overline{\varpi}_i}) ~\leq t~{\ell_{\tiny max}(\widehat{\varpi}_i)}/{\ell_{\tiny min}^{1/2}(\widehat{\varpi}_i)}
  $ so that
  $$
 \iota_{0}(\overline{\varpi}_0)\geq 
\frac{1}{t}~\frac{\ell_{\tiny min}^{1/2}(\widehat{\varpi}_0)}{\ell_{\tiny max}(\widehat{\varpi}_0)}
\quad \mbox{and}\quad \iota_{1}(\overline{\varpi}_1)\geq \frac{1}{t}~\frac{\ell_{\tiny min}^{1/2}(\widehat{\varpi}_1)}{\ell_{\tiny max}(\widehat{\varpi}_1)}.
  $$
  We conclude that
$$
1+\varepsilon^{-1}~\iota(\overline{\varpi})\geq 1+t~\widehat{\varepsilon}^{\,-1}~\left(
\frac{\ell_{\tiny min}^{1/2}(\widehat{\varpi}_0)}{\ell_{\tiny max}(\widehat{\varpi}_0)}\vee
\frac{\ell_{\tiny min}^{1/2}(\widehat{\varpi}_1)}{\ell_{\tiny max}(\widehat{\varpi}_1)}\right).
$$
On the other hand, we have
  \begin{eqnarray*}
 (1+\phi(\varepsilon))^2
&=&1+
t^2~\widehat{\varepsilon}^{\,-1}~\left(1+\frac{2}{\sqrt{1+4 t^2~\widehat{\varepsilon}^{\,-1}}+1}\right)\leq 
1+2
t^2~\widehat{\varepsilon}^{\,-1}.
\end{eqnarray*}
Choosing
$$
2t\leq 
\frac{\ell_{\tiny min}^{1/2}(\widehat{\varpi}_0)}{\ell_{\tiny max}(\widehat{\varpi}_0)}
\vee
\frac{\ell_{\tiny min}^{1/2}(\widehat{\varpi}_1)}{\ell_{\tiny max}(\widehat{\varpi}_1)},
$$
we check that
$$
 (1+\phi(\varepsilon))^2\leq  1+t~\widehat{\varepsilon}^{\,-1}~\left( \frac{\ell_{\tiny min}^{1/2}(\widehat{\varpi}_0)}{\ell_{\tiny max}(\widehat{\varpi}_0)}\vee
\frac{\ell_{\tiny min}^{1/2}(\widehat{\varpi}_1)}{\ell_{\tiny max}(\widehat{\varpi}_1)}\right)\leq 1+\varepsilon^{-1}~\iota(\overline{\varpi}).
$$
 \end{examp}

\subsection{Entropic maps}\label{entropy-maps-sec}
By Theorem 4.2 in~\cite{nutz} Schr\"odinger bridges takes the form (\ref{sbr-ref})
with the potential functions $\WW_{\mu,\eta}$ and $\WW^{\flat}_{\eta,\mu}$  defined by
\begin{equation}
  U(x)+\WW_{\mu,\eta}(x,y)=\UU(x)+W(x,y)+\VV(y)=  V(y)+\WW^{\flat}_{\eta,\mu}(y,x).\label{log-densities-fp}
\end{equation}
In the above display $(\UU,\VV)$ stands for the solution of the Schr\" odinger system
\begin{equation}
\UU=U+\log{K_W(e^{-\VV})}\quad \mbox{\rm and}\quad
\VV=V+\log{K_{W^{\flat}}(e^{-\UU})}.
\label{prop-schp-fp}
\end{equation}
Note that $(\UU,\VV)$ is the fixed point of the integral equations discussed in (\ref{prop-schp-fp}) and the potential functions $(\WW_{\mu,\eta},\WW^{\flat}_{\eta,\mu},)$ are defined as $(W_{2n},W_{2n+1}^{\flat})$ by replacing in  in (\ref{log-densities}) the potential functions $(U_{n},V_n)$ by $(\UU,\VV)$.
 A description of the barycentric projections (\ref{bary-ref})  in terms of the entropic potential functions $(\UU,\VV)$ is also given in (\ref{ref-nablaUVng-bg}).

We further assume that the potential functions $(U,V)$
satisfy the following estimates
\begin{equation}\label{allc}
u_+^{-1}\leq \nabla^2 U\leq u^{-1}_-\quad\mbox{\rm and}\quad
v_+^{-1}\leq \nabla^2 V\leq v^{-1}_-,
\end{equation}
 for  some positive semi-definite matrices  $u_{-},v_{-}\in\Sa_d^+\cup\{0\}$ and $u_{+}^{-1},v_{+}^{-1}\in\Sa_d^+\cup\{0\}$.
 When $s=0$ we use the convention $s^{-1}=\infty\times I$. 
 
 We underline the following situations:
\begin{itemize}
\item When $u_-=0$ and $v_->0$ and  $u_+^{-1}>0$ condition (\ref{allc}) takes the form (\ref{convention-uv}).
In addition, when $v_+^{-1}=0$  condition (\ref{allc}) reduces to (\ref{bary-ref-uv}).
\item  When $v_-=0$ and $u_->0$ and $v_+^{-1}>0$ condition (\ref{allc}) takes the form
(\ref{convention-vu}).
In addition, when $u_+^{-1}=0$ condition (\ref{allc}) reduces to (\ref{bary-ref-vu}).

\end{itemize}

\begin{theo}\label{last-est-th}
For any $x\in\RR^d$ we have the estimates
\begin{equation}\label{last-est}
   \begin{array}{rcccl}
\cchi^{\prime}~ v_-^{1/2} ~r_{\varpi_0}~v_-^{1/2}~\cchi&\leq &\nabla K_{\WW_{\mu,\eta}}(I)(x)~\cchi&\leq&
\cchi^{\prime}~v_+^{1/2} ~r_{\overline{\varpi}_0}~v_+^{1/2}~\cchi,\\
&&&&\\
\cchi~u_-^{1/2}~r_{\varpi_1}~u_-^{1/2}~\cchi^{\prime}&\leq& \nabla K_{\WW^{\flat}_{\eta,\mu}}(I)(x)~\cchi^{\prime}&\leq&
\cchi~u_+^{1/2}~r_{\overline{\varpi}_1}~u_+^{1/2}~\cchi^{\prime} .
\end{array}
\end{equation}

\end{theo}
The proof of the above theorem is provided in Appendix~\ref{sec-a-stm} on page~\pageref{last-est-th-proof}.

 We underline the following situations:
\begin{itemize}
\item When $u_-=0$ and $v_->0$ and  $u_+^{-1}>0$, by (\ref{convention-u}) the estimate   (\ref{last-est}) reduces to
$$
   \begin{array}{rcccl}
\cchi^{\prime}~ v_-^{1/2} ~r_{\varpi_0}~v_-^{1/2}~\cchi&\leq &\nabla K_{\WW_{\mu,\eta}}(I)~\cchi&\leq&
\cchi^{\prime}~v_+~\cchi,\\
&&&&\\
&& \nabla K_{\WW^{\flat}_{\eta,\mu}}(I)~\cchi^{\prime}&\leq&
\cchi~u_+^{1/2}~r_{\overline{\varpi}_1}~u_+^{1/2}~\cchi^{\prime}.
\end{array}
$$
In addition, when $v_+^{-1}=0$ we have the estimate (\ref{ch-uv}).
\item When $v_-=0$ and $u_->0$ and $v_+>0$, by (\ref{convention-v}) the estimate (\ref{last-est}) reduces to
$$
   \begin{array}{rcccl}
& &\nabla K_{\WW_{\mu,\eta}}(I)~\cchi&\leq&
\cchi^{\prime}~v_+^{1/2} ~r_{\overline{\varpi}_0}~v_+^{1/2}~\cchi,\\
&&&&\\
\cchi~u_-^{1/2}~r_{\varpi_1}~u_-^{1/2}~\cchi^{\prime}&\leq& \nabla K_{\WW^{\flat}_{\eta,\mu}}(I)~\cchi^{\prime}&\leq&
\cchi~u_+~\cchi^{\prime}.
\end{array}
$$
In addition, when $u_+^{-1}=0$ we have the estimate (\ref{ch-vu}).
\end{itemize}

As in (\ref{sc-reff}) the convexity properties of the bridge
potential functions $(\WW_{\mu,\eta},\WW^{\flat}_{\eta,\mu})$ can be estimated in term of the ones of the potential functions $(U,V)$. For instance (see Example~\ref{l-exam-Vlspi}) we have
$$
 \begin{array}{l}
 (x_1-x_2)^{\prime}  \left(\nabla_2W^{\flat}_{\eta,\mu}(y,x_1)-    \nabla_2W^{\flat}_{\eta,\mu}(y,x_2)\right)\geq 
 (x_1-x_2)^{\prime} \left(\nabla U(x_1)-\nabla U(x_2)\right).
\end{array} 
$$
As in (\ref{ref-transf}) the above estimate can be used to transfer the convexity properties of the potential functions $U$ to the ones of bridge transition. For instance, when the potential function $U$ satisfies (\ref{ex-lg-cg}) the bridge transitions $K_{W^{\flat}_{\eta,\mu}}$  satisfies  the log-Sobolev inequality. 
In addition, when the more restrictive condition (\ref{bary-ref-uv}) holds we have
\begin{equation}\label{reffin}
 \begin{array}{l}
 (x_1-x_2)^{\prime}  \left(\nabla_2W^{\flat}_{\eta,\mu}(y,x_1)-    \nabla_2W^{\flat}_{\eta,\mu}(y,x_2)\right)\\
 \\
\displaystyle\geq 
 (x_1-x_2)^{\prime} \left(\nabla U(x_1)-\nabla U(x_2)\right) + (x_1-x_2)^{\prime}~\cchi^{\prime}~v_-^{1/2} ~r_{\varpi_0}~v_-^{1/2}~\cchi~(x_1-x_2).
\end{array} \end{equation}
Following the discussion given in Section~\ref{sec-lingauss} the above estimate is based on condition
(\ref{bary-ref-uv}) which requires  the potential function $U$ to be strongly convex.
In this situation, by (\ref{nabla2-WN-ex-b}) (and without assuming the uniform upper bound $\nabla^2 V\leq v_-^{-1}$) we already have the uniform estimate $\nabla_2^2W^{\flat}_{\eta,\mu}\geq u_+^{-1}$. Thus the strong convexity of $U$ ensures directly  that $K_{W^{\flat}_{\eta,\mu}}$ satisfies the log-Sobolev inequality $LS(\Vert u_+\Vert)$. 

\subsection{Sinkhorn maps}\label{sinkmaps-sec}

Sinkhorn maps $(\Za_{2n},\Za_{2n+1})$ are the random maps defined by
\begin{equation}\label{sinkhorn-maps}
\Ka_{2n}(x,dy)=\PP(\Za_{2n}(x)\in dy)\quad\mbox{and}\quad
\Ka_{2n+1}(x,dy)=\PP(\Za_{2n+1}(x)\in dy).
\end{equation}
For clarity of exposition, we write $(\Sigma_{2n},\Sigma_{2n+1})$ instead of 
$
(\Sigma_{W_{2n}},\Sigma_{W^{\flat}_{2n+1}})$ the conditional covariances (\ref{cov-mat-W}) of the Sinkhorn maps  associated with Sinkhorn transitions
$
 (\Ka_{2n},
  \Ka_{2n+1})=(K_{W_{2n}},K_{W^{\flat}_{2n+1}})
$.

In terms of the entropic maps (\ref{e-maps}), by (\ref{bridge-form-intro}) and (\ref{sbr-ref}) we also have 
$$
\begin{array}{l}
\Za_{2n}=\ZZ_{\mu,\pi_{2n}}\quad \mbox{\rm and}\quad
\Za_{2n+1}=\ZZ^{\flat}_{\eta,\pi_{2n+1}}\\
\\
\Longrightarrow\quad \Sigma_{2n}=\Sigma_{\mu,\pi_{2n}}
\quad \mbox{\rm and}\quad \Sigma_{2n+1}=\Sigma^{\flat}_{\eta,\pi_{2n+1}}.
\end{array}
$$
A closed form expression of Sinkhorn maps for  Gaussian distributions (\ref{def-U-V})  is provided in Appendix~\ref{sec-a-stm} on page~\pageref{lin-gauss-sec} (see also~\cite{adm-24}).  

Following (\ref{bary-ref}), the barycentric projections (a.k.a. the entropic maps) of  Sinkhorn bridges are given by 
\begin{equation}\label{bary-ref-sink}
\Ka_{2n}(I)(x)=\EE(\Za_{2n}(x))
\quad\mbox{\rm and}\quad
\Ka_{2n+1}(I)(x)=\EE(\Za_{2n+1}(x)).
\end{equation}

By Theorem~\ref{theo-baryp}
we have the conditional covariance formulae
\begin{equation}\label{conditional-cov-sink}
\cchi^{\prime}~\Sigma_{2n}~\cchi=\nabla\Ka_{2n}(I)~\cchi\quad\mbox{and}\quad
\cchi~\Sigma_{2n+1}~\cchi^{\prime}=\nabla\Ka_{2n+1}(I)~\cchi^{\prime}.
\end{equation}
An alternative proof of (\ref{conditional-cov-sink}) based on elementary differential calculus is provided in Appendix~\ref{sec-grad-hess-a}, see formula (\ref{finex-n}) in Example~\ref{ex-gauss-UVn}. As in Section~\ref{entropy-maps-sec}, replacing the fixed point equations (\ref{log-densities-fp}) by Sinkhorn recursions (\ref{prop-schp})
several Sinkhorn maps estimates can be obtained in terms of Riccati matrix difference equations (see Proposition~\ref{prop-BL-CR} in Appendix~\ref{sec-a-stm}). 
\subsection*{Notes and references}\label{comparisons-1}

Continuity properties of Sinkhorn conditional means and conditional covariances, expressed in terms of the relative entropy between Sinkhorn and Schrödinger bridges, are investigated in~\cite{dp-ent-26}. Related perspectives based on stochastic control representations of entropic optimal transport, as well as backward propagation techniques along Hamilton–Jacobi equations, can be found in~\cite{chiarini, conforti-ptrf, greco-25}.

When specialized to the model considered in Example~\ref{examp-tau-t} and choosing $\rho=\|v_+\|_2$, the estimate (ii) in Proposition~1.3 of~\cite{chiarini} coincides, for sufficiently small values of $t$, with the exponential decay rate obtained in~(\ref{f2-c-t}). In contrast, as $t\to\infty$, the decay rate in part (i) of Theorem~1.2 in~\cite{chiarini} converges to $1/2$, whereas our estimates~(\ref{f1-c-t}) and~(\ref{f2-c-t}) vanish in the same limit, as shown in Section~\ref{reg-transp-sec}.

We further emphasize that the exponential decay estimates in~(\ref{expo-decays-imp}), established in Theorem~\ref{theo-thc-s}, are strictly sharper than those given in~(\ref{f1-c-t}) and~(\ref{f2-c-t}). In particular, they rely on an improved convergence rate~(\ref{def-vareps-intro}); see also the comparison in~(\ref{mrf}).

More refined contraction estimates, such as part (i) of Proposition~1.3 in~\cite{chiarini}, remain valid under the additional assumption $\nabla^2 V\leq c_v I$ for some constant $c_v>0$. In that case, the resulting decay rates depend explicitly on $c_v$ and deteriorate as $c_v\to\infty$. When the Hessians of $U$ or $V$ are uniformly bounded, sharper estimates can be derived using Brascamp–Lieb and Cramér–Rao inequalities. These situations are discussed in detail in Section~\ref{scvx-sec}, which is devoted to strongly convex models; see also~\cite{adm-24} for related results in the linear–Gaussian setting.

Finally, in the setting of Example~\ref{examp-tau-t}, the Wasserstein contraction estimate~(\ref{d2-prox}) coincides with the bound established in Lemma~2 of~\cite{lee}. Similarly, the entropy estimates~(\ref{ent-prox-vv2}) are closely related to Lemma~4 in~\cite{lee} and to Theorem~6 in~\cite{guan}, which are derived in the context of proximal samplers on Riemannian manifolds.

\section*{Acknowledgments}

I would like to thank the anonymous reviewers for their excellent suggestions for improving the paper. Their detailed comments greatly improved the presentation of the article.

\appendix

\section{Covariance inequalities}\label{BL-CR-ineq-sec}
Consider a Markov transition $K_W$ with a log-density function $W$ such that
\begin{equation}\label{strong-convex}
\tau^{-1}\leq \nabla_2^2W(x,y)\leq\varsigma^{-1}\quad \mbox{\rm for some $\tau,\varsigma\in \Sa_d^+$.}
\end{equation} 
  By the strong convexity condition (\ref{strong-convex}), the Brascamp-Lieb inequality ensures that
$$
  \Sigma_{W}(x)\leq \int K_W(x,dy)~\left(\nabla^2_2W(x,y)\right)^{-1}\leq \tau,
  $$
  with the conditional covariance function  $  \Sigma_{W}$ defined in (\ref{cov-mat-W}).

In this context,  the Cramer-Rao inequality also ensures that
$$
  \Sigma_{W}(x)^{-1}\leq \int K_W(x,dy)~\nabla^2_2W(x,y)\leq\varsigma^{-1}.
$$
We summarize the above discussion with  the covariance estimate
\begin{equation}\label{bl-cr}
(\ref{strong-convex})\Longrightarrow
\varsigma\leq  \Sigma_{W}(x)\leq \tau.
\end{equation}

\section{A review on Sinkhorn semigroups}\label{review-sinkhorn}

In this section, we discuss some of the theory behind Sinkhorn semigroups and present some results that are of use throughout the paper. This section is mainly taken from~\cite{adm-24}.\\

\subsection*{Conjugate formulae}
Sinkhorn transitions  $(\Ka_{2n},\Ka_{2n+1})$ can alternatively be defined sequentially by the conjugate (a.k.a. forward-backward) formulae 
\begin{equation}
\pi_{2n}\times \Ka_{2n+1}=(\mu\times \Ka_{2n})^{\flat}\quad \mbox{\rm and}\quad
(\pi_{2n+1}\times \Ka_{2(n+1)})^{\flat}=\eta\times \Ka_{2n+1}.
\label{s-2}
\end{equation}
with the flow of distributions $\pi_n$ introduced in (\ref{def-pin}).  Using the transport equations (\ref{def-pin})  the entropic transport problems (\ref{sinhorn-entropy-form}) take the following form
$$
\Pa_{2n+1}= \argmin_{Q\in \Pi(\pi_{2n+1},\eta)}\Ha(Q~|~\mu\times\Ka_{2n})
\quad \mbox{\rm and}\quad
\Pa_{2(n+1)}= \argmin_{Q\in  \Pi(\mu,\pi_{2(n+1)})}\Ha(Q~|~(\eta\times\Ka_{2n+1})^{\flat}).
$$
In terms of the potential functions $(U_n,V_n)$ defined in (\ref{prop-schp}), for any $n\geq 0$ we have
$$
\Pa_{n}(d(x,y))=e^{-U_n(x)}~e^{-W(x,y)}~e^{-V_n(y)}~\lambda^{\otimes 2}(d(x,y)).
$$
Recalling that  Sinkhorn algorithm starts from the reference distribution
$\Pa_0=P=\mu\times K_W$, for $P$-a.e.  we have
$$
\frac{d\Pa_{2n+1}}{dP}(x,y)=e^{-(U_{2n}-U)(x)}~e^{-V_{2n}(y)}
\quad \mbox{\rm and}\quad
\frac{d\Pa_{2(n+1)}}{dP}(x,y)=e^{-(U_{2(n+1)}-U)(x)}~e^{-V_{2(n+1)}(y)}.
$$
Applying Theorem 2.1 in~\cite{nutz} or more simply by  (\ref{def-entropy-pb-v2-2}), this implies that
$$
\Pa_{2n+1}= \argmin_{Q\in \Pi(\pi_{2n+1},\eta)}\Ha(Q~|~P)
\quad \mbox{\rm and}\quad
\Pa_{2n}=\argmin_{Q\in  \Pi(\mu,\pi_{2n})}\Ha(Q~|~P).
$$
This yields the bridge formulae (\ref{bridge-form-intro}). Applying  (\ref{ref-back-bridge}) we also have
\begin{eqnarray*}
\eta\times\Ka_{2n+1}&=&\Pa_{2n+1}^{\flat}
=
(P_{\pi_{2n+1},\eta})^{\flat}=P^{\flat}_{\eta,\pi_{2n+1}},\\
\mu\times\Ka_{2n}&=&\Pa_{2n}=P_{\mu,\pi_{2n}}.
\end{eqnarray*}
This ends the proof of  (\ref{bridge-form-2}).\label{ref-kk-pi-proof}

$\bullet$ At every time step $n\geq 0$ by (\ref{log-densities}) we have joint probability measure 
$$
\Pa_{2n}(d(x,y))=\mu(dx)~\Ka_{2n}(x,dy)=e^{-U(x)-W_{2n}(x,y)}~\lambda(dx)\lambda(dy).
$$
The conditional distribution of the first coordinate given the second coincides with the backward transition
$$
\Ka_{2n+1}(y,dx)=\mu(dx)~\frac{d\delta_x\Ka_{2n}}{d\mu \Ka_{2n}}(y)=e^{-W^{\flat}_{2n+1}(y,x)}~\lambda(dx).
$$
For any given index $n\geq 0$ the Markov chain with transitions
$S_{2n+1}:=\Ka_{2n}\Ka_{2n+1}$ alternates forward samples from $\Ka_{2n}$ and backward samples from $\Ka_{2n+1}$ and we have
$$
\mu(dx_1)\Ka_{2n}(x_1,dx_2)\Ka_{2n+1}(x_2,dx_3)=
\mu(dx_3)\Ka_{2n}(x_3,dx_2)~\Ka_{2n+1}(x_2,dx_1).
$$
In summary, this reversible chain 
that coincides the forward/backward Gibbs sampler with target measure $\Pa_{2n}$. 

$\bullet$ At every time step $n\geq 0$ by (\ref{log-densities}) we also have the probability measure 
$$
\Pa_{2n+1}(d(x,y))=\eta(dy)~\Ka_{2n+1}(y,dx)=e^{-V(y)-W_{2n+1}^{\flat}(y,x)}~\lambda(dx)\lambda(dy).
$$
The conditional distribution of the second coordinate given the first coincides with the foward transition
$$
\Ka_{2(n+1)}(x,dy)=\eta(dy)~\frac{d\delta_y\Ka_{2n+1}}{d\eta \Ka_{2n+1}}(x)=e^{-W_{2(n+1)}(x,y)}~\lambda(dy).
$$
For any given index $n\geq 0$ the Markov chain with transitions
$S_{2(n+1)}:=\Ka_{2n+1}\Ka_{2(n+1)}$ alternates backward samples from $\Ka_{2n+1}$  and forward samples from $\Ka_{2(n+1)}$ and we have
$$
\eta(dy_1)\Ka_{2n+1}(y_1,dy_2)\Ka_{2(n+1)}(y_2,dy_3)=
\eta(dy_3)
\Ka_{2n+1}(y_3,dy_2)~\Ka_{2(n+1)}(y_2,dy_1).$$
In summary, this reversible chain 
that coincides the  backward/forward Gibbs sampler with target measure $\Pa_{2n+1}$. 

This yields for any $n\geq 0$   the fixed point equations (\ref{def-pin})

Using (\ref{prop-schp}) and (\ref{log-densities}) we check that
\begin{eqnarray}
\Ka_{2n}(x,dy)&=&e^{-(W_{2n}-W)(x,y)}~K_W(x,dy)\nonumber\\
&=&
e^{-(U_{2n}-U)(x)}~K_W(x,dy)~e^{-V_{2n}(y)}=\frac{K_W(x,dy)~e^{-V_{2n}(y)}}{K_W(e^{-V_{2n}})(x)}.\label{ktos-1}
\end{eqnarray}
In the same vein, using (\ref{prop-schp}) and (\ref{log-densities}) we check that
\begin{eqnarray}
\Ka_{2n+1}(y,dx)&=&e^{-(W^{\flat}_{2n+1}-W^{\flat})(y,x)}~K_{W^{\flat}}(y,dx)\nonumber\\
&=&e^{-(V_{2n+1}-V)(y)}~K_{W^{\flat}}(y,dx)~e^{-U_{2n+1}(x)}=\frac{K_{W^{\flat}}(y,dx)~e^{-U_{2n+1}(x)}}{K_{W^{\flat}}(e^{-U_{2n+1}})(y)}.~~~~\label{ktos-2}
\end{eqnarray}
Using the  Schr\" odinger system (\ref{prop-schp-fp}) and (\ref{log-densities-fp}) we also have
\begin{eqnarray}
K_{W_{\mu,\eta}}(x,dy)&=&e^{-(W_{\mu,\eta}-W)(x,y)}~K_W(x,dy)\nonumber\\
&=&
e^{-(\UU-U)(x)}~K_W(x,dy)~e^{-\VV(y)}=\frac{K_W(x,dy)~e^{-\VV(y)}}{K_W(e^{-\VV})(x)}.\label{ktos-1-b}
\end{eqnarray}
as well as
\begin{eqnarray}
K_{W^{\flat}_{\eta,\mu}}(y,dx)&=&e^{-(W^{\flat}_{\eta,\mu}-W^{\flat})(y,x)}~K_{W^{\flat}}(y,dx)\nonumber\\
&=&e^{-(\VV-V)(y)}~K_{W^{\flat}}(y,dx)~e^{-\UU(x)}=\frac{K_{W^{\flat}}(y,dx)~e^{-\UU(x)}}{K_{W^{\flat}}(e^{-\UU})(y)}.\label{ktos-2-b}
\end{eqnarray}

\subsection*{Sinkhorn potential functions}
Observe that
$$
\begin{array}{l}
\displaystyle\frac{d\Pa_{2n}}{d\Pa_{2n+1}}(x,y)=\frac{d\pi_{2n}}{d\lambda_V}(y)=e^{V_{2n+1}(y)-V_{2n}(y)}\\
\\
\displaystyle\Longrightarrow \pi_{2n}=\lambda_{V^{\pi}_{2n}}\quad \mbox{\rm with}\quad
V^{\pi}_{2n}:=V+(V_{2n}-V_{2n+1}).
\end{array}$$
In the same vein, we have
$$
\begin{array}{l}
\displaystyle\frac{d\Pa_{2n+1}}{d\Pa_{2(n+1)}}(x,y)=\frac{d\pi_{2n+1}}{d\lambda_U}(y)=e^{U_{2(n+1)}(y)-U_{2n+1}(y)}\\
\\
\displaystyle\Longrightarrow \pi_{2n+1}=\lambda_{U^{\pi}_{2n+1}}\quad \mbox{\rm with}\quad
U^{\pi}_{2n+1}:=U+(U_{2n+1}-U_{2(n+1)}).
\end{array}$$
 Using the above formula we readily check (\ref{pi-gibbs}) as well as the entropy equations (\ref{ent-r1}).

Note that
$$
\begin{array}{l}
(\pi_{2n},\eta)= (\mu \Ka_{2n},\pi_{2n-1}\Ka_{2n})
\\
\\
\Longrightarrow
  \Ha(\pi_{2n}~|~\eta)\leq   
    \Ha(\mu~|~\pi_{2n-1})\quad \mbox{\rm and}\quad
    \Ha(\eta~|~\pi_{2n})\leq 
 \Ha(\pi_{2n-1}~|~\mu).
 \end{array}
$$
In the same vein, we have
$$\label{e2}
\begin{array}{l}
(\pi_{2n+1},\mu)= (\eta \Ka_{2n+1},\pi_{2n}\Ka_{2n+1})
\\
\\
\Longrightarrow
  \Ha(\pi_{2n+1}~|~\mu)\leq   
    \Ha(\eta~|~\pi_{2n})\quad \mbox{\rm and}\quad
    \Ha(\mu~|~\pi_{2n+1})\leq 
 \Ha(\pi_{2n}~|~\eta).
 \end{array}
$$
 These inequalities are direct consequence of the well known non expansive properties of the Markov transport 
 map with respect to the relative entropy, see for instance~\cite{dm-03}.

\section{Gradient and Hessian formulae}\label{sec-grad-hess-a}

The gradient and the Hessian of the reference potential function $W$ defined in (\ref{ref-KW}) w.r.t. the first coordinate are denoted by
$$
W_x(y):=\nabla_1W(x,y)=\nabla_2W^{\flat}(y,x)\quad \mbox{\rm and}\quad
W_{x,x}(y):=\nabla_1^2W(x,y)=\nabla_2^2W^{\flat}(y,x).
$$
In this notation, we also have the dual functions
$$
\nabla_2 W(x,y)=\nabla_1W^{\flat}(y,x)=W^{\flat}_y(x)\quad \mbox{\rm and}\quad
\nabla_2^2 W(x,y)=\nabla_1^2W^{\flat}(y,x)=W^{\flat}_{y,y}(x).
$$
\subsection*{Covariance formulae}\label{cov-hess-sec}
Consider the conditional  covariances
functions
\begin{eqnarray*}
 C_{2n}(x)&:=&\frac{1}{2}\int~
  \Ka_{2n}(x,dy_1)~  \Ka_{2n}(x,dy_2)~(W_x(y_1)-W_x(y_2))~(W_x(y_1)-W_x(y_2))^{\prime},
  \\
 C_{2n+1}(y)&:=&\frac{1}{2}\int~
  \Ka_{2n+1}(y,dx_1)~  \Ka_{2n+1}(y,dx_2)~(W^{\flat}_y(x_1)-W^{\flat}_y(x_2))~(W^{\flat}_y(x_1)-W^{\flat}_y(x_2))^{\prime}.
  \end{eqnarray*}
We also denote by $C_{\mu,\eta}$ and $C^{\flat}_{\eta,\mu}$ the conditional covariance functions
functions
\begin{eqnarray*}
C_{\mu,\eta}(x)&:=&\frac{1}{2}\int~
  \Ka_{\WW_{\mu,\eta}}(x,dy_1)~  \Ka_{\WW_{\mu,\eta}}(x,dy_2)~(W_x(y_1)-W_x(y_2))~(W_x(y_1)-W_x(y_2))^{\prime},
  \\
 C^{\flat}_{\eta,\mu}(y)&:=&\frac{1}{2}\int~
  \Ka_{\WW^{\flat}_{\eta,\mu}}(y,dx_1)~  \Ka_{\WW^{\flat}_{\eta,\mu}}(y,dx_2)~(W^{\flat}_y(x_1)-W^{\flat}_y(x_2))~(W^{\flat}_y(x_1)-W^{\flat}_y(x_2))^{\prime}.
  \end{eqnarray*}
  
\begin{examp}\label{ref-hessian-W-f-examp}
Consider the linear-Gaussian  potential $W$ of the form (\ref{def-W}). In this situation we have the gradient and Hessian formulae
\begin{equation}\label{ref-hessian-W-f}
\begin{array}{l}
 W_x(y)=\cchi^{\prime}((\alpha+\beta x)-y)\quad\mbox{and}\quad
W^{\flat}_y(x)=\tau^{-1}(y-(\alpha+\beta x)) \\
\\
\Longrightarrow
W_{x,x}(y)= \cchi^{\prime}\beta
\quad\mbox{and}\quad  W^{\flat}_{y,y}(x)=\tau^{-1}.
\end{array}
\end{equation}
Also note that
$$
W^{\flat}_{y_1}(x)-W^{\flat}_{y_2}(x)=\tau^{-1}(y_1-y_2)
\quad\mbox{and}\quad
 W_{x_1}(y)- W_{x_2}(y)=\cchi^{\prime}\beta (x_1-x_2).
$$
This implies
\begin{eqnarray*}
 W_x(y)-
 \Ka_{2n}\left( W_x\right)(x)&=&-\cchi^{\prime}~(I(y)- \Ka_{2n}(I)(x)),\\
W^{\flat}_y(x)-
 \Ka_{2n+1}( W^{\flat}_y)(y)&=&-\cchi~(I(x)- \Ka_{2n+1}(I)(y))\quad\mbox{\rm with}\quad
I(x):=x.
\end{eqnarray*}
from which we check that
\begin{equation}\label{ref-cov-wg-f}
C_{2n}=\cchi^{\prime}~\Sigma_{2n}~\cchi
\quad\mbox{\rm and}\quad
C_{2n+1}=\cchi~\Sigma_{2n+1}~\cchi^{\prime},
\end{equation}
with the conditional covariances functions $\Sigma_n$ defined in Section~\ref{sinkmaps-sec}.

In the same vein, we check that
\begin{equation}\label{ref-bcov-wg-f}
C_{\mu,\eta}=\cchi^{\prime}~\Sigma_{\mu,\eta}~\cchi
\quad\mbox{\rm and}\quad
C^{\flat}_{\eta,\mu}=\cchi~\Sigma^{\flat}_{\eta,\mu}~\cchi^{\prime}.
\end{equation}
\end{examp}

\subsection*{Sinkhorn bridge potentials}
We set
$$
 \Delta C_{2n+1}:=C_{2n+1}-C_{2n-1}
 \quad \mbox{\rm and}\quad
 \Delta C_{2(n+1)}:=C_{2(n+1)}-C_{2n}.
$$
Note that $\Ka_{0} (W_x)=0$ and $C_{W_{0}}=0$. We also consider the signed operators
 $$
 \Delta \Ka_{2n+1}:=
 \Ka_{2n+1}-\Ka_{2n-1}\quad \mbox{\rm and}\quad
 \Delta \Ka_{2(n+1)}:=\Ka_{2(n+1)}-\Ka_{2n}.
 $$
 Using (\ref{prop-schp}), (\ref{prop-schp-fp}) and Sinkhorn transition formulae (\ref{ktos-1}), (\ref{ktos-2}), (\ref{ktos-1-b}) and (\ref{ktos-2-b})  we check the following lemma.
 \begin{lem}
 For any $n\geq 0$ we have
 \begin{eqnarray}
\nabla U_{2n}(x)&=&\nabla U(x)-\Ka_{2n}( W_x)(x),\nonumber\\
\nabla V_{2n+1}(y)&=&\nabla V(y)-\Ka_{2n+1}( W^{\flat}_y)(y).\label{ref-nablaUVng}
 \end{eqnarray}
 In addition, we have
 \begin{eqnarray}
\nabla \UU(x)&=&\nabla U(x)-K_{W_{\mu,\eta}}( W_x)(x),\nonumber\\
\nabla \VV(y)&=&\nabla V(y)-K_{W^{\flat}_{\eta,\mu}}( W^{\flat}_y)(y).\label{ref-nablaUVng-b}
 \end{eqnarray} 
 \end{lem}
Combining (\ref{log-densities}) with
 (\ref{ref-nablaUVng}) we check that
  \begin{eqnarray}
   \nabla_1 W_{2n}(x,y)&=&W_x(y)-\Ka_{2n}( W_x)(x),\nonumber\\
   \nabla_1 W_{2n+1}^{\flat}(y,x)&=&W^{\flat}_y(x) -\Ka_{2n+1}( W^{\flat}_y)(y). \label{ref-nabla1W}
 \end{eqnarray} 
Taking another differential we obtain the following lemma.
\begin{lem}
For any $n\geq 0$ we have
 \begin{eqnarray}
\nabla^2 U_{2n}(x)&=&\nabla^2 U(x)-\Ka_{2n}( W_{x,x})(x)+C_{2n}(x),\nonumber\\
\nabla^2 V_{2n+1}(y)&=&\nabla^2 V(y)-\Ka_{2n+1}( W^{\flat}_{y,y})(y)+C_{2n+1}(y),\label{ref-nablaUVn1g}
 \end{eqnarray}
 as well as the entropic potential Hessian formulae
 \begin{eqnarray}
\nabla^2 \UU(x)&=&\nabla^2 U(x)-K_{W_{\mu,\eta}}( W_{x,x})(x)+C_{\mu,\eta}(x),\nonumber\\
\nabla^2 \VV(y)&=&\nabla^2 V(y)-K_{W^{\flat}_{\eta,\mu}}( W^{\flat}_{y,y})(y)+C^{\flat}_{\eta,\mu}(y).\label{ref-nablaUVn1g-b}
 \end{eqnarray}
 \end{lem} 
 Next lemma provides some Hessian formulae for the  Boltzmann-Gibbs potential functions (\ref{pi-gibbs}).
 \begin{lem}
For any $n\geq 1$ we have
 \begin{eqnarray}
\nabla^2U^{\pi}_{2n-1}(x)&=&\nabla^2U(x)+ \Delta \Ka_{2n}( W_{x,x})(x)- \Delta C_{2n}(x),\nonumber\\
\nabla^2V^{\pi}_{2n}(y)&=&\nabla^2V(y)+ \Delta \Ka_{2n+1}( W^{\flat}_{y,y})(y)- \Delta C_{2n+1}(y).\label{HessUpi}
 \end{eqnarray}
 For $n=0$ we have
 \begin{equation}\label{Hess-Vpi-0}
 \nabla^2V^{\pi}_{0}(y)=
\Ka_{1}( W^{\flat}_{y,y})(y)-C_{1}(y).
\end{equation}
\end{lem}
\proof
 Recalling that $V_{2n}=V_{2n-1}$ and using (\ref{pi-gibbs-p}) for any $n\geq 1$ we have
 \begin{eqnarray*}
\nabla^2V^{\pi}_{2n}(y)&=&\nabla^2V(y)+(\nabla^2V_{2n-1}-\nabla^2V_{2n+1})(y).
 \end{eqnarray*}
 Similarly, recalling that $U_{2n}=U_{2n+1}$, for any $n\geq 0$ we check that
  \begin{eqnarray*}
\nabla^2U^{\pi}_{2n+1}(x)&=&\nabla^2U(x)+(\nabla^2U_{2n+1}-\nabla^2U_{2(n+1)})(x).
 \end{eqnarray*}
 We end the proof of (\ref{HessUpi}) using (\ref{ref-nablaUVn1g}).
 Finally note that
 $$
\begin{array}{l}
V^{\pi}_{0}:=V+(V_{0}-V_{1})=V-V_1\\
\\
\Longrightarrow
\nabla^2V^{\pi}_{0}=\nabla^2V-\nabla^2V_1\Longrightarrow (\ref{Hess-Vpi-0}).
\end{array}$$
 This ends the proof of the lemma.
 \cqfd
 \begin{examp}\label{ex-gauss-UVn}
Consider a linear-Gaussian  potential $W$ of the form (\ref{def-W}). In this situation using (\ref{ref-hessian-W-f}) and (\ref{ref-cov-wg-f}) for any $n\geq 0$ the Hessian formulae (\ref{ref-nablaUVn1g}) simplify to
 \begin{eqnarray}
\nabla^2 U_{2n}&=&\nabla^2 U- \cchi^{\prime}\beta+\cchi^{\prime}~\Sigma_{2n}~\cchi,\nonumber\\
\nabla^2 V_{2n+1}&=&\nabla^2 V-\tau^{-1}+\cchi~\Sigma_{2n+1}~\cchi^{\prime}.\label{ref-nablaUVn1gex}
 \end{eqnarray}
Using  (\ref{ref-hessian-W-f}) (\ref{ref-nablaUVng}) and (\ref{ref-nablaUVng-b}) we also have
  \begin{eqnarray}
\nabla U_{2n}(x)&=&\nabla U(x)-\cchi^{\prime}(\alpha+\beta x)+\cchi^{\prime}~\Ka_{2n}(I)(x),\nonumber\\
\nabla V_{2n+1}(y)&=&\nabla V(y)-\tau^{-1}~y+\tau^{-1}\left(\alpha+\beta~\Ka_{2n+1}(I)(y)\right),\label{ref-nablaUVng-bg-n}
 \end{eqnarray} 
 as well as
  \begin{eqnarray}
\nabla \UU(x)&=&\nabla U(x)-\cchi^{\prime}(\alpha+\beta x)+\cchi^{\prime}~K_{\WW_{\mu,\eta}}(I)(x),\nonumber\\
\nabla \VV(y)&=&\nabla V(y)-\tau^{-1}~y+\tau^{-1}\left(\alpha+\beta~K_{\WW^{\flat}_{\eta,\mu}}(I)(y)\right).\label{ref-nablaUVng-bg}
 \end{eqnarray} 
By (\ref{ref-hessian-W-f}) also note that (\ref{ref-nabla1W}) resumes to
   \begin{eqnarray}
   \nabla_1 W_{2n}(x,y)&=&-\cchi^{\prime}~(y-\Ka_{2n}(I)(x)),\nonumber\\
   \nabla_1 W_{2n+1}^{\flat}(y,x)&=&-\cchi~(x- \Ka_{2n+1}(I)(y)). \label{ref-nabla1W-bis}
 \end{eqnarray} 
 Taking the differential in (\ref{ref-nablaUVng-bg-n}) (\ref{ref-nabla1W-bis}) and using (\ref{rule-ah}) and (\ref{ref-nablaUVn1gex}) we check that
  \begin{eqnarray}
\nabla^2 (U_{2n}- U)(x)+ \cchi^{\prime}\beta=\cchi^{\prime}~\Sigma_{2n}(x)~\cchi&=&\nabla \Ka_{2n}(I)(x)~\cchi=   \nabla_1^2 W_{2n}(x,y),\nonumber\\
\nabla^2 (V_{2n+1}- V)(y)+\tau^{-1}=\cchi~\Sigma_{2n+1}(y)~\cchi^{\prime}&=&\nabla K_{2n+1}(I)(y)~\cchi^{\prime}= \nabla_1^2 W_{2n+1}^{\flat}(y,x).~~~~~~~~~\label{finex-n}\end{eqnarray}
In  the same vein, taking the differential in (\ref{ref-nablaUVng-bg}) and using (\ref{ref-bcov-wg-f}) we have
  \begin{eqnarray}
\nabla^2 (\UU- U)+ \cchi^{\prime}\beta&=&\cchi^{\prime}~\Sigma_{\mu,\eta}~\cchi=\nabla K_{\WW_{\mu,\eta}}(I)~\cchi,\nonumber\\
\nabla^2 (\VV- V)+\tau^{-1}&=&\cchi~\Sigma^{\flat}_{\eta,\mu}~\cchi^{\prime}=\nabla K_{\WW^{\flat}_{\eta,\mu}}(I)~\cchi^{\prime}.\label{finex}\end{eqnarray}
 In the same vein, formulae (\ref{HessUpi}) for $n\geq 1$ simplify to
 \begin{eqnarray}
\nabla^2U^{\pi}_{2n-1}&=&\nabla^2U+ \cchi^{\prime}(\Sigma_{2(n-1)}-\Sigma_{2n})
\cchi,\nonumber\\
\nabla^2V^{\pi}_{2n}&=&\nabla^2V+\cchi~(\Sigma_{2n-1}-\Sigma_{2n+1})
 \cchi^{\prime}\quad\mbox{and}\quad
  \nabla^2V^{\pi}_{0}=
\tau^{-1}-\cchi~\Sigma_{1}~\cchi^{\prime}.
 \label{HessUpiex}
 \end{eqnarray}
\end{examp}
\subsection*{Transition potentials}

\begin{lem}\label{lem-appendix-Hess}
For any $n\geq 0$ we have  (\ref{nabla1-WN}) as well as the Hessian formulae
\begin{eqnarray}
  \nabla_2^2W_{2(n+1)}(x,y)
  &=&\nabla^2 V(y)+\left[W^{\flat}_{y,y}(x)-\Ka_{2n+1}( W^{\flat}_{y,y})(y)\right]+C_{2n+1}(y)\nonumber\\
   \nabla_2^2W^{\flat}_{2n+1}(y,x)
    &=&\nabla^2 U(x)+\left[W_{x,x}(y)-\Ka_{2n}( W_{x,x})(x)\right]+C_{2n}(x)\label{nabla2-WN}
\end{eqnarray}
In addition, for any $n\geq 0$ we have
\begin{eqnarray}
  \nabla_2W_{2n}(x_1,y)-  \nabla_2W_{2n}(x_2,y)&=&
   W^{\flat}_y(x_1)- W^{\flat}_y(x_2)\nonumber\\
     \nabla_2W^{\flat}_{2n+1}(y_1,x)-  \nabla_2W^{\flat}_{2n+1}(y_2,x)&=& W_x(y_1)-W_x(y_2)\label{nabla2-WN-rr}
\end{eqnarray}
as well as
\begin{eqnarray}
  \nabla_2W_{2n}(x,y_1)-  \nabla_2W_{2n}(x,y_2)&=&  \nabla V_{2n}(y_1)-\nabla V_{2n}(y_2)+W^{\flat}_{y_1}(x)-W^{\flat}_{y_2}(x)\nonumber\\
    \nabla_2W^{\flat}_{2n+1}(y,x_1)-    \nabla_2W^{\flat}_{2n+1}(y,x_2)&=&  \nabla U_{2n}(x_1)- \nabla U_{2n}(x_2)+  W_{x_1}(y)-W_{x_2}(y)~~~~~~\label{nabla2-WN-rr-2}
\end{eqnarray}
\end{lem}

\proof
 Recalling that $V_{2n+1}=V_{2(n+1)}$ and using (\ref{log-densities}) we check that
  \begin{eqnarray*}
  \nabla_2W_{2(n+1)}(x,y)&=&  W^{\flat}_y(x)+  \nabla V_{2n+1}(y)=\nabla V(y)+ W^{\flat}_y(x)-\Ka_{2n+1}( W^{\flat}_y)(y), \\
  \nabla_2^2W_{2(n+1)}(x,y)&=&  W^{\flat}_{y,y}(x)+  \nabla^2V_{2n+1}(y).
  \end{eqnarray*}
In the same vein, recalling that $U_{2n}=U_{2n+1}$ taking the differentials of the log-densities formulae (\ref{log-densities})  we check that
 \begin{eqnarray*}
  \nabla_2W^{\flat}_{2n+1}(y,x)&=&  \nabla U_{2n}(x)+  W_x(y)=\nabla U(x)+W_x(y)-\Ka_{2n}( W_x)(x),\\
    \nabla_2^2W^{\flat}_{2n+1}(y,x)&=&  \nabla^2 U_{2n}(x)+  W_{x,x}(y).
\end{eqnarray*}
This implies that
\begin{eqnarray*}
  \nabla_2W_{2n}(x_1,y)-  \nabla_2W_{2n}(x_2,y)&=&\nabla_2 W(x_1,y)-\nabla_2 W(x_2,y)=
   W^{\flat}_y(x_1)- W^{\flat}_y(x_2),\nonumber\\
     \nabla_2W^{\flat}_{2n+1}(y_1,x)-  \nabla_2W^{\flat}_{2n+1}(y_2,x)
     &=&\nabla_2W^{\flat}(y_1,x)-\nabla_2W^{\flat}(y_2,x)\\
    & =&\nabla_1W(x,y_1)-\nabla_1W(x,y_2)= W_x(y_1)-W_x(y_2).
\end{eqnarray*}
Formulae (\ref{nabla1-WN}) are checked using the gradient formulae stated above.
Formulae (\ref{nabla2-WN}) are now consequences of (\ref{ref-nablaUVn1g}).
This ends the proof of the lemma.
\cqfd

Next lemma is a consequence of the fixed point equations (\ref{log-densities-fp}).
\begin{lem}
We have the gradient formulae
\begin{eqnarray}
  \nabla_2\WW_{\mu,\eta}(x,y_1)-  \nabla_2\WW_{\mu,\eta}(x,y_2)&=&    \nabla \VV(y_1)- \nabla \VV(y_2)
  +W^{\flat}_{y_1}(x)-W^{\flat}_{y_2}(x),\nonumber\\
      \nabla_2W^{\flat}_{\eta,\mu}(y,x_1)-  \nabla_2W^{\flat}_{\eta,\mu}(y,x_2)&=&  \nabla \UU(x_1)- \nabla \UU(x_2)+  W_{x_1}(y)-
      W_{x_2}(y).\label{nabla2-WWN-rr-2}
\end{eqnarray}
\end{lem}
\proof
Using (\ref{log-densities-fp}) we check that
$$
  \nabla_2\WW_{\mu,\eta}(x,y)=    \nabla \VV(y)+W^{\flat}_y(x)\quad\mbox{and}\quad
    \nabla_2W^{\flat}_{\eta,\mu}(y,x)=  \nabla \UU(x)+  W_x(y).
  $$
  This ends the proof of the lemma.
\cqfd

 \begin{examp}\label{l-exam-DeltaUV}
Consider the linear-Gaussian  potential $W$ of the form (\ref{def-W}). In this situation using 
(\ref{ref-nablaUVng-bg-n}) and (\ref{nabla2-WN-rr-2}) for any $n\geq 0$ we check that
$$
  \begin{array}{l}
   \nabla_2W^{\flat}_{2n+1}(y,x_1)-    \nabla_2W^{\flat}_{2n+1}(y,x_2)\\
   \\
=\nabla U_{2n}(x_1)-\nabla U_{2n}(x_2)+\cchi^{\prime}\beta( x_1-x_2)\\
\\
=\nabla U(x_1)-\nabla U(x_2)+\cchi^{\prime}~(\Ka_{2n}(I)(x_1)-\Ka_{2n}(I)(x_2)).
 \end{array} 
 $$
For a column vector
 function $h=(h^i)_{1\leq i\leq d}$ we have
$$
h(x_1)-h(x_2)=\int_0^1~(\nabla h(x_2+s(x_1-x_2)))^{\prime}ds~(x_1-x_2).
$$
Using (\ref{rule-ah}) and (\ref{finex-n}) we check that
 $$
 \cchi^{\prime}~(\Ka_{2n}(I)(x_1)-\Ka_{2n}(I)(x_2))=\int_0^1\cchi^{\prime}~\Sigma_{2n}(x_2+s(x_1-x_2))~\cchi~ds~(x_1-x_2).
 $$
Therefore for any $n\geq 0$ we have
\begin{equation}\label{nbU}
 \begin{array}{l}
 (x_1-x_2)^{\prime}  \left(\nabla_2W^{\flat}_{2n+1}(y,x_1)-    \nabla_2W^{\flat}_{2n+1}(y,x_2)\right)\\
 \\
\displaystyle\geq 
 (x_1-x_2)^{\prime} \left(\nabla U(x_1)-\nabla U(x_2)\right) + (x_1-x_2)^{\prime}~\cchi^{\prime}~\Sigma_{2n}(x_1,x_2)~\cchi~(x_1-x_2),
\end{array} 
\end{equation}
with
$$
\Sigma_{2n}(x_1,x_2):=
 \int_0^1~\Sigma_{2n}(x_2+s(x_1-x_2))~ds.
$$
In the same vein, for any $n\geq 1$ we have
$$
  \begin{array}{l}
   \nabla_2W_{2n}(x,y_1)-  \nabla_2W_{2n}(x,y_2)\\
   \\
=\nabla V_{2n}(y_1)-\nabla V_{2n}(y_2)+\tau^{-1}~(y_1-y_2)\\
\\
=\nabla V(y_1)-\nabla V(y_2)+\cchi~(\Ka_{2n-1}(I)(y_1)-\Ka_{2n-1}(I)(y_2)).
 \end{array} 
 $$
 Therefore for any $n\geq 1$ we have
\begin{equation}\label{nbV}
   \begin{array}{l}
 (y_1-y_2)^{\prime}\left(   \nabla_2W_{2n}(x,y_1)-  \nabla_2W_{2n}(x,y_2)\right)\\
 \\
\displaystyle \geq 
 (y_1-y_2)^{\prime}\left(\nabla V(y_1)-\nabla V(y_2) \right)+
 (y_1-y_2)^{\prime}~\cchi~\Sigma_{2n-1}(y_1,y_2)~\cchi^{\prime}~ (y_1-y_2),
  \end{array} 
  \end{equation}
with
$$
\Sigma_{2n-1}(y_1,y_2):=
\int_0^1\int_0^1\Sigma_{2n-1}(y_2+s(y_1-y_2))~ds.
$$

\end{examp}

Using (\ref{log-densities-fp}) we also check the following lemma.
\begin{lem}\label{lem-appendix-Hess-b}
We have the gradient formulae
\begin{eqnarray}
  \nabla_2\WW_{\mu,\eta}(x_1,y)-  \nabla_2\WW_{\mu,\eta}(x_2,y)&=&
   \nabla_2W(x_1,y)- \nabla_2W(x_2,y),\nonumber
\\
  \nabla_2\WW^{\flat}_{\eta,\mu}(y_1,x)-  \nabla_2\WW^{\flat}_{\eta,\mu}(y_2,x)&=&
   \nabla_2W^{\flat}(y_1,x)- \nabla_2W^{\flat}(y_2,x),
\label{nabla1-WN-b}
\end{eqnarray}
as well as the Hessian formulae
\begin{eqnarray}
  \nabla_2^2\WW_{\mu,\eta}(x,y)
  &=&\nabla^2 V(y)+\left[W^{\flat}_{y,y}(x)-K_{\WW^{\flat}_{\eta,\mu}}( W^{\flat}_{y,y})(y)\right]+C^{\flat}_{\eta,\mu}(y),\nonumber\\
   \nabla_2^2\WW^{\flat}_{\eta,\mu}(y,x)
    &=&\nabla^2 U(x)+\left[W_{x,x}(y)-K_{\WW_{\mu,\eta}}( W_{x,x})(x)\right]+C_{\mu,\eta}(x).\label{nabla2-WN-b}
\end{eqnarray}
\end{lem}

 \begin{examp}\label{l-exam-Vpi}
Consider the linear-Gaussian  potential $W$ of the form (\ref{def-W}). In this situation using (\ref{ref-hessian-W-f}) and (\ref{ref-cov-wg-f}) for any $n\geq 0$ the Hessian formulae (\ref{nabla2-WN})  simplify to
\begin{eqnarray}
  \nabla_2^2W_{2(n+1)}(x,y)
  &=&\nabla^2 V(y)+\cchi~\Sigma_{2n+1}(y)~\cchi^{\prime},\nonumber\\
   \nabla_2^2W^{\flat}_{2n+1}(y,x)
    &=&\nabla^2 U(x)+\cchi^{\prime}~\Sigma_{2n}(x)~\cchi.\label{nabla2-WN-ex}
\end{eqnarray}
In the same same vein, (\ref{nabla2-WN-b}) reduces to
\begin{eqnarray}
  \nabla_2^2\WW_{\mu,\eta}(x,y)
  &=&\nabla^2 V(y)+\cchi~\Sigma^{\flat}_{\eta,\mu}(y)~\cchi^{\prime},\nonumber\\
   \nabla_2^2\WW^{\flat}_{\eta,\mu}(y,x)
    &=&\nabla^2 U(x)+\cchi^{\prime}~\Sigma_{\mu,\eta}(x)~\cchi.\label{nabla2-WN-ex-b}
\end{eqnarray}
In addition, using the gradient formulae (\ref{ref-hessian-W-f}) and (\ref{nabla2-WN-rr}) 
we also check that
\begin{eqnarray*}
  \nabla_2W_{2n}(x_1,y)-  \nabla_2W_{2n}(x_2,y)&=&
W_y^{\flat}(x_1)- W_y^{\flat}(x_2)=\cchi ~(x_2-x_1),\nonumber
\\
  \nabla_2W^{\flat}_{2n+1}(y_1,x)-  \nabla_2W^{\flat}_{2n+1}(y_2,x)&=&
  W_x(y_1)-    W_x(y_2)=\cchi^{\prime}~(y_2-y_1).
\end{eqnarray*}
Note that
\begin{equation}\label{ref-flat-1}
  \nabla_2^2W^{\flat}_{1}(y,x)
    =\nabla^2 U(x)+\cchi^{\prime}~\tau~\cchi\geq  u_+^{-1}+\cchi^{\prime}~\tau~\cchi=
    u_+^{-1}+\beta^{\prime}\tau^{-1}\beta.
\end{equation}
 Applying Brascamp-Lieb and Cram\'er-Rao inequalities (\ref{bl-cr}) we check that
 \begin{eqnarray*}
\cchi~ \Sigma_{1}~\cchi^{\prime}
&\leq&  \cchi~(u_+^{-1}+\beta^{\prime}\tau^{-1}\beta)^{-1}~ \cchi^{\prime}\\
&=&
  \tau^{-1}\beta~(u_+^{-1}+\beta^{\prime}\tau^{-1}\beta)^{-1}\beta^{\prime}  \tau^{-1}
  = \tau^{-1}-\left( \tau+\beta~ u_+~\beta^{\prime}\right)^{-1}.
\end{eqnarray*}
The last assertion comes from the matrix inversion lemma. Using (\ref{HessUpiex}) we conclude that
$$
\nabla^2V^{\pi}_{0}=
\tau^{-1}-\cchi~\Sigma_{1}~\cchi^{\prime}\geq \left( \tau+\beta~ u_+~\beta^{\prime}\right)^{-1}.
$$
\end{examp}

\section{Strongly convex models}\label{sec-a-stm}

Consider the
linear-Gaussian reference transitions $K_W$ associated with some potential $W$ of the form (\ref{def-W}) for some parameters $(\alpha,\beta,\tau)$.
Assume that the potential functions $(U,V)$
satisfy condition (\ref{UV-d2}).

\subsection*{Covariance estimates}
We further assume condition (\ref{UV-d2}) is satisfied.
Consider the positive matrices $\varsigma_n,\tau_n$ defined 
by $ \varsigma_{0}=\tau=\tau_0$ and the recursive  formulae
\begin{equation}\label{cov-ess}
 \begin{array}{rcl}
\varsigma_{2n+1}^{-1}&=&u_-^{-1}+\cchi^{\prime}~\tau_{2n}~\cchi\\
\tau_{2n+1}^{-1} &=&u_+^{-1}+\cchi^{\prime}~\varsigma_{2n}~\cchi
 \end{array} \quad \&\quad \begin{array}{rcl}
\varsigma_{2(n+1)}^{-1}&=&v^{-1}_-+\cchi~\tau_{2n+1}~\cchi^{\prime},\\ 
\tau_{2(n+1)}^{-1}&=&v^{-1}_++\cchi~\varsigma_{2n+1}~\cchi^{\prime}.
 \end{array}
\end{equation}
In this context,
the Brascamp-Lieb and Cram\'er-Rao inequalities presented in Section~\ref{BL-CR-ineq-sec} yield the following conditional covariance estimates.
\begin{prop}\label{prop-BL-CR}
For any $n\geq 0$ we have the covariance estimates
$$
\varsigma_{n}\leq \Sigma_{n}\leq \tau_{n}, 
$$
as well as the Hessian estimates
\begin{equation}
\tau_{2n}^{-1}
\leq  \nabla_2^2W_{2n}\leq
\varsigma_{2n}^{-1}\quad \mbox{and}\quad
\tau_{2n+1}^{-1} 
\leq 
  \nabla_2^2W^{\flat}_{2n+1}
    \leq \varsigma_{2n+1}^{-1}.
\end{equation}
\end{prop}

\proof
Note that
$$
 W_{0}(x,y)=W(x,y)\Longrightarrow
 \nabla^2_yW_0(x,y)=\nabla^2_yW(x,y)=\tau^{-1}
 \Longrightarrow \Sigma_{0}(x)=\tau.
 $$
By (\ref{nabla2-WN-ex}), this implies that
$$
 u_+^{-1}+\cchi^{\prime}\tau\cchi\leq
  \nabla_2^2W^{\flat}_{1}
    \leq u_-^{-1}+\cchi^{\prime}\tau\cchi.
$$
Applying Brascamp-Lieb and Cram\'er-Rao inequalities (\ref{bl-cr}) we check that
$$
\varsigma_{1}:=(u_-^{-1}+\cchi^{\prime}\tau\cchi)^{-1}
\leq \Sigma_{1}\leq \tau_{1}:= ( u_+^{-1}+\cchi^{\prime}\tau\cchi)^{-1}.
$$
Assume the inequalities are satisfied at rank $(2n-1)$. In this case, we have
$$
v^{-1}_++\cchi~\varsigma_{2n-1}~\cchi^{\prime}
\leq  \nabla_2^2W_{2n}\leq
v^{-1}_-+\cchi~\tau_{2n-1}~\cchi^{\prime}.
$$
Applying Brascamp-Lieb and Cram\'er-Rao inequalities (\ref{bl-cr})  we check that
$$
\varsigma_{2n}=(v^{-1}_-+\cchi~\tau_{2n-1}~\cchi^{\prime})^{-1}
\leq \Sigma_{2n}\leq \tau_{2n}=(v^{-1}_++\cchi~\varsigma_{2n-1}~\cchi^{\prime})^{-1}.
$$
By (\ref{nabla2-WN-ex}) this implies that
$$
u_+^{-1}+\cchi^{\prime}~\varsigma_{2n}~\cchi
\leq 
  \nabla_2^2W^{\flat}_{2n+1}
    \leq u_-^{-1}+\cchi^{\prime}~\tau_{2n}~\cchi.
$$
Applying once more Brascamp-Lieb and Cram\'er-Rao inequalities (\ref{bl-cr})  we check that
$$
\varsigma_{2n+1}=(u_-^{-1}+\cchi^{\prime}~\tau_{2n}~\cchi)^{-1}\leq\Sigma_{2n+1}\leq\tau_{2n+1} =(u_+^{-1}+\cchi^{\prime}~\varsigma_{2n}~\cchi)^{-1}.
$$
This ends the proof of the proposition.
\cqfd

We set
$$
\Phi_{v}(s):=\left(v^{-1}+\cchi~s~ \cchi^{\prime}
\right)^{-1}\quad \mbox{\rm and}\quad
\Phi^{\prime}_{v}(s):=
\left(v^{-1}+\cchi^{\prime}~s~\cchi
\right)^{-1}.
$$
Consider the positive definite matrices defined by

\begin{equation}\label{bbPhi}
 \begin{array}{rcl}
v_{n+1}&=&
\Phi_{v_-}(\omega_n)\\
 \sigma_{n+1}&=&\Phi_{v_+}( s_{n})
 \end{array} \quad \&\quad \begin{array}{rcl}
s_{n+1}&=&\Phi^{\prime}_{u_-}(\sigma_n),\\ 
{\omega}_{n+1}&=&\Phi^{\prime}_{u_+}(v_{n}).
 \end{array}
\end{equation}
with the initial conditions
$$
( v_{0},\sigma_0)=(0,v_+)\quad \mbox{\rm and}\quad
( s_{0},\omega_0)=(0,u_+).
$$
\begin{lem}
For any $n\geq 0$ we have
\begin{equation}\label{bcov}
 v_{n}\leq \Sigma_{\mu,\eta}\leq \sigma_n\quad \mbox{and}\quad
 s_{n}\leq  \Sigma^{\flat}_{\eta,\mu}\leq \omega_n.
\end{equation}
In addition, for any $n\geq 1$ we have
\begin{equation}\label{bcov-hess}
 \sigma_{n}^{-1}\leq \nabla_2^2\WW_{\mu,\eta}
  \leq  v_{n}^{-1}
  \quad \mbox{and}\quad
\omega_n^{-1}\leq \nabla_2^2\WW^{\flat}_{\eta,\mu}
  \leq s_n^{-1}.
\end{equation}
\end{lem}
\proof
By (\ref{nabla2-WN-ex-b}) we have
$$
v_+^{-1}\leq 
  \nabla_2^2\WW_{\mu,\eta}\quad \mbox{\rm and}\quad
u_+^{-1}\leq  \nabla_2^2\WW^{\flat}_{\eta,\mu}.
$$
Applying Brascamp-Lieb and Cram\'er-Rao inequalities (\ref{bl-cr}) we check that
$$
 v_{0}:=0\leq \Sigma_{\mu,\eta}\leq \sigma_0:=v_+\quad \mbox{\rm and}\quad
 s_{0}:=0\leq  \Sigma^{\flat}_{\eta,\mu}\leq \omega_0:=u_+.
$$
Assume that (\ref{bcov}) is met at rank $n$.
In this situation, we have
$$
\Phi_{v_+}( s_{n})^{-1}\leq \nabla_2^2\WW_{\mu,\eta}(x,y)
  =\nabla^2 V(y)+\cchi~\Sigma^{\flat}_{\eta,\mu}(y)~\cchi^{\prime}\leq (\Phi_{v_-}(\omega_n))^{-1},
$$
as well as
$$
\Phi^{\prime}_{u_+}( v_{n})^{-1}\leq \nabla_2^2\WW^{\flat}_{\eta,\mu}(y,x)
    =\nabla^2 U(x)+\cchi^{\prime}~\Sigma_{\mu,\eta}(x)~\cchi\leq \Phi^{\prime}_{u_-}(\sigma_n)^{-1}.
$$
Applying Brascamp-Lieb and Cram\'er-Rao inequalities (\ref{bl-cr}) we check that
$$
v_{n+1}=
\Phi_{v_-}(\omega_n)
\leq \Sigma_{\mu,\eta}\leq \sigma_{n+1}=\Phi_{v_+}( s_{n}),
$$
and
$$
s_{n+1}=\Phi^{\prime}_{u_-}(\sigma_n)\leq \Sigma^{\flat}_{\eta,\mu}\leq 
{\omega}_{n+1}=\Phi^{\prime}_{u_+}(v_{n}).
$$
This ends the proof of the lemma.
\cqfd

\subsection*{Riccati matrix difference equations}

For a given invertible matrix $\gamma\in \Ga l_d$ we set
 $$
 \begin{array}{lccl}
\Psi_{\gamma}: &\Sa^0_d &\mapsto &\Sa^+_d\\
&v &\leadsto &\Psi_{\gamma}(v):=(I+\gamma v \gamma^{\prime})^{-1}\leq \Psi_{\gamma}(0)=I. 
\end{array}
$$
Note the decreasing property
\begin{equation}\label{ricc-maps-decincr}
v_1\leq v_2\Longrightarrow \Psi_{\gamma}(v_1)\geq \Psi_{\gamma}(v_2).
\end{equation}
\begin{lem}[\cite{adm-24}]\label{varphitoRicc}
We have the factorization formula
\begin{equation}\label{ff-1}
\Psi_{\gamma}\circ\Psi_{\gamma^{\prime}}=
\mbox{\rm Ricc}_{\psi(\gamma)}\quad\mbox{with}\quad
\psi(\gamma):=(\gamma\gamma^{\prime})^{-1}.
\end{equation}
In addition, for any $n\geq 0$ we have
\begin{equation}\label{ff-2}
\Psi_{\gamma^{\prime}}\circ 
\mbox{\rm Ricc}^n_{\psi(\gamma)}=\mbox{\rm Ricc}_{\psi(\gamma^{\prime})}^n
\circ \Psi_{\gamma^{\prime}}\quad\mbox{and}\quad
\Psi_{\gamma^{\prime}}(r_{\psi(\gamma)})=r_{\psi(\gamma^{\prime})}.
\end{equation}
 \end{lem}
The proof of the factorization formula is used in the proof of  Theorem 4.3 in~\cite{adm-24}. For completeness a detailed proof is provided in the appendix on page~\pageref{varphitoRicc-proof}.

 \begin{defi}
Consider the rescaled matrices
$$
 \begin{array}{rcl}
\overline{\varsigma}_{2n+1}&:=&u_-^{-1/2}~
\varsigma_{2n+1}~u_-^{-1/2}\\
\overline{\tau}_{2n+1}&=&u_+^{-1/2}~
\tau_{2n+1}~u_+^{-1/2}
 \end{array} \qquad \&\qquad \begin{array}{rcl}
\overline{\varsigma}_{2n}&:=&v^{-1/2}_-~\varsigma_{2n}~
v^{-1/2}_-,\\
\overline{\tau}_{2n}&=&v^{-1/2}_+~{\tau}_{2n}~v^{-1/2}_+.
  \end{array}
  $$
 \end{defi}

Consider the matrices
 $$
 \gamma_0:=v^{1/2}_-~\cchi~u_+^{1/2}
\qquad \gamma_1 :=u_-^{1/2}~\cchi^{\prime}~v^{1/2}_+\quad\mbox{\rm and}\quad
( \overline{\gamma}_0, \overline{\gamma}_1):=(\gamma_1^{\prime},\gamma_0^{\prime}).
 $$
In this notation (\ref{cov-ess}) takes the following form
\begin{equation}\label{ovvt}
 \begin{array}{rcl}
\overline{\varsigma}_{2n+1}&=&\Psi_{ \gamma_1}(\overline{\tau}_{2n})\\
\overline{\tau}_{2n+1} &=&\Psi_{\overline{\gamma}_1}(\overline{\varsigma}_{2n})
 \end{array} \qquad \&\qquad \begin{array}{rcl}
\overline{\varsigma}_{2(n+1)}&=&\Psi_{ \gamma_0}(\overline{\tau}_{2n+1}),\\ 
\overline{\tau}_{2(n+1)}&=&\Psi_{  \overline{\gamma}_0}(\overline{\varsigma}_{2n+1}).\end{array}
\end{equation}

Consider the positive definite matrices
  $$
\left( \varpi_0,\overline{\varpi}_0\right):=\left(\psi(\gamma_0), \psi(\overline{\gamma}_0)\right)
 \quad \mbox{and}\quad
\left(  \varpi_1, \overline{\varpi}_1\right):=\left( \psi(\gamma_1),\psi(\overline{\gamma}_1)\right).
$$

Using Lemma~\ref{varphitoRicc} we check
\begin{eqnarray*}
\Psi_{ \gamma_1}\circ\Psi_{  \overline{\gamma}_0}&=&
\Psi_{ \gamma_1}\circ\Psi_{ \gamma_1^{\prime}}=\mbox{\rm Ricc}_{\psi(\gamma_1)},\\
\Psi_{\overline{\gamma}_1}\circ \Psi_{ \gamma_0}&=&
\Psi_{\overline{\gamma}_1}\circ \Psi_{\overline{\gamma}_1^{\prime}}=
\mbox{\rm Ricc}_{\psi(\overline{\gamma}_1)},\\
\Psi_{ \gamma_0}\circ \Psi_{\overline{\gamma}_1}&=&
\Psi_{ \gamma_0}\circ \Psi_{\gamma^{\prime}_0}=\mbox{\rm Ricc}_{\psi(\gamma_0)},\\
\Psi_{  \overline{\gamma}_0}\circ \Psi_{ \gamma_1}&=&
\Psi_{  \overline{\gamma}_0}\circ \Psi_{ \overline{\gamma}_0^{\prime}}=
\mbox{\rm Ricc}_{\psi(\overline{\gamma}_0)}.
 \end{eqnarray*}
This  yields for any $n\geq 1$  the Riccati matrix difference equations (\ref{ricc-interm}) as well as 
\begin{equation}\label{tmde}
 \begin{array}{rclcccrcl}
\overline{\varsigma}_{2n}&=&\mbox{\rm Ricc}_{\varpi_0}\left(\overline{\varsigma}_{2(n-1)}\right) &&\quad&\overline{\varsigma}_{2n+1}&=&\mbox{\rm Ricc}_{\varpi_1}\left(\overline{\varsigma}_{2n-1}\right).
 \end{array}
 \end{equation}
 with the initial conditions
$$
\begin{array}{rcl}
\overline{\varsigma}_{0}&:=&v^{-1/2}_-~\tau~
v^{-1/2}_-\\
\overline{\tau}_{0}&=&v^{-1/2}_+~\tau~v^{-1/2}_+
  \end{array}
 \qquad \&\qquad 
  \begin{array}{rcl}
\overline{\varsigma}_{1}&=&\Psi_{ \gamma_1}(\overline{\tau}_{0}),\\
\overline{\tau}_{1} &=&\Psi_{\overline{\gamma}_1}(\overline{\varsigma}_{0}).
 \end{array}
  $$
\begin{rmk}\label{rmk-0v}
\begin{itemize}
\item Assume (\ref{bary-ref-vu}) holds for some 
$u_->0$ and $v_+>0$. 
Following word-for-word the proof of Proposition~\ref{prop-BL-CR} for any $n\geq 0$ we check that
$$
\varsigma_{2n+1}=u_-^{1/2}~\overline{\varsigma}_{2n+1}~u_-^{1/2}\leq \Sigma_{2n+1}\quad \mbox{and}\quad
\Sigma_{2(n+1)}\leq \tau_{2(n+1)}=v^{1/2}_+~\overline{\tau}_{2(n+1)}~v^{1/2}_+,
$$
with the flow of Riccati matrices
 $$
 \overline{\tau}_{2(n+1)}=\mbox{\rm Ricc}_{\overline{\varpi}_0}\left(\overline{\tau}_{2n}\right)  \quad\mbox{and}\quad
 \overline{\varsigma}_{2n+1}=\mbox{\rm Ricc}_{\varpi_1}\left(\overline{\varsigma}_{2n-1}\right)
 $$
 starting at  
 $\tau_0=\tau$ and $\varsigma_{1}=(u_-^{-1}+\cchi^{\prime}\tau\cchi)^{-1}$.
By (\ref{nabla2-WN-ex}) for any $n\geq 1$ we also have
\begin{eqnarray}
  \nabla_2^2W_{2n}
 &\geq& v_+^{-1}+
  \cchi~\varsigma_{2n-1}~\cchi^{\prime}=\tau_{2n}^{-1},\nonumber\\
   \nabla_2^2W^{\flat}_{2n+1}
   &\leq & u_-^{-1}+
    \cchi^{\prime}~ \tau_{2n}~\cchi=\varsigma_{2n+1}^{-1}.
    \label{nabla2-WN-ex-rmk}
\end{eqnarray}
\item Assume (\ref{bary-ref-uv}) holds for some 
$v_->0$ and $u_+>0$.
Following word-for-word the proof of Proposition~\ref{prop-BL-CR} for any $n\geq 0$ we check that
\begin{equation}\label{cov-ev}
\varsigma_{2(n+1)}=v_-^{1/2}~\overline{\varsigma}_{2(n+1)}~v_-^{1/2}\leq \Sigma_{2(n+1)}\quad \mbox{and}\quad
\Sigma_{2n+1}\leq \tau_{2n+1}=u^{1/2}_+~\overline{\tau}_{2n+1}~u^{1/2}_+.
\end{equation}
with the flow of Riccati matrices
 $$
 \overline{\tau}_{2n+1}=\mbox{\rm Ricc}_{\overline{\varpi}_1}\left(\overline{\tau}_{2n-1}\right)  \quad\mbox{and}\quad
 \overline{\varsigma}_{2(n+1)}=\mbox{\rm Ricc}_{\varpi_0}\left(\overline{\varsigma}_{2n}\right)
 $$
 starting at  
 $\tau_1=(\tau+\cchi^{\prime}\tau\cchi)^{-1}$ and $\varsigma_{0}=\tau$. By (\ref{nabla2-WN-ex}) for any $n\geq 1$ we also have
\begin{eqnarray}
  \nabla_2^2W^{\flat}_{2n+1}
 &\geq& u_+^{-1}+
  \cchi^{\prime}~\varsigma_{2n}~\cchi=\tau_{2n+1}^{-1},\nonumber\\
   \nabla_2^2W_{2n}
   &\leq & v_-^{-1}+
    \cchi~ \tau_{2n-1}~\cchi^{\prime}=\varsigma_{2n}^{-1}.
    \label{nabla2-WN-ex-rmk-2}
\end{eqnarray}
\end{itemize}
\end{rmk}

\begin{rmk}\label{rmk-TLS-sd}
By (\ref{HessUpiex}) we have
$$
\nabla^2U^{\pi}_{1}\geq  u_+^{-1}+ \cchi^{\prime}\left(
\tau-
v^{1/2}_+~
~\mbox{\rm Ricc}_{\overline{\varpi}_0}\left(\overline{\tau}_0\right)
~v^{1/2}_+\right)~\cchi.
$$
Note that 
$$
 v^{1/2}_+~
~\mbox{\rm Ricc}_{\overline{\varpi}_0}\left(\overline{\tau}_0\right)
~v^{1/2}_+=\left(v_+^{-1}+\left((\cchi u_-\cchi^{\prime})^{-1}+\tau\right)^{-1}\right)^{-1}\leq v_+.
$$
This yields
$$
\nabla^2U^{\pi}_{1}\geq 
 u_+^{-1}+ \beta^{\prime}\left(\tau^{-1}
-\tau^{-1}v_+\tau^{-1}\right)~\beta.
$$
Similarly, we check that
\begin{eqnarray*}
\nabla^2V^{\pi}_{2}&\geq&v_+^{-1}+\cchi~\left(\left(u_-^{-1}+ \beta^{\prime}~\tau^{-1}~\beta\right)^{-1}-u_+\right)
 \cchi^{\prime}.
\end{eqnarray*}
The last assertion comes from the fact that
$$
u_-^{1/2}~\overline{\varsigma}_{1}~u_-^{1/2}=\left(u_-^{-1}+ \cchi^{\prime}~v^{1/2}_+~\overline{\tau}_{0}~v^{1/2}_+\cchi\right)^{-1}=
\left(u_-^{-1}+ \beta^{\prime}~\tau^{-1}~\beta\right)^{-1}.
$$
More generally,
for any $n\geq 1$ we have the rather crude estimates
 $$
\nabla^2U^{\pi}_{2n+1}\geq u_+^{-1}- \beta^{\prime}\tau^{-1}v_+\tau^{-1}\beta
\quad \mbox{\rm and}\quad
\nabla^2V^{\pi}_{2(n+1)}\geq v_+^{-1}-\tau^{-1}\beta~u_+
 \beta^{\prime}\tau^{-1}.
$$
\end{rmk}

 \begin{defi}
Consider the rescaled matrices
$$
 \begin{array}{rcl}
\overline{\sigma}_{n}&:=&v_+^{-1/2} ~\sigma_{n}~v_+^{-1/2},\\
\overline{s}_{n}&:=&u_-^{-1/2}~
s_{n}~u_-^{-1/2}
 \end{array} \qquad \&\qquad \begin{array}{rcl}
\overline{v}_{n}&:=&v_-^{-1/2}~v_{n+1}~v_-^{-1/2},\\
\overline{\omega}_{n}&:=&
u_+^{-1/2}~{\omega}_{n}~u_+^{-1/2}.
  \end{array}
  $$

\end{defi}

In this notation (\ref{bbPhi}) takes the following form

\begin{equation}\label{ovvt}
 \begin{array}{rcl}
\overline{\sigma}_{n+1}&=&\Psi_{\gamma_1^{\prime}}\left(\overline{s}_n\right)\\
\overline{s}_{n+1}&=&\Psi_{\gamma_1}\left(\overline{\sigma}_n\right)
 \end{array} \qquad \&\qquad \begin{array}{rcl}
\overline{v}_{n+1}&=&\Psi_{\gamma_0}\left(\overline{\omega}_{n}\right),\\ 
\overline{\omega}_{n+1}&=&\Psi_{\gamma_0^{\prime}}\left(\overline{v}_n\right).\end{array}
\end{equation}
This implies that
$$
 \begin{array}{rcl}
\overline{\sigma}_{n+2}&=&
\mbox{\rm Ricc}_{\overline{\varpi}_0}\left(\overline{\sigma}_{n}\right),
\\
\overline{s}_{n+2}&=&\mbox{\rm Ricc}_{{\varpi}_1}(\overline{s}_n)
 \end{array} \qquad \&\qquad \begin{array}{rcl}
\overline{v}_{n+2}&=&\mbox{\rm Ricc}_{{\varpi}_0}(\overline{v}_n),
\\
\overline{\omega}_{n+2}&=&\mbox{\rm Ricc}_{\overline{\varpi}_1}\left(\overline{\omega}_{n}\right).
\end{array}
$$

\begin{rmk}
\begin{itemize}
\item When $u_-=0$  we use the convention $\Phi^{\prime}_{u_-}(\sigma_n)=0=s_{n+1}$. In this case we have
$$
\sigma_{n+2}=\Phi_{v_+}(0)=v_+.$$
By (\ref{convention-u})  for any $s\geq 0$ we have 
$$
\overline{\sigma}_{n+2}=\Psi_{\gamma_1^{\prime}}(0)=I=\mbox{\rm Ricc}_{\overline{\varpi}_0}(s)=r_{\overline{\varpi}_0}\quad
\mbox{and}\quad\mbox{\rm Ricc}_{{\varpi}_1}(s)=I=r_{\varpi_1}=\overline{s}_{n+2}.
$$
\item When $v_-=0$ we also use the convention $\Phi_{v_-}(\omega_n)=0=v_{n+1}$. In this case we have
$$
{\omega}_{n+2}=
\Phi^{\prime}_{u_+}(0)=u_+.
$$
By (\ref{convention-v})  for any $s\geq 0$ we have 
$$
\overline{\omega}_{n+2}=\Psi_{\gamma_0^{\prime}}(0)=I=\mbox{\rm Ricc}_{\overline{\varpi}_1}(s)=r_{\overline{\varpi}_1}
\quad
\mbox{and}\quad\mbox{\rm Ricc}_{{\varpi}_0}(s)=I=r_{\varpi_0}=\overline{v}_{n+2}.
$$
\end{itemize}\end{rmk}

Using (\ref{bcov}) we conclude that
\begin{equation}\label{infty-cov}
\begin{array}{rcccl}
v_-^{1/2} ~r_{\varpi_0}~v_-^{1/2}&\leq& \Sigma_{\mu,\eta}&\leq &v_+^{1/2} ~r_{\overline{\varpi}_0}~v_+^{1/2},
\\
u_-^{1/2}~r_{\varpi_1}~u_-^{1/2}&\leq&  \Sigma^{\flat}_{\eta,\mu}&\leq &u_+^{1/2}~r_{\overline{\varpi}_1}~u_+^{1/2}.
\end{array}
\end{equation}

 \begin{examp}\label{l-exam-Vlspi}
Consider the linear-Gaussian  potential $W$ of the form (\ref{def-W}).
Using (\ref{ref-nablaUVng-bg}) and (\ref{nabla2-WWN-rr-2}) and following
word-for-word the arguments given in Example~\ref{l-exam-DeltaUV} we check that
$$
\begin{array}{l}
\nabla \UU(x_1)-\nabla \UU(x_2)+\cchi^{\prime}\beta(x_1-x_2)\\
\\
=\nabla U(x_1)-\nabla U(x_2)+\cchi^{\prime}~(K_{\WW_{\mu,\eta}}(I)(x_1)-
K_{\WW_{\mu,\eta}}(I)(x_2))\\
\\
=   \nabla_2W^{\flat}_{\eta,\mu}(y,x_1)-  \nabla_2W^{\flat}_{\eta,\mu}(y,x_2).
\end{array}$$
Using (\ref{finex}) we also check that
\begin{equation}\label{nbUU}
 \begin{array}{l}
 (x_1-x_2)^{\prime}  \left(\nabla_2W^{\flat}_{\eta,\mu}(y,x_1)-    \nabla_2W^{\flat}_{\eta,\mu}(y,x_2)\right)\\
 \\
\displaystyle\geq 
 (x_1-x_2)^{\prime} \left(\nabla U(x_1)-\nabla U(x_2)\right) + (x_1-x_2)^{\prime}~\cchi^{\prime}~\Sigma_{\mu,\eta}(x_1,x_2)~\cchi~(x_1-x_2),
\end{array} 
\end{equation}
with
$$
\Sigma_{\mu,\eta}(x_1,x_2):=
 \int_0^1~\Sigma_{\mu,\eta}(x_2+s(x_1-x_2))~ds.
$$
Using (\ref{infty-cov}) we conclude that
$$
 \begin{array}{l}
 (x_1-x_2)^{\prime}  \left(\nabla_2W^{\flat}_{\eta,\mu}(y,x_1)-    \nabla_2W^{\flat}_{\eta,\mu}(y,x_2)\right)\\
 \\
\displaystyle\geq 
 (x_1-x_2)^{\prime} \left(\nabla U(x_1)-\nabla U(x_2)\right) + (x_1-x_2)^{\prime}~\cchi^{\prime}~v_-^{1/2} ~r_{\varpi_0}~v_-^{1/2}~\cchi~(x_1-x_2).
\end{array} $$
\end{examp}

Now we come to the proof of Theorem~\ref{last-est-th}.\\

\noindent{\bf Proof of Theorem~\ref{last-est-th}:}\label{last-est-th-proof}
For a linear-Gaussian  potential $W$ of the form (\ref{def-W}), 
 using (\ref{finex}) and (\ref{bcov}) we check the estimates
   \begin{eqnarray*}
\cchi^{\prime}~ v_-^{1/2} ~r_{\varpi_0}~v_-^{1/2}~\cchi\leq \nabla K_{\WW_{\mu,\eta}}(I)~\cchi&=&\cchi^{\prime}~\Sigma_{\mu,\eta}~\cchi\leq
\cchi^{\prime}~v_+^{1/2} ~r_{\overline{\varpi}_0}~v_+^{1/2}~\cchi,\\
\cchi~u_-^{1/2}~r_{\varpi_1}~u_-^{1/2}~\cchi^{\prime}\leq \nabla K_{\WW^{\flat}_{\eta,\mu}}(I)~\cchi^{\prime}&=&\cchi~\Sigma^{\flat}_{\eta,\mu}~\cchi^{\prime}\leq
\cchi~u_+^{1/2}~r_{\overline{\varpi}_1}~u_+^{1/2}~\cchi^{\prime}. \end{eqnarray*}
This ends the proof of (\ref{last-est}). The proof of the theorem is completed.\cqfd

\subsection*{Proof of Lemma~\ref{lem-HW}}\label{lem-HW-proof}
By (\ref{ricc-maps-incr}) for any $n\geq p\geq 1$ we have
  $$
\overline{\tau}_{2(n+1)}\leq \mbox{\rm Ricc}_{\overline{\varpi}_0}^{n}(I)\leq \mbox{\rm Ricc}_{\overline{\varpi}_0}^{p}(I)\leq I.
  $$
  Using (\ref{ovvt}) and (\ref{ricc-maps-incr}) for any $n\geq p\geq 1$ we also have
\begin{eqnarray*}
\overline{\tau}_{2n+1}&=&\Psi_{\overline{\gamma}_1}(\overline{\varsigma}_{2n})\leq 
\Psi_{\overline{\gamma}_1}\left( \mbox{\rm Ricc}_{\varpi_0}^n(0)\right)\leq  
\Psi_{\overline{\gamma}_1}\left( \mbox{\rm Ricc}_{\varpi_0}^p(0)\right).\end{eqnarray*}

Using
 Lemma~\ref{varphitoRicc} we check that
\begin{eqnarray*} 
 \Psi_{\overline{\gamma}_1}\circ  \mbox{\rm Ricc}_{\varpi_0}^p&=&
  \Psi_{\gamma_0^{\prime}}\circ  \mbox{\rm Ricc}_{\psi(\gamma_0)}^p=
   \mbox{\rm Ricc}_{\psi(\gamma_0^{\prime})}^p\circ   \Psi_{\gamma_0^{\prime}}=   \mbox{\rm Ricc}_{\psi(\overline{\gamma}_1)}^p\circ   \Psi_{\overline{\gamma}_1}.
 \end{eqnarray*}
This yields the formula
$$
 \Psi_{\overline{\gamma}_1}\left( \mbox{\rm Ricc}_{\varpi_0}^p(0)\right)=
 \mbox{\rm Ricc}_{\overline{\varpi}_1}^p\left(   \Psi_{\overline{\gamma}_1}(0)\right)= \mbox{\rm Ricc}_{\overline{\varpi}_1}^p\left(  I\right).
$$
The proof of the lemma is completed.
\cqfd

\subsection*{Linear-Gaussian models}\label{lin-gauss-sec}

This section provides a brief overview of the
closed form expressions of the Sch\" odinger bridges (\ref{bridge-form-intro}) and the distributions (\ref{transport-g}) for general classes of linear-Gaussian models presented in the article~\cite{adm-24}.

In the context of Gaussian distributions (\ref{def-U-V}), we recall from (\ref{def-over-varpi-g}) that
 $$
 ({\varpi}_0,
{\varpi}_1)= (\overline{\varpi}_0,
\overline{\varpi}_1)=\left(v^{-1/2}~\left(\cchi~u~\cchi^{\prime}\right)^{-1}~v^{-1/2}, u^{-1/2}~(\cchi^{\prime}~v~\cchi)^{-1}~u^{-1/2}
\right).
 $$
 In this situation, $\Ka_{2n}$ and $\Ka_{2n+1} $ are linear Gaussian transitions. Their covariance matrices are constant functions given by
$$
\Sigma_{2n}=\tau_{2n}=v^{1/2}~\overline{\tau}_{2n}~v^{1/2}\quad\mbox{\rm and}\quad
\Sigma_{2n+1}=\tau_{2n+1}=u^{1/2}~\overline{\tau}_{2n+1}~u^{1/2},
$$
with the solutions of the Riccati equations
$$
 \begin{array}{rclcccrcl}
\overline{\tau}_{2n}&=&\mbox{\rm Ricc}_{\overline{\varpi}_0}\left(\overline{\tau}_{2(n-1)}\right)  &&~\mbox{and}~\quad&
 \overline{\tau}_{2n+1} &=&\mbox{\rm Ricc}_{\overline{\varpi}_1}\left(\overline{\tau}_{2n-1}\right),
 \end{array}
$$
and the initial conditions
$$
\begin{array}{rcl}
\overline{\tau}_{0}&=&v^{-1/2}~\tau~v^{-1/2}
  \end{array}
\quad\mbox{\rm and}\quad
  \begin{array}{rcl}
\overline{\tau}_{1} &=&u^{-1/2}~\tau_1~u^{-1/2}
 \end{array}
 \quad\mbox{\rm with}\quad
\tau_1^{-1}=  u^{-1}+~\cchi^{\prime}~\tau~~\cchi.
  $$
  In the linear Gaussian situation, $\pi_{n}$ are also Gaussian distributions with mean $m_n$ and covariance $\sigma^{\pi}_n>0$ for any $n\geq 0$. 
  Using (\ref{ref-WW-Gauss-eom}) and (\ref{Sigma-Gauss}) we also have
\begin{eqnarray}
\ZZ_{\mu,\pi_{2n}}(x)&=&m_{2n}+\tau_{2n}~\cchi~ (x-m)+\tau_{2n}^{1/2}~G,
\nonumber\\
\ZZ_{\eta,\pi_{2n+1}}^{\flat}(x)&=&m_{2n+1}+\tau_{2n+1}~\cchi^{\prime}~ (x-\overline{m})+\tau_{2n+1}^{1/2}~G,
\end{eqnarray}
with the covariance matrices
\begin{equation}\label{ref-fpp}
\sigma_{\mu,\pi_{2n}}=\tau_{2n}
\quad\mbox{\rm and}\quad
\sigma_{\eta,\pi_{2n+1}}^{\flat}=\tau_{2n+1}.
\end{equation}
By (\ref{bridge-form-intro}) and (\ref{sbr-ref})
\begin{eqnarray*}
\Pa_{2n}&=&P_{\mu,\pi_{2n}}=\mu\times \Ka_{2n}=\mu~\times~K_{\WW_{\mu,\pi_{2n}}},
\\
\Pa_{2n+1}^{\flat}&=&P_{\eta,\pi_{2n+1}}=\eta\times\Ka_{2n+1} =\eta~\times~K_{\WW^{\flat}_{\eta,\pi_{2n+1}}}.
\end{eqnarray*}
The parameters $m_n$ satisfy the recursion
\begin{eqnarray*}
m_{2n+1}&=&m+
\tau_{2n+1}~\cchi^{\prime}~ (\overline{m}-m_{2n}),\\
m_{2(n+1)}&=&\overline{m}+\tau_{2(n+1)}~\cchi~(m-m_{2n+1}).
\end{eqnarray*}
with the initial condition
$
m_0=\alpha_0+\beta_0 m$. In addition, we have
$$
\sigma_{2n}^{\pi}:=
\tau_{2n}~\cchi~u~\cchi^{\prime}~\tau_{2n}+\tau_{2n}
\quad \mbox{and}\quad\sigma_{2n+1}^{\pi}:=\tau_{2n+1}~\cchi^{\prime}~v~\cchi~\tau_{2n+1}+\tau_{2n+1}.
$$
Formulae (\ref{ref-fpp}) can also be checked using the matrix inversion lemma. For instance, set
$$
\varpi^{\pi}_{2n}:=(\sigma_{2n}^{\pi})^{-1/2}~(\cchi ~u~\cchi^{\prime})^{-1}~(\sigma_{2n}^{\pi})^{-1/2}.
$$
To check that
$$
r_{\varpi^{\pi}_{2n}}:=
(\sigma_{2n}^{\pi})^{-1/2}\tau_{2n}~(\sigma_{2n}^{\pi})^{-1/2}=\mbox{\rm Ricc}_{\varpi^{\pi}_{2n}}\left(r_{\varpi^{\pi}_{2n}}\right),
$$
observe that
$$
(\sigma_{2n}^{\pi})^{-1}=\tau_{2n}^{-1}-\left(\tau_{2n}+
(\cchi u\cchi^{\prime})^{-1}\right)^{-1}=
\tau_{2n}^{-1}-\tau_{2n}^{-1}\left(\tau_{2n}^{-1}+\tau_{2n}
(\cchi u\cchi^{\prime})^{-1}\tau_{2n}\right)^{-1}\tau_{2n}^{-1}.
$$
\section{Some technical proofs}\label{sec-tech-proofs}
\subsection*{Proof of (\ref{prop-Hent})}\label{prop-Hent-proof}
The proof of the l.h.s. formula (a.k.a. the data processing inequality) is well known, see for instance~\cite{dm-03}. For instance in terms of the dual  transition
(\ref{dual-transition-r}) we have that
\begin{eqnarray*}
\Ha(\nu~|~\mu)&=&\Ha(\nu\times\Ka~|~\mu\times\Ka )\\
&=&\Ha((\nu\Ka)\times\Ka^{\ast}_{\nu}~|~(\mu\Ka)\times \Ka^{\ast}_{\mu} )\\
&=&
\Ha(\nu\Ka~|~\mu\Ka )+\int(\nu K)(dx)~\Ha(\delta_x\Ka^{\ast}_{\nu}~|~\delta_x\Ka^{\ast}_{\mu})\geq \Ha(\nu\Ka~|~\mu\Ka ).
\end{eqnarray*}
The r.h.s. estimate comes from  the convexity of the function $\varphi(u)=u\log u$. Indeed we have 
\begin{eqnarray*}
 \varphi\left(\frac{d\mu\Ka}{d\nu}(y)\right)
 &=&\varphi\left(\int~\mu(dx)~\frac{d\delta_{x}\Ka}{d\nu}(y)\right)\leq \int \mu(dx)~\varphi\left(\frac{d\delta_x\Ka}{d\nu}(y)\right).
\end{eqnarray*}
This implies that
\begin{eqnarray*}
\Ha(\mu\Ka~|~\nu)&=&\int~\nu(dy) ~\varphi\left(\frac{d\mu\Ka}{d\nu}(y)\right)\\
&\leq &\int~~\mu(dx) \int~\nu(dy) ~\varphi\left(\frac{d\delta_x\Ka}{d\nu}(y)\right) =
\int~\mu(dx)~\Ha(\delta_{x}\Ka~|~\nu).
\end{eqnarray*}
This ends the proof of (\ref{prop-Hent}).\cqfd

\subsection*{Proof of Lemma~\ref{varphitoRicc}}\label{varphitoRicc-proof}

Applying the matrix inversion lemma we have

$$
\Psi_{\gamma^{\prime}}(v):=(I+\gamma^{\prime} v \gamma)^{-1}=I-
\gamma^{\prime} (v^{-1}+ \gamma \gamma^{\prime} )^{-1}\gamma.
$$
This yields the formula
\begin{eqnarray*}
(\Psi_{\gamma}(\Psi_{\gamma^{\prime}}(v)))^{-1}&=&I+\gamma~\Psi_{\gamma^{\prime}}(v) ~\gamma^{\prime}\\
&=&I+\gamma \gamma^{\prime}-
\gamma\gamma^{\prime} (v^{-1}+ \gamma \gamma^{\prime} )^{-1}\gamma
\gamma^{\prime}=I+\left((\gamma\gamma^{\prime})^{-1}+v\right)^{-1}.
\end{eqnarray*}
This ends the proof of (\ref{ff-1}).  We check the l.h.s. of (\ref{ff-2}) by induction w.r.t. the parameter $n\geq 0$. The formula (\ref{ff-2}) is clearly met for $n=0$. 
By (\ref{ff-1}) also note that
$$
\Psi_{\gamma^{\prime}}\circ 
\mbox{\rm Ricc}_{\psi(\gamma)}=
(\Psi_{\gamma^{\prime}}\circ \Psi_{\gamma})\circ\Psi_{\gamma^{\prime}}
=\mbox{\rm Ricc}_{\psi(\gamma^{\prime})}\circ\Psi_{\gamma^{\prime}}.
$$
Assume l.h.s. of  (\ref{ff-2})  is met at rank $n$. In this case, we have
\begin{eqnarray*}
\Psi_{\gamma^{\prime}}\circ 
\mbox{\rm Ricc}^{n+1}_{\psi(\gamma)}&=&\left(\Psi_{\gamma^{\prime}}\circ 
\mbox{\rm Ricc}_{\psi(\gamma)}\right)\circ\mbox{\rm Ricc}^{n}_{\psi(\gamma)}\\
&=&\mbox{\rm Ricc}_{\psi(\gamma^{\prime})}\circ\left(\Psi_{\gamma^{\prime}}\circ\mbox{\rm Ricc}^{n}_{\psi(\gamma)}\right)=\mbox{\rm Ricc}_{\psi(\gamma^{\prime})}\circ\left(
\mbox{\rm Ricc}_{\psi(\gamma^{\prime})}^n
\circ \Psi_{\gamma^{\prime}}\right).
\end{eqnarray*}
This shows that l.h.s. of  (\ref{ff-2})  is met at rank $(n+1)$. The inductive proof is now completed. Finally observe that
$$
\Psi_{\gamma^{\prime}}\left(
\mbox{\rm Ricc}^n_{\psi(\gamma)}(r_{\psi(\gamma)})\right)=\Psi_{\gamma^{\prime}}\left(
r_{\psi(\gamma)}\right)=\mbox{\rm Ricc}_{\psi(\gamma^{\prime})}^n
\left( \Psi_{\gamma^{\prime}}(r_{\psi(\gamma)})\right).
$$
By uniqueness of the fixed point $r_{\psi(\gamma^{\prime})}$ of $\mbox{\rm Ricc}_{\psi(\gamma^{\prime})}$ we conclude that
$\Psi_{\gamma^{\prime}}(r_{\psi(\gamma)})=r_{\psi(\gamma^{\prime})}$.

This ends the proof of the lemma.
\cqfd

\subsection*{Proof of Theorem~\ref{theo-imp}}\label{theo-imp-proof}

Using (\ref{e12}) and (\ref{bridge-form-v2-inside}) we check that
\begin{eqnarray*}
 \Ha(P_{\mu,\eta}~|~\Pa_{2n})&=&   \Ha(P_{\mu,\eta}|~\Pa_{2n-1})-\Ha(\mu~|~\pi_{2n-1})\\
 &\leq &  \Ha(P_{\mu,\eta}|~\Pa_{2n-1})- \Ha(\pi_{2n}~|~\eta)\\
 &\leq & \Ha(P_{\mu,\eta}|~\Pa_{2n-1})- \epsilon_{2n+1}^{-1}~ \Ha(P_{\mu,\eta}~|~\Pa_{2n+1}).
\end{eqnarray*}
We set 
$$
I_{n}:= \Ha(P_{\mu,\eta}~|~\Pa_{n}).
$$
This implies that for any $n\geq 1$ we have
$$
(1+\epsilon_{\star}^{-1})~I_{2n+1}~\leq ~I_{2n}+\epsilon_{\star}^{-1}~I_{2n+1}
\leq I_{2n-1}.
$$
In the same vein, using (\ref{bridge-form-2-V2-inside}) we check that
\begin{eqnarray*}
   \Ha(P_{\mu,\eta}~|~\Pa_{2n+1})&=& 
    \Ha(P_{\mu,\eta}~|~\Pa_{2n})-   \Ha(\eta~|~\pi_{2n})\\
    &\leq &   \Ha(P_{\mu,\eta}~|~\Pa_{2n})-     \Ha(\pi_{2n+1}~|~\mu)\\
    &\leq & \Ha(P_{\mu,\eta}~|~\Pa_{2n})-\epsilon_{2(n+1)}^{-1}    \Ha(P_{\mu,\eta}~|~\Pa_{2(n+1)}).
\end{eqnarray*}
which implies that for any $n\geq 0$ we have
$$
(1+\epsilon_{\star}^{-1}) ~   I_{2(n+1)}\leq I_{2n+1}+\epsilon_{\star}^{-1}    I_{2(n+1)}
\leq  I_{2n}.$$
In summary we have proved that for any $n\geq 0$ the second order linear homogeneous recurrence relation
$$
 I_{n}\geq I_{n+1}+  \epsilon_{\star}^{-1} I_{n+2} \geq (1+\epsilon_{\star}^{-1})~I_{n+2}.
$$

For any $n\geq 0$ and $p\geq 1$ we have
$$
 I_{n}\geq a_p~I_{n+p}+  b_p~ I_{n+(p+1)}\geq (a_p+b_p)~I_{n+(p+1)},
$$
with the parameters $z_p:=\left(
\begin{array}{l}
a_{p}\\
b_{p}
\end{array}
\right)$
defined by $z_1:=\left(
\begin{array}{l}
1\\
\epsilon_{\star}^{-1}
\end{array}
\right)
$ and the recursive equation
$$
z_{p+1}=M
z_p\quad\mbox{\rm with}\quad M:=\left(
\begin{array}{cc}
1&1\\
 \epsilon_{\star}^{-1} &0
\end{array}
\right).
$$
The eigenvalues of $M$ are defined by
$$
\ell_1:=\frac{1}{2}~(1-\sqrt{1+4  \epsilon_{\star}^{-1} })<0<1< 
\ell_2:=\frac{1}{2}~(1+\sqrt{1+4  \epsilon_{\star}^{-1} }).
$$
The corresponding eigenvectors are given by
$$
v_1=\left(
\begin{array}{l}
\ell_1\\
\epsilon_{\star}^{-1}
\end{array}
\right)\quad\mbox{\rm and}\quad
v_2=\left(
\begin{array}{l}
\ell_2\\
\epsilon_{\star}^{-1}
\end{array}
\right)\Longrightarrow
\left(
\begin{array}{cc}
\ell_1&\ell_2\\
 \epsilon_{\star}^{-1} & \epsilon_{\star}^{-1} 
\end{array}
\right)^{-1}=\frac{1}{\ell_2-\ell_1}
\left(
\begin{array}{cc}
-1&\ell_2~ \epsilon_{\star}\\
1 &-\ell_1~ \epsilon_{\star} 
\end{array}
\right).
$$ 
This yields the diagonalisation formulae 
$$
A=\frac{1}{\ell_2-\ell_1}~\left(
\begin{array}{cc}
\ell_1&\ell_2\\
 \epsilon_{\star}^{-1} & \epsilon_{\star}^{-1} 
\end{array}
\right)
\left(
\begin{array}{cc}
\ell_1&0\\
0 &\ell_2
\end{array}
\right)
\left(
\begin{array}{cc}
-1&\ell_2~ \epsilon_{\star}\\
1 &-\ell_1~ \epsilon_{\star} 
\end{array}
\right).
$$
We conclude that
$$
z_{p+1}=A^pz_1=\frac{1}{\ell_2-\ell_1}~\left(
\begin{array}{cc}
\ell_1&\ell_2\\
 \epsilon_{\star}^{-1} & \epsilon_{\star}^{-1} 
\end{array}
\right)
\left(
\begin{array}{cc}
\ell_1^p&0\\
0 &\ell_2^p
\end{array}
\right)
\left(
\begin{array}{cc}
-1&\ell_2~ \epsilon_{\star}\\
1 &-\ell_1~ \epsilon_{\star} 
\end{array}
\right)\left(
\begin{array}{l}
1\\
\epsilon_{\star}^{-1}
\end{array}
\right).
$$
This yields the formulae
$$
a_{p+1}=\ell_1^{p+1}~\frac{\ell_2-1}{\ell_2-\ell_1}+
\ell_2^{p+1}~\frac{1-\ell_1}{\ell_2-\ell_1}
\quad\mbox{\rm 
and }
\quad
b_{p+1}=\epsilon_{\star}^{-1}~a_p.
$$
Using the formulae
$$
\begin{array}{l}
\displaystyle\ell_2-\ell_1=\sqrt{1+4  \epsilon_{\star}^{-1} }\quad \mbox{\rm and}\quad
\ell_1+\ell_2=1,
\end{array}
$$
we check that
\begin{eqnarray*}
a_{p}&=&\frac{\ell_2^{p+1}}{\sqrt{1+4  \epsilon_{\star}^{-1} }}\left(
1-\left(\frac{1-\ell_2}{\ell_2}\right)^{p+1}\right)
\\
&\geq& 
\frac{\ell_2^{p+1}}{\sqrt{1+4  \epsilon_{\star}^{-1} }}\left(
1-\left(\frac{\ell_2-1}{\ell_2}\right)^{p+1}\right)
\geq \frac{\ell_2^{p+1}}{\sqrt{1+4  \epsilon_{\star}^{-1} }}\left(
1-\frac{\ell_2-1}{\ell_2}\right)=
 \frac{\ell_2^{p}}{\sqrt{1+4  \epsilon_{\star}^{-1} }}.
\end{eqnarray*}
This yields the estimate
$$
a_{p+1}+b_{p+1}\geq \left(
 \frac{\ell_2}{\sqrt{1+4  \epsilon_{\star}^{-1} }}+
 \epsilon_{\star}^{-1} \frac{1}{\sqrt{1+4  \epsilon_{\star}^{-1} }}\right)~\ell_2^{p}\geq 
 \ell_2^{p}.
$$
The last assertion comes from the fact that
$$
\frac{\ell_2}{\sqrt{1+4  \epsilon_{\star}^{-1} }}+ \epsilon_{\star}^{-1} \frac{1}{\sqrt{1+4  \epsilon_{\star}^{-1} }}=\frac{1}{2}~\left(1+\frac{1+2\epsilon_{\star}^{-1}}{\sqrt{1+4  \epsilon_{\star}^{-1} }}
\right)\geq 1.
$$
Finally switching the roles of the parameters $(p,n)$
observe that for any $p\geq 0$ and $n\geq 1$ we have
$$
I_{(n+1)+p}\leq \frac{1}{a_n+b_n}
 I_{p}\leq \ell_2^{-((n+1)-2)}~I_p.
$$
Observe that
$$
\frac{1+\sqrt{1+4  \epsilon_{\star}^{-1} }}{2}=1+\phi(\epsilon_{\star})\quad
\mbox{\rm with}\quad
\phi(\epsilon_{\star}):=\frac{\epsilon_{\star}^{-1} }{\sqrt{1/4+ \epsilon_{\star}^{-1} }+1/2}.
$$

This ends the proof of the theorem.\cqfd

\subsection*{Proof of Theorem~\ref{theo-impp}}\label{theo-impp-proof}

Conditions $A_2(\Ka_{n})$ and $A_2(\Sa)$ are met with $\kappa$ as in (\ref{def-kappa})
and the parameters $\rho_{n}$ given for any $n\geq p\geq  1$ by 
  \begin{eqnarray*}
  (\rho_0,\rho_1,\rho_2)&=&(\Vert\tau\Vert_2,\Vert\tau_1\Vert_2,\Vert\tau_2\Vert_2),\\
 \rho_{2n+1} &=& \rho_{2p+1},\\
 &=&\Vert u_+\Vert_2/\xi_{2p+1}\quad \mbox{\rm and}\quad
 \rho_{2(n+1)} = \rho_{2(p+1)}=\Vert v_+\Vert_2/\xi_{2p}.
  \end{eqnarray*}
In this context, the parameters  $\epsilon_{2n+1}$ defined in (\ref{def-varepsi-SW}) takes the form
  \begin{eqnarray*}
(\epsilon_0, \epsilon_1,\epsilon_2 )&=&(\varepsilon_{\kappa}(\rho_{0},\Vert u_+\Vert_2),\varepsilon_{\kappa}(\rho_{1},\Vert v_+\Vert_2),\varepsilon_{\kappa}(\rho_{2},\Vert u_+\Vert_2)),\\
 \epsilon_{2n+1}&=& \epsilon_{2p+1}\\
 &=&\varepsilon_{\kappa}(\Vert u_+\Vert_2,\Vert v_+\Vert_2)/\xi_{2p+1}
\quad \mbox{\rm and}\quad
 \epsilon_{2(n+1)}= \epsilon_{2(p+1)}=\varepsilon_{\kappa}(\Vert u_+\Vert_2,\Vert v_+\Vert_2)/\xi_{2p}.
  \end{eqnarray*}

  Using (\ref{f2-c}) we conclude that
$$
 \Ha(P_{\mu,\eta}~|~\Pa_{2(n+1)})\leq   \Ha(P_{\mu,\eta}~|~\Pa_{2n+1})\leq  \left(1+  \varepsilon^{-1}_{2p+1}\right)^{-(n-p)} \underbrace{ \Ha(P_{\mu,\eta}~|~\Pa_{2p-1})}_{ \leq \Ha(P_{\mu,\eta}~|~\Pa_{2(p-1)})}.
$$
Using (\ref{f1-c}) we conclude that
$$
   \Ha(P_{\mu,\eta}~|~\Pa_{2(n+1)+1})\leq
 \Ha(P_{\mu,\eta}~|~\Pa_{2(n+1)})\leq \left(1+  \varepsilon^{-1}_{2(p+1)}
\right)^{-(n-p)}~\underbrace{ \Ha(P_{\mu,\eta}~|~\Pa_{2p}) }_{\leq \Ha(P_{\mu,\eta}~|~\Pa_{2p-1}) }.
$$
This yields for any $n\geq p\geq 1$ the estimate (\ref{ineq-theo-impp}).
This ends the proof of the theorem.
\cqfd

\subsection*{Proof of Corollary~\ref{cor-impp}}\label{cor-impp-proof}
  Using (\ref{cv-ricc-intro}) there exists some constants $(c_{\overline{\varpi}_0},c_{\overline{\varpi}_1})$ such that for any $p\geq 1$ and $i\in\{0,1\}$ we have
$$
  \iota_{i}(\overline{\varpi}_i)~\left(1+
c_{\overline{\varpi}_i}~\delta_{\overline{\varpi}_i}^p\right)^{-1}\leq \xi_{2p+i}\leq \iota_{i}(\overline{\varpi}_i),
$$
with the exponential decay parameter $\delta_{\varpi}$ defined in (\ref{cv-ricc-intro}).
This implies that
$$
  \left(\iota_{0}(\overline{\varpi}_0)\vee \iota_{1}(\overline{\varpi}_1)\right)~\left(1+
c_{\overline{\varpi}_0,\overline{\varpi}_1}~\delta_{\overline{\varpi}_0,\overline{\varpi}_1}^p\right)^{-1}\leq \xi_{2p}\vee \xi_{2p+1}.
$$
This ends the proof of the theorem.
\cqfd

\section*{Fundings and Competing interests}

I confirm that no funding was received. 
I also confirm that I have no competing interests, no relevant financial or non-financial competing interests to report.

\end{document}